\documentclass[11pt, a4paper]{amsart}
\usepackage{graphicx} 

\usepackage[T1]{fontenc}
\usepackage{microtype}

\usepackage{amsmath}								
\usepackage{amssymb}
\usepackage{amsthm}
\usepackage{amscd}
\usepackage{amsfonts}
\usepackage{bbm}
\usepackage{stmaryrd}
\usepackage{appendix}

\usepackage{extarrows}

\usepackage[bookmarks=true, colorlinks=true, linkcolor=blue!50!black,
citecolor=orange!50!black, urlcolor=orange!50!black, pdfencoding=unicode]{hyperref}
\usepackage{color}
\usepackage{etoolbox}
\usepackage{zref-clever} 
\zcsetup{cap} 
\usepackage{comment}
\usepackage{mathrsfs}

\usepackage{tikz}									
\usetikzlibrary{matrix}
\usetikzlibrary{patterns}
\usetikzlibrary{matrix}
\usetikzlibrary{positioning}
\usetikzlibrary{decorations.pathmorphing}
\usetikzlibrary{cd}
\usetikzlibrary{trees}
\usetikzlibrary{calc}
\usetikzlibrary{intersections}

\usepackage{fullpage}
\usepackage{enumerate}

\makeatletter
\def\l@section{\@tocline{1}{0pt}{1pc}{}{}}
\def\l@subsection{\@tocline{2}{0pt}{3pc}{5pc}{}}
\def\l@subsubsection{\@tocline{3}{0pt}{5pc}{7pc}{}}
\makeatother

\zcRefTypeSetup{conjecture}{Name-sg=Conjecture, name-sg=conjecture,
  Name-pl=Conjectures, name-pl=conjectures}
\zcRefTypeSetup{question}{Name-sg=Question, name-sg=question,
  Name-pl=Questions, name-pl=questions}
\zcRefTypeSetup{construction}{Name-sg=Construction, name-sg=construction,
  Name-pl=Constructions, name-pl=constructions}

\ExplSyntaxOn
\newcommand{\zcrefglobalstringname}[2]{
  \__zrefclever_opt_varname_lang_type:nnnn{\languagename}{#1}{#2}{tl}
}
\newcommand{\zcreflocalstringname}[2]{
  \__zrefclever_opt_varname_type:een{#1}{#2}{tl}
}
\newcommand{\zcrefgetstring}[2]{
  \__zrefclever_provide_langfile:e { \languagename }
  \__zrefclever_process_language_settings:
  \ifcsvoid{\zcreflocalstringname{#1}{#2}}
    {\csname \zcrefglobalstringname{#1}{#2} \endcsname}
    {\csname \zcreflocalstringname{#1}{#2} \endcsname}
}
\ExplSyntaxOff
\newcommand{\zcreftypename}[1]{\zcrefgetstring{#1}{Name-sg}}
\newtheorem{sharedtheoremcounter}{}[section]
\undef\sharedtheoremcounter 

\NewDocumentCommand{\newzctheorem}{mO{#1}O{\zcreftypename{#2}}}{
  \newtheorem{#1}[sharedtheoremcounter]{#3}
    \AddToHook{env/#1/begin}{%
      \zcsetup{countertype={sharedtheoremcounter=#2}}}
}

\newtheorem*{theorem*}{Theorem}

\newtheorem{maintheorem}{Theorem}[section]

\zcsetup{countertype={maintheorem=theorem}}
\newtheorem{guidingtheorem}{Theorem}
\zcsetup{countertype={guidingtheorem=theorem}}
\AddToHook{env/maincorollary/begin}{%
  \zcsetup{countertype={maintheorem=corollary}}}

\newzctheorem{theorem}
\newzctheorem{lemma}
\newzctheorem{proposition}
\newzctheorem{corollary}
\newzctheorem{conjecture}
\newzctheorem{question}

\theoremstyle{definition}
\newzctheorem{definition}
\newzctheorem{example}
\newzctheorem{construction}

\newzctheorem{remark}
\newzctheorem{remarks}[remark][Remarks]

\newtheoremstyle{myitemstyle}						
	{}			
	{}			
	{}			
	{}			
	{}			
	{.}			
	{ }			
	{}			
\theoremstyle{myitemstyle}
\newtheorem{myitemthm}{}

\newcommand{\R}{\mathbb{R}}

\newcommand{\Z}{\mathbb{Z}}

\newcommand{\Q}{\mathbb{Q}}
\newcommand{\C}{\mathbb{C}}

\renewcommand{\P}{\mathbb{P}}
\newcommand{\PP}{\mathbb{P}}

\newcommand{\RR}{\mathbb{R}}
\newcommand{\A}{\mathbb{A}}
\newcommand{\G}{\mathbb{G}}

\newcommand{\F}{\mathbb{F}}

\newcommand{\Gm}{\mathbb{G}_m}
\newcommand{\Gr}{\mathrm{Gr}}
\newcommand{\Groo}{\mathrm{Gr}^{\circ}}

\newcommand{\Trop}{\operatorname{Trop}}

\newcommand{\Dr}{\operatorname{Dr}}
\newcommand{\inw}{\operatorname{in}_w}

\newcommand{\sign}{\operatorname{sign}}
\DeclareMathOperator{\cromb}{cr}
\newcommand{\link}{\operatorname{link}}
\newcommand{\gr}{\operatorname{gr}}
\newcommand{\Res}{\operatorname{Res}}
\newcommand{\PGL}{\operatorname{PGL}}
\newcommand{\Sing}{\operatorname{Sing}}
\newcommand{\dlog}{\operatorname{dlog}}
\newcommand{\relint}{\operatorname{relint}}

\newcommand{\diag}{\operatorname{diag}}

\newcommand{\cX}{\mathcal{X}}

\title{Positive Geometry of Oriented Matroid Strata}
\author{Arne Kuhrs}
\address{Max Planck Institute for Mathematics in the Sciences, Leipzig, Germany}
\email{arne.kuhrs@mis.mpg.de}

\author{Artyom Lisitsyn}
\address{Department of Physics and Astronomy, University of California, Davis, USA}
\email{alisitsyn@ucdavis.edu}
\date{\today}

\begin{document}

\maketitle

\begin{abstract}
	We study the positive geometry of Grassmannians whose boundary is a union of Pl\"ucker
	divisors. For the full Pl\"ucker boundary the regions are closures of realization spaces of
	uniform oriented matroids. In rank three we show that the pair $(\Gr_{3,n},Z_{3,n})$ is a
	positive arrangement in the sense of Sturmfels and Telen for $n\le 9$ points and compute its
	combinatorial rank from point counts of the configuration space $\cX_{3,n}$ over finite
	fields. For seven points the canonical forms of the regions span all logarithmic forms. For
	eight points they do not, because the M\"obius--Kantor configuration is not realizable over
	the reals. This determines for seven and eight points whether the integrands of Cachazo,
	Early and Zhang for the generalized biadjoint amplitudes of Cachazo, Early, Guevara and
	Mizera span the space of logarithmic forms. From ten points on, the genus is positive
	because a K3 surface appears in the boundary. For smaller boundaries and arbitrary rank we
	give a recursion that adds one point at a time and preserves the genus and the combinatorial
	rank. It recovers the positive Grassmannian and shows that realization spaces of stacked
	polytopes are regions of positive arrangements of combinatorial rank one.
\end{abstract}

\section{Introduction}

Positive geometry studies semialgebraic sets that carry a distinguished differential form, the \emph{canonical form}, whose poles lie on the boundary and whose residues are again canonical forms \cite{Arkani-Hamed:2017tmz}. Arkani-Hamed, Bai and Lam define these forms recursively through their residues. Brown and Dupont recast the theory in terms of the mixed Hodge structure of a pair $(X,Y)$ of a complex variety and a divisor \cite{Brown:2025jjg}: every relative cycle has a unique canonical form when an invariant of the pair, the \emph{genus}, vanishes, and these forms span a space whose dimension is the \emph{combinatorial rank}. This turns positivity into a property of complex varieties that can be tested systematically. Sturmfels and Telen introduce positive arrangements and their real combinatorial rank \cite{Sturmfels:2026nzc}.

The Grassmannian is a central source of examples, and many of them come from physics. The totally nonnegative Grassmannian $\Gr^{\ge0}_{k,n}$ of Postnikov \cite{Postnikov:2006kva} is a positive geometry \cite{Arkani-Hamed:2017tmz}, and its image, the amplituhedron of Arkani-Hamed and Trnka, computes amplitudes in planar $\mathcal N=4$ super Yang--Mills theory \cite{Arkani-Hamed:2012zlh,Arkani-Hamed:2013jha}. For $k=2$, the regions of $\Gr_{2,n}(\R)$ cut out by all Pl\"ucker divisors are, up to a torus, the regions of $M_{0,n}(\R)$. Their canonical forms are the Parke--Taylor forms \cite{Arkani-Hamed:2017mur}, which give the biadjoint scalar amplitudes through the scattering equations \cite{Cachazo:2013iea}. Cachazo, Early, Guevara and Mizera (CEGM) replaced $M_{0,n}$ by the configuration space $\cX_{k,n}$ of $n$ points in $\PP^{k-1}$ to define generalized biadjoint amplitudes \cite{Cachazo:2019ngv}. Their study led to connections with tropical Grassmannians and cluster algebras \cite{Cachazo:2019apa,Drummond:2019qjk,Henke:2019hve} and with the positive part of $\cX_{k,n}$ \cite{Arkani-Hamed:2020cig,Arkani-Hamed:2019rds}. In all of these examples except the amplituhedron, the boundary is made of Pl\"ucker divisors.

For a set $S\subseteq\binom{[n]}{k}$ we put
\begin{equation}\label{eq:ZS}
	Z_S\ :=\ \bigcup_{I\in S}\{p_I=0\}\ \subset\ \Gr_{k,n},
	\qquad
	Z_{k,n}:=Z_{\binom{[n]}{k}},
\end{equation}
and we ask when the regions of $\Gr_{k,n}(\R)$ cut out by $Z_S$ have canonical forms, and how large the space is that these forms span. Two extreme cases are well understood. The first is the full Pl\"ucker boundary for $k=2$: the open part of the pair $(\Gr_{2,n},Z_{2,n})$ is, up to a torus, the moduli space $M_{0,n}$ (\zcref{thm:guide-M0n}). The second is the positive Grassmannian, which exists for every $k$ and $n$ but has only the $n$ cyclic Pl\"ucker divisors in its boundary (\zcref{thm:guide-posgr}). This paper pursues two directions for finding new positive arrangements $(\Gr_{k,n},Z_S)$ between them. To state the two guiding theorems and our results, we first recall the framework.

\subsection*{Positive arrangements and the two combinatorial ranks}

We work in the framework of Brown and Dupont \cite{Brown:2025jjg} and of Sturmfels and Telen \cite{Sturmfels:2026nzc}. Let $X$ be an irreducible complex projective variety of dimension $d$ and $Y\subsetneq X$ a divisor with $U=X\setminus Y$ smooth. The \emph{genus} $g(X,Y)$ and the \emph{combinatorial rank} $\cromb(X,Y)$ are Hodge numbers of the relative homology $H_d(X,Y)$, and both depend only on $U$ \cite[(8), (16)]{Brown:2025jjg}. If the genus vanishes, every relative cycle $\sigma$ has a canonical form $\omega(\sigma)$ in the space $\Omega^d_{\log}(U)$ of logarithmic top forms on $U$, which has dimension $\cromb(X,Y)$. The canonical form $\omega(\sigma)$ is determined by its iterated residues at the corners of a normal crossing compactification \cite[\S2, Cor.~2.18]{Brown:2025jjg}. The pair $(X,Y)$ is a \emph{positive arrangement} if it has genus zero, $X$ and every irreducible component of $Y$ are defined over $\R$ and have a smooth real point, and $\Sing(X)\subseteq Y$ \cite[Def.~6.2]{Sturmfels:2026nzc}. Its \emph{regions} are the closures of the orientable connected components of $X(\R)\setminus Y(\R)$, and the \emph{real combinatorial rank} $\cromb_\R(X,Y)$ is the dimension of the span of their canonical forms \cite[Def.~6.4, Def.~6.10]{Sturmfels:2026nzc}. Always $\cromb_\R\le\cromb$, and equality means that every logarithmic form is a combination of canonical forms of regions.

In this language the two extreme cases read as follows, and these two known results guide our investigation.

\begin{guidingtheorem}[Moduli of genus zero curves; {\cite[Cor.~2.2.16, Cor.~2.3.11]{Koefler}, \cite[Prop.~6.2]{Brown:2025jjg}}]\label{thm:guide-M0n}
	For every $n\ge4$ the pair $(\Gr_{2,n},Z_{2,n})$ is a positive arrangement, and
	\[
	\cromb_\R\bigl(\Gr_{2,n},Z_{2,n}\bigr)=\cromb\bigl(\Gr_{2,n},Z_{2,n}\bigr)=\cromb(\cX_{2,n})=\dim\widetilde H_{n-4}\bigl(\link\Trop\cX_{2,n}\bigr)=(n-2)!\,.
	\]
\end{guidingtheorem}

\noindent Here $\cX_{2,n}=M_{0,n}$ is the complement of a real hyperplane arrangement, in which case the two ranks agree. This follows from  \cite[Prop.~6.2, \S 6.8]{Brown:2025jjg}. The torus factor of $\Groo_{2,n}$ changes neither the genus nor the combinatorial rank \cite[Prop.~2.3.10, Prop.~2.2.19]{Koefler}, and by the same argument as in \zcref{prop:descent-real} it does not change the real combinatorial rank. We write $\link\Trop\cX$ for the link of a tropical variety, its intersection with a sphere around the origin. That it computes the combinatorial rank is shown in \cite[Cor.~3.3]{Hacking} and \cite[Prop.~2.2.12]{Koefler}.

\begin{guidingtheorem}[Positive Grassmannian;
	{\cite[Prop.~3.2]{Koefler:2026vkb}}]\label{thm:guide-posgr}
	Let $\Gr^{\ge0}_{k,n}$ be the nonnegative Grassmannian, the set of real points of
	$\Gr_{k,n}$ whose Pl\"ucker coordinates are all nonnegative. Let $S^+_{k,n}$ be the set of
	cyclic intervals $\{i,i+1,\dots,i+k-1\}$ in $[n]$. Then $Z_{S^+_{k,n}}$ is the algebraic
	boundary of $\Gr^{\ge0}_{k,n}$, that is, the Zariski closure of its Euclidean boundary. The
	pair $(\Gr_{k,n},Z_{S^+_{k,n}})$ is a positive arrangement for all $k<n$, and
	\[
	\cromb_\R\bigl(\Gr_{k,n},Z_{S^+_{k,n}}\bigr)=\cromb\bigl(\Gr_{k,n},Z_{S^+_{k,n}}\bigr)=1
	.
	\]
\end{guidingtheorem}

\noindent Genus zero is shown in \cite[Prop.~3.2]{Koefler:2026vkb}. We prove the full statement, including the two ranks, by a recursion on $n$ in \zcref{cor:positive-grassmannian-rank}.

Our two main directions of study in this paper start from these two theorems. The first keeps the full Pl\"ucker boundary of \zcref{thm:guide-M0n} and passes from rank two to rank three. The second keeps arbitrary $k$ and $n$ as in \zcref{thm:guide-posgr} and enlarges the class of boundaries.

\subsection*{Full Pl\"ucker boundary for $k=3$}

The complement $\Groo_{k,n}$ of the full Pl\"ucker boundary is a trivial torus bundle over the configuration space $\cX_{k,n}$ of $n$ points in linearly general position in $\PP^{k-1}$. A real configuration has a chirotope, the sign vector of its maximal minors, and grouping configurations by their chirotope up to reorientation decomposes $\cX_{k,n}(\R)$ into \emph{oriented matroid strata}, the realization spaces of uniform oriented matroids. By Mn\"ev's universality theorem these strata can be arbitrarily complicated \cite{mnev}, \cite[\S8.6]{BLVSWZ}. In rank three and for $n\le9$, however, every nonempty stratum is connected \cite{richter1988kombinatorische}, so the regions are the closures of the strata. Canonical forms of oriented matroids also appear in work of Eur and Lam \cite{EurLam}, in a
different sense: there an oriented matroid is fixed, and its topes, the regions of a single arrangement, receive canonical forms in the Orlik--Solomon algebra. Here each oriented
matroid indexes a region of the configuration space, namely its realization space, and the canonical form lives on the configuration space.
For five points nothing new happens, since the duality $\Gr_{3,5}\cong\Gr_{2,5}$ identifies $\cX_{3,5}$ with $\cX_{2,5}=M_{0,5}$. From $n=6$ points on, $\cX_{3,n}$ is no longer the complement of a hyperplane arrangement, but of an arrangement of hyperplanes and quadrics.
 
 The same regions appear in physics. For $n\le8$, Cachazo, Early and Zhang attach to each region
 of $\cX_{3,n}(\R)$ a rational form, the integrand of a generalized biadjoint amplitude in the formulation of \cite{Cachazo:2019ngv}. They leave open whether the regions other than the positive one are positive geometries, and hence how to define their canonical forms
 \cite[\S5.2]{Cachazo:2023ltw}. They also ask for the dimension of the space spanned by the integrands. Further, citing Skorobogatov \cite{Skorobogatov}, they note that the numbers of points of $\cX_{3,n}$ over finite fields, evaluated at $q=0$, give $126$ for six and $7470$ for seven points \cite[\S7.2]{Cachazo:2023ltw}.
 Our first three theorems address these questions.

\begin{maintheorem}\label{thm:main-posarr}
	For $6\le n\le9$ the pair $(\Gr_{3,n},Z_{3,n})$ is a positive arrangement, and its combinatorial rank is
	\[
	\cromb(\Gr_{3,n},Z_{3,n})=\cromb(\cX_{3,n})=\dim\Omega^{2n-8}_{\log}(\cX_{3,n})=
	\begin{cases}
		126 & \text{for } n=6,\\
		7470 & \text{for } n=7,\\
		1\,116\,480 & \text{for } n=8,\\
		392\,936\,880 & \text{for } n=9.
	\end{cases}
	\]
\end{maintheorem}

So the regions are positive geometries in the sense of \cite[Def.~6.4]{Sturmfels:2026nzc}, and for $n\le8$ their canonical forms are, up to sign and a common factor, the integrands of Cachazo, Early and Zhang (\zcref{remark:cegm_integrands}). The proof of \zcref{thm:main-posarr} in \zcref{ss:finitefield} uses point counts over finite fields. For a smooth affine variety the genus and the combinatorial rank are read off from the $E$-polynomial (\zcref{prop:ep0}). Then by a theorem of Katz \cite[Thm.~6.1.2(3)]{Katz}, if the number of points over $\F_q$ is a polynomial in $q$ for all prime powers $q\equiv1\pmod m$, then the $E$-polynomial is this polynomial evaluated at $xy$. So the genus vanishes, and the combinatorial rank is the value of the count at $q=0$. The counts of $\cX_{3,n}(\F_q)$ are due to Glynn \cite{Gly88} and to Iampolskaia, Skorobogatov and Sorokin \cite{ISS95}. Away from characteristics two and three they are polynomials for six and seven points and quasi-polynomials, that is, polynomials on residue classes of $q$, for eight and nine points. Lam evaluates them at $q=1$ to obtain the Euler characteristic of $\cX_{3,n}$ \cite[Thm.~A.1]{Agostini:2021rze}. At $q=0$ they give the combinatorial rank, which explains the observation of Cachazo, Early and Zhang. For eight points the count contains a term $840\,b(q)$, where $b(q)$ is the number of primitive cube roots of unity in $\F_q$. It comes from the M\"obius--Kantor configuration, which reappears in the analysis of the real combinatorial rank of $\cX_{3,8}$ in \zcref{thm:main38}.

It remains to compare the span of the canonical forms of the regions with the space of all logarithmic forms. For $\cX_{3,6}$ the two ranks are both $126$ \cite[\S12]{Sturmfels:2026nzc}, \cite[Prop.~5.3]{HPP}. For seven points they agree as well, which confirms the expected dimension $7470$ (see
\cite[\S7.2]{Cachazo:2023ltw} and \cite[\S5.2]{HPP}).

\begin{maintheorem}\label{thm:main37}
	The real combinatorial rank of seven points equals the combinatorial rank:
	\[
	\cromb_{\R}\bigl(\Gr_{3,7},Z_{3,7}\bigr)=\cromb\bigl(\Gr_{3,7},Z_{3,7}\bigr)=\cromb\bigl(\cX_{3,7}\bigr)=\dim\widetilde H_5\bigl(\link\Trop\cX_{3,7}\bigr)=7470 .
	\]
\end{maintheorem}

With \zcref{remark:cegm_integrands} it follows that the integrands for seven points span the whole space of logarithmic forms.
We prove \zcref{thm:main37} in \zcref{ss:real37} with real tropical data. After passing from the Grassmannian to the configuration space (\zcref{prop:descent-real}), we use a tropical compactification of $\cX_{3,7}$. By a theorem of Corey every initial degeneration of $\cX_{3,7}$ is smooth and irreducible \cite[Thm.~1.2, Lemma~7.1]{Corey}, so there is one corner over every maximal cone of the tropical variety, and every corner is real. Near a corner the complement of the boundary consists of $64$ orthants, which lie in $64$ different strata. Using the theorem $\Dr^\chi(3,7)=\Trop^\chi\Gr(3,7)$ of Antolini and Early \cite[Thm.~1.1]{Early:2024asu}, we show that the strata that reach a corner are those whose chirotropical Dressian contains the cone of the corner (\zcref{prop:orthants-dressian}). So modulo $2$ the iterated boundaries of the regions form the incidence matrix of the $27\,240$ chirotropical Dressians and the $211\,365$ Pl\"ucker cones of dimension six, which we build from the fans of Antolini and Early. Its rank over $\F_2$ is $7470$. 

For eight points the two ranks differ, and the reason is a configuration that exists over $\C$ but not over $\R$.

\begin{maintheorem}\label{thm:main38}
	The real combinatorial rank of eight points is smaller than the combinatorial rank:
	\[
	\cromb_{\R}\bigl(\Gr_{3,8},Z_{3,8}\bigr)+840\ \le\ \cromb\bigl(\Gr_{3,8},Z_{3,8}\bigr)=1\,116\,480 .
	\]
\end{maintheorem}

By Corey and Luber, the initial degenerations of $\cX_{3,8}$ have two components over the $840$ cones of the $S_8$-orbit of the M\"obius--Kantor cone, one for each of the two complex conjugate realizations of the M\"obius--Kantor matroid \cite[Thm.~1.3, Thm.~1.4, Thm.~5.2]{CoreyLuber}. The corresponding corners of a tropical compactification are not real, so no region reaches them. The difference of the two corners over each cone gives a logarithmic form, and the span of these $840$ forms meets the span of the canonical forms of the regions only in zero. So for eight points the integrands span a subspace of codimension at least $840$ in $\Omega^8_{\log}(\cX_{3,8})$. The exact value of the real combinatorial rank of eight points remains open.

From ten points on, the genus is positive.

\begin{maintheorem}\label{thm:main310}
	For every $n\ge10$ we have $g(\Gr_{3,n},Z_{3,n})=g(\cX_{3,n})\ge1$. In particular,
	$(\Gr_{3,n},Z_{3,n})$ is not a positive arrangement.
\end{maintheorem}

The reason is a K3 surface in the boundary for ten points. For a certain configuration of ten lines through ten points, the realization space of ten points in which exactly these triples are collinear is a dense open subset of a K3 surface \cite[Thm.~1.2]{Sink-K3surface}, building on \cite{Sturmfels1991}. In \zcref{sec:main310} we connect $\cX_{3,10}$ to this realization space by a chain of strata, each obtained from the previous one by imposing one collinearity. Each stratum together with the next one is a smooth affine variety, so the genus can only drop in each step \cite[Prop.~3.10]{Brown:2025jjg}, and the genus of $\cX_{3,10}$ is at least that of the K3 surface, which is one. Adding a point cannot lower the genus (\zcref{prop:genus-monotone}), and with the duality of
Grassmannians we obtain $g(\cX_{k,n})\ge1$ whenever $3\le k\le n-7$ or $7\le k\le n-3$
(\zcref{cor:genus-range}). By \zcref{cor:pointcount}, the number of $n$-arcs in $\PP^2(\F_q)$ is
therefore not a quasi-polynomial in $q$ for $n\ge10$. For ten points this was announced, and for
more points conjectured, in \cite[\S1.1]{ISHAM2022102006}.

\subsection*{Smaller boundaries and polytope realization spaces}

By \zcref{thm:main310} the full Pl\"ucker boundary does not always give a positive arrangement. In the second direction we keep $k$ and $n$ arbitrary and look for smaller boundaries $Z_S\subsetneq Z_{k,n}$ that do. We build them by adding one point at a time, and for this we need one combinatorial notion.

A $k$-nesting $\Gamma=(I_1,\dots,I_k)$ on $[n-1]$ is a tuple of $(k-1)$-subsets that satisfies a combinatorial condition (\zcref{def:nesting}). Each $I_r$ spans a hyperplane $H_r$ through $k-1$ of the first $n-1$ points. The set $N(\Gamma)\subseteq\binom{[n-1]}{k}$ lists Pl\"ucker coordinates whose nonvanishing makes $H_1,\dots,H_k$ independent. The set $D(\Gamma)=\{I_r\cup n:r\in[k]\}$ lists those whose nonvanishing says that the new point $n$ lies off every $H_r$ (see \zcref{fig:k3-extensions}). The new point is then recorded by $k$ nonzero coordinates relative to these hyperplanes, which gives a torus factor.

\begin{maintheorem}\label{thm:mainrecursion}
	Let $2\le k<n-1$, let $\Gamma$ be a $k$-nesting on $[n-1]$, and let $S\subseteq\binom{[n-1]}{k}$ with $N(\Gamma)\subseteq S$. If $(\Gr_{k,n-1},Z_S)$ is a positive arrangement, then so is $(\Gr_{k,n},Z_{S'})$ for every $S'$ with $(S\setminus N(\Gamma))\cup D(\Gamma)\subseteq S'\subseteq S\cup D(\Gamma)$, and
	\[
	\cromb(\Gr_{k,n},Z_{S'})=\cromb(\Gr_{k,n-1},Z_S).
	\]
\end{maintheorem}

The proof in \zcref{sec:smaller-boundaries} has two steps. First, every point of $\Gr_{k,n-1}\setminus Z_{N(\Gamma)}$ has a unique matrix representative with prescribed zeros and ones, so that the last column can be chosen freely in $(\C^*)^k$. This gives an isomorphism $(\Gr_{k,n-1}\setminus Z_S)\times(\C^*)^k\cong\Gr_{k,n}\setminus Z_{S\cup D(\Gamma)}$ (\zcref{lem:nesting-isomorphism}), and a torus factor changes neither the genus nor the combinatorial rank \cite[Prop.~2.3.10, Prop.~2.2.19]{Koefler}. Second, the divisors $Z_{N(\Gamma)}$ can be removed from the boundary. Shrinking the boundary can only lower the genus \cite[Prop.~3.10]{Brown:2025jjg}, so the genus stays zero. Pushing forward a logarithmic form along the isomorphism multiplies it by $\prod_{r=1}^{k-1}p_{I_r\cup j_r}/\prod_{I\in D(\Gamma)}p_I$, where the sets $I_r\cup j_r$ are the elements of $N(\Gamma)$ (\zcref{lem:form-recursion}). The numerator cancels every pole along $Z_{N(\Gamma)}$, so every logarithmic form extends to the larger space. For the nesting formed by the cyclic intervals next to the new point, this recursion recovers \zcref{thm:guide-posgr}.

Further boundaries of this kind arise from convex geometry. For a simplicial $(k-1)$-polytope $K$ with $n$ labeled vertices, the Grassmannian model of its realization space is the semialgebraic set $\Gr_K^{>0}\subset\Gr_{k,n}(\R)$ of points whose columns, read as homogenized vertices, realize $K$ \cite{GMW-Polytope}. It is cut out by sign conditions on the Pl\"ucker coordinates $p_{F\cup v}$, for a facet $F$ and a vertex $v\notin F$, which record on which side of $F$ the vertex $v$ lies. For the $n$-gon with vertices labeled in cyclic order it is the positive Grassmannian $\Gr^{>0}_{3,n}$. We write $\Gr_K^{\ge0}$ for its closure. The \emph{polytopal divisor} $Z_K$ is the union of the divisors $Z_{\{F\cup F'\}}$ over pairs of adjacent facets $F,F'$, along which the two facets become coplanar. A polytope is \emph{stacked} if it is obtained from a simplex by repeatedly adding a new vertex just beyond one facet.

\begin{maintheorem}\label{thm:mainpolytope}
	Let $K$ be a stacked $(k-1)$-polytope with $n$ vertices. Then $(\Gr_{k,n},Z_K)$ is a positive arrangement, $\Gr_K^{\geq 0}$ is one of its regions, and
	\[
	\cromb_\R(\Gr_{k,n},Z_K)\ =\ \cromb(\Gr_{k,n},Z_K)\ =\ 1 .
	\]
\end{maintheorem}

That $(\Gr_{k,n},Z_K)$ is a positive arrangement with $\cromb=1$ follows from \zcref{thm:mainrecursion} by stacking one vertex at a time, starting from $\Gr_{k,k+1}\cong\PP^k$ with its coordinate hyperplanes. 

\subsection*{Organization}

\zcref{sec:background} recalls the framework. \zcref{ss:finitefield} proves \zcref{thm:main-posarr}, \zcref{sec:realrank} proves \zcref{thm:main37,thm:main38}, and \zcref{sec:main310} proves \zcref{thm:main310}. \zcref{sec:smaller-boundaries} introduces nestings and proves \zcref{thm:mainrecursion}, and \zcref{sec:polytopes} treats polytopes and proves \zcref{thm:mainpolytope}.

\subsection*{Code and data}

The code and data for the computations in this paper are available at \cite{KL26data}. They consist of the construction of the matrix $B$ and the computation of its rank in the proof of \zcref{thm:main37}, and of the verification of the iterated residue of the explicit $12$-form on $\cX_{3,10}$ in \zcref{rem:X310-form}.

\subsection*{Acknowledgements}
We thank Bernd Sturmfels for suggesting this project. We also thank Dario Antolini, Nick Early, Ben Hollering, Joris Koefler, Dante Luber,
Melvyn Nabavi, Ana\"elle Pfister, Lizzie Pratt, Simon Telen, Martin Winter and Claudia Yun
for very useful discussions.

\subsection*{AI use disclosure}
This paper was prepared with substantial assistance from LLM tools, primarily Claude Opus 5.5 used with the Claude Code harness. It was used for literature search, provided mathematical suggestions, helped to simplify arguments, to draft exposition and to produce some of the figures, and wrote most of the code in \cite{KL26data}. The authors are responsible for all assertions made in the paper.

\section{Background}\label{sec:background}

In this section we set up the framework. Our objects are the pairs $(\Gr_{k,n}, Z_S)$ of a Grassmannian and a union of Plücker divisors. The open part of $(\Gr_{k,n}, Z_{k,n})$ is a trivial torus bundle over the configuration space $\cX_{k,n}$ of $n$ points in $\P^{k-1}$. We make this explicit in coordinates. The real points of $\cX_{k,n}$ are organized by oriented matroids. We then recall the genus and the combinatorial rank of Brown and Dupont and the positive arrangements of Sturmfels and Telen. 

\subsection{The Grassmannian, the configuration space and oriented matroid strata}\label{ss:objects}
 
Let $\Gr_{k,n}$ be the Grassmannian of $k$-planes in $\C^n$, embedded in $\PP^{\binom nk-1}$ by the Pl\"ucker coordinates $(p_I)_{I\in\binom{[n]}{k}}$. Let $Z_I = \{p_I = 0\} \cap \Gr_{k,n}$, and for $S \subseteq \binom{[n]}{k}$,
\[
  Z_S:= \bigcup_{I \in S }Z_I,\ Z_{k,n} := Z_{\binom{[n]}{k}},\   \Groo_{k,n}\ =\ \Gr_{k,n}\setminus Z_{k,n}.
\]

A point of $\Gr_{k,n}$ can be represented by the row span of a $k\times n$ matrix $A$ of rank $k$ and its Plücker coordinates are the maximal minors $p_I(A) = \det A_I$.  Two matrices have the same row span exactly when they differ by left multiplication by an element of $\mathrm{GL}_k$. The torus $(\Gm)^n$ acts on $\Gr_{k,n}$ by scaling columns. For $\tau \in \G_m^n$, 
\begin{equation}\label{eq:plucker-scaling}
    p_I(A \diag(\tau)) = \bigl(\prod_{i \in I} \tau_i\bigr)p_I(A),
\end{equation}
so the diagonal torus $\G_m$ inside $\G_m^n$ acts trivially. Further, the torus $T=(\Gm)^n/\Gm\cong(\Gm)^{n-1}$ acts freely on $\Groo_{k,n}$. The quotient
\[
  \cX_{k,n}\ =\ \Groo_{k,n}/T
\]
is the \emph{configuration space} of $n$ points in $\PP^{k-1}$ in linearly general position, that is, no $k$ of them lie on a hyperplane, modulo $\PGL_k$. The $n$ points are the columns of a $k\times n$ matrix $A$ whose row span is a point of $\Groo_{k,n}$. The space $\cX_{k,n}$ has dimension $d = (k-1)(n-k-1)$, and for $k=2$ it is the moduli space $M_{0,n}$ of $n$ distinct points on $\P^1$. 
Let $T_P = \G_m^{\binom{n}{k}} / \G_m$ be the Plücker torus inside $\P^{\binom{n}{k}-1}$. Then $\cX_{k,n}$ is a closed subvariety of the torus $T_P / T$. In particular, it is very affine. 
 
In this article we focus on the case $k=3$, and we make the relation between $\Groo_{3,n}$ and $\cX_{3,n}$ explicit only in this case; the same construction works for every $k\ge2$, with the first $k+1$ columns normalized to the standard projective frame. Every point of $\cX_{3,n}$ has a unique representative of the form \cite[Equation (2.6)]{Cachazo:2023ltw}
 \begin{equation}\label{eq:normalform}
 	M(y) =\begin{pmatrix}1&0&0&1&y_{15}&y_{16}&\dots&y_{1n}\\0&1&0&1&y_{25}&y_{26}&\dots&y_{2n}\\0&0&1&1&1&1&\dots&1\end{pmatrix}.
 \end{equation}
 Here a projective transformation moves the first four points to the standard projective frame, and the remaining columns are scaled to have last entry $1$. 
 The $d = 2n-8$ entries are the cross-ratios of Plücker coordinates,
 \[
 y_{1a}=\frac{p_{23a}\,p_{124}}{p_{12a}\,p_{234}},\qquad
 y_{2a}=\frac{p_{13a}\,p_{124}}{p_{12a}\,p_{134}}\qquad(5\le a\le n).
 \]
 Every index occurs the same number of times in the numerator and in the denominator, so by \eqref{eq:plucker-scaling} they are $T$-invariant. They are regular and nowhere zero on $\Groo_{3,n}$, hence they descend to $\cX_{3,n}$. Further, they identify $\cX_{3,n}$ with the open subset of $\A^d$ on which all $3 \times 3$ minors of $M(y)$ are nonzero. In these coordinates the quotient map $\pi: \Groo_{3,n} \to \cX_{3,n}$ is $g \mapsto (y_{1a}(g),y_{2a}(g))_{5 \leq a \leq n}$ and $y \mapsto M(y)$ is a section of $\pi$.
 
 For the torus we use the $(n-1)$ Pl\"ucker ratios
 \[
 t_1=\frac{p_{124}}{p_{234}},\qquad t_2=-\frac{p_{124}}{p_{134}},\qquad
 t_a=\frac{p_{12a}}{p_{123}}\quad(4\le a\le n).
 \]
 Each is a ratio of two Pl\"ucker coordinates, so it does not depend on the matrix representative, and it is a unit on $\Groo_{3,n}$. By \eqref{eq:plucker-scaling} they record the scaling factors,
 \begin{equation}\label{eq:t-equivariance}
 	t_i(\tau\cdot g)=\frac{\tau_i}{\tau_3}\,t_i(g)\qquad(i\in[n]\setminus\{3\})
 \end{equation}
 where $\tau \in \G_m^n, g \in \Groo_{3,n}$. On the normal form we have $p_{123}=p_{124}=p_{234}=p_{12a}=1$ and $p_{134}=-1$, so that $t_i(M(y))=1$ for all
 $i\ne3$. 
 
 These coordinates give an explicit trivialization: the morphism
 \begin{equation}\label{eq:trivialisation}
 	\Phi:=(\pi,\mathbf t)\colon\Groo_{3,n}\longrightarrow\cX_{3,n}\times\Gm^{n-1},
 \end{equation}
 where $\mathbf t=(t_i)_{i\ne3}$, is an isomorphism. We here identify $T$ with $\Gm^{n-1}$ via $\tau\mapsto\tau_i/\tau_3$. The inverse of $\Phi$ sends $(y,t)$ to $t\cdot M(y)$, the normal form with its $i$-th column scaled by $t_i$ for $i\ne3$. On the level of real points, the isomorphism $\Phi$ identifies $\Groo_{3,n}(\R)$ with $\cX_{3,n}(\R) \times (\R^\times)^{n-1}$.

The sign vector $(\sign p_I(g))_I$ of a point $g \in \Groo_{3,n}(\R)$, defined up to an overall sign, is called a uniform \emph{chirotope} of rank $3$ on $[n]$ and it is constant on each connected component of $\Groo_{3,n}(\R)$. The torus $T(\R)$ changes it by reorientations. For a uniform chirotope $\chi$ of rank $3$ on $[n]$ let $S_\chi \subseteq \cX_{3,n}(\R)$ be the set of configurations whose chirotope lies in the reorientation class of $\chi$, the \emph{oriented matroid stratum} of $\chi$. Every connected component of $\cX_{3,n}(\R)$ lies in a single stratum. The strata are the realization spaces of uniform oriented matroids.  For $k = 3$ and $n \leq 9$ conversely every nonempty stratum is connected \cite{richter1988kombinatorische}. These oriented matroid strata are the real objects whose positive geometry we study.

\subsection{Positive arrangements, genus, and the combinatorial rank}\label{ss:posarr}

We use the mixed Hodge theoretic framework of Brown and Dupont \cite{Brown:2025jjg}. Let $X$ be
a complex variety of dimension $d$ and $Y\subseteq X$ a closed subvariety. The relative homology $H_d(X,Y)$ with rational
coefficients carries a mixed Hodge structure whose Hodge numbers govern canonical forms.
\begin{definition}[{\cite[Def.~2.1, Def.~2.2]{Brown:2025jjg}}]\label{def:genus-cr}
	The \emph{genus} and the \emph{combinatorial rank} of the pair $(X,Y)$ are
	\[
	g(X,Y)=\sum_{p>0}h^{-p,0}\bigl(H_d(X,Y)\bigr),
	\qquad
	\cromb(X,Y)=h^{0,0}\bigl(H_d(X,Y)\bigr).
	\]
\end{definition}
In the rest of this subsection $X$ is compact and
irreducible, $Y\subsetneq X$ is a divisor, and $U:=X\setminus Y$ is smooth.

Let $\Omega^d_{\log}(U)$ be the space of holomorphic $d$-forms on $U$ with logarithmic poles along the boundary of a smooth compactification of $U$ with normal crossing boundary. This does not depend on the choice of compactification \cite[Def. 1.7, Prop. 1.8]{Brown:2025jjg}. Brown and Dupont construct a surjective linear map
\begin{equation}\label{eqn:Rmap}
R: \Omega^d_{\log}(U) \to \gr^W_0H_d(X,Y)_\C
\end{equation}
with kernel of dimension $g(X,Y)$, so
\[
\dim \Omega^d_{\log}(U) = \cromb(X,Y) + g(X,Y).
\]

If $g(X,Y) = 0$, then $R$ is an isomorphism, and composing its inverse with the quotient map $H_d(X,Y) \to \gr^W_0H_d(X,Y)$ gives the \emph{canonical form} map
\[
\mathrm{can}: H_d(X,Y) \to \Omega^d_{\log}(U).
\]
The relative homology of the pair is the locally finite homology of $U$ \cite[Eq. (8)]{Brown:2025jjg}, and Poincar\'e duality is an isomorphism of mixed Hodge structures $H_d(X,Y) \cong H^d(U)(d)$. Hence the genus, the combinatorial rank and the canonical form map depend only on $U$. We write $g(U):=g(X,Y)$ and $\cromb(U):=\cromb(X,Y)$, where $X$ is any compactification of $U$ and $Y=X\setminus U$. Note that this differs from the notation of \cite{Brown:2025jjg}, where $g(U)$ stands for $g(U,\varnothing)$.

Real regions enter through the positive arrangements of Sturmfels and Telen \cite[\S 6]{Sturmfels:2026nzc}.
 
\begin{definition}[{\cite[Def.~6.2]{Sturmfels:2026nzc}}]\label{def:posarr}
Let $X$ be an irreducible complex projective variety and $Y$ a divisor on $X$. The pair $(X,Y)$ is a \emph{positive arrangement} if
\begin{enumerate}
  \item $X$ is defined over $\R$ and has a smooth real point,
  \item each irreducible component of $Y$ is defined over $\R$ and has a smooth real point,
  \item $X\setminus Y\subseteq X_{\mathrm{reg}}$, equivalently $\Sing(X)\subseteq Y$,
  \item the pair $(X,Y)$ has genus zero.
\end{enumerate}
\end{definition}
 
A \emph{region} of a positive arrangement is the closure $\sigma$ in $X(\R)$ of an orientable connected component $\sigma^{\circ}$ of $X(\R)\setminus Y(\R)$. An orientation of $\sigma^{\circ}$ defines a relative cycle with a class $[\sigma] \in H_d(X,Y)$, and its canonical form is $\omega(\sigma) := \mathrm{can} [\sigma]$. The \emph{algebraic boundary} $\partial_a\sigma$ of a region $\sigma$ is the Zariski closure in $X$ of its Euclidean boundary.  

For the Grassmannian with its full Plücker boundary, both invariants can be computed on the configuration space.

\begin{proposition}\label{cor:genus_and_cr_descent}
	For all $2 \leq k \leq n-2$ we have $g(\Gr_{k,n}, Z_{k,n}) =  g(\Groo_{k,n}) =   g(\cX_{k,n})$ and $\cromb(\Gr_{k,n},Z_{k,n}) = \cromb(\Groo_{k,n})  =  \cromb(\cX_{k,n}).$
\end{proposition}

\begin{proof}
	Since $\Gr_{k,n}$ is compact and  $\Gr_{k,n} \setminus Z_{k,n} = \Groo_{k,n}$, the first equality in each chain holds because the invariants only depend on the open part. We have an isomorphism $\Groo_{k,n} \cong \cX_{k,n} \times \G_m^{n-1}$, thus by \cite[Prop. 2.2.19, Prop. 2.3.10]{Koefler} the genus and combinatorial rank do not change when a torus factor is removed.
\end{proof}

\section{Genus and combinatorial rank via finite field counts} \label{ss:finitefield}

In this section we compute the genus and the combinatorial rank by using point counts over finite fields for $\cX_{3,n}$. We use the $E$-polynomial of a complex variety, which is a polynomial in two variables $x,y$ whose coefficients $e_{p,q}$ in front of $x^py^q$ record a signed sum of Hodge numbers $h^{p,q}$. We observe that for a smooth affine variety the coefficients $e_{p,0}$ of the monomials $x^py^0$ are the Hodge numbers $h^{-p,0}$ of $H_d(X,Y)$. So they define the combinatorial rank $(p=0)$ and the genus $(p > 0)$. If the variety has strongly polynomial count over finite fields, Katz's theorem identifies its $E$-polynomial with the counting polynomial evaluated at $xy$. Then the genus vanishes and the combinatorial rank is the constant term of the count, up to sign. The known point count formulas for $\cX_{3,n}(\F_q)$ which satisfy the assumptions of Katz's theorem then show that $(\Gr_{3,n}, Z_{3,n})$ is a positive arrangement for $6 \leq n \leq 9$ and compute its combinatorial rank.

\subsection{The $E$-polynomial}

Let $U$ be a complex variety. Its \emph{E-polynomial} is 

\[
E(U;x,y)=\sum_{p,q}e_{p,q}(U)\,x^py^q,\qquad
e_{p,q}(U)=\sum_k(-1)^k\,h^{p,q}\bigl(\gr^W_{p+q}H^k_c(U;\C)\bigr),
\]
where $h^{p,q}$ are the Hodge numbers of Deligne's mixed Hodge structure on compactly supported cohomology. For further background we refer to \cite[\S 6]{Katz}. 

The $E$-polynomial knows the genus and combinatorial rank if $U$ is smooth and affine.

\begin{proposition}\label{prop:ep0}
	Let $U$ be a smooth affine irreducible complex variety of dimension $d$, and let $(X,Y)$ be a compactification of $U$. Then for every integer $p \geq 0$
	\[
	e_{p,0}(U) = (-1)^d h^{-p,0}(H_d(X,Y)).	
	\]
	In particular, 
	\[
	\cromb(U)=(-1)^d\,E(U;0,0),\qquad
	g(U)=(-1)^d\sum_{p\geq1}e_{p,0}(U).
	\]
\end{proposition}

\begin{proof}
	Since $X$ is compact, the mixed Hodge structure on $H^k_c(U)$ is that of the relative cohomology $H^k(X,Y)$ (see \cite[(8)]{Brown:2025jjg}). Since $H^k(X,Y)$ is the dual of $H_k(X,Y)$ as a mixed Hodge structure, we have $h^{p,q}(H^k(X,Y))=h^{-p,-q}(H_k(X,Y))$ for all $p,q$. Hence
	\[
	e_{p,0}(U)=\sum_k(-1)^k\,h^{-p,0}\bigl(H_k(X,Y)\bigr),
	\]
	and it suffices to show that $h^{-p,0}(H_k(X,Y))=0$ for all $k\neq d$ and $p\ge0$.

    For $k>d$, Poincar\'e duality gives $H_c^k(U)\cong H^{2d-k}(U)^\vee(-d).$ The Hodge types $(a,b)$ of $H^{2d-k}(U)$ satisfy
    $0\leq a,b\leq 2d-k$. Consequently, the Hodge types
    $(p,q)=(d-a,d-b)$ of $H_c^k(U)$ satisfy
    $p,q\geq k-d>0$. Since $q>0$, we get $h^{p,0}(H_c^k(U))=0$.
    For $k<d$ we use that $U$ is smooth and affine. Poincar\'e duality gives an isomorphism of mixed Hodge structures $H_k(X,Y)\cong H^{2d-k}(U)(d)$. Since $2d-k>d$ and $U$ is a smooth affine variety of dimension $d$, the Andreotti--Frankel theorem  \cite[Thm. 1]{AF59}  gives $H^{2d-k}(U)=0$. So $H_k(X,Y)=0$ for $k<d$, and
	\[
	e_{p,0}(U)=(-1)^d\,h^{-p,0}\bigl(H_d(X,Y)\bigr).
	\]
	By definition, $\cromb(X,Y)=h^{0,0}(H_d(X,Y))$ and $g(X,Y)=\sum_{p>0}h^{-p,0}(H_d(X,Y))$. The first formula follows from $E(U;0,0)=e_{0,0}(U)$.
\end{proof}

We can make this effective if we have a polynomial point count, possibly after restricting to a residue class of $q$.

\begin{definition}[{\cite[\S1.1]{ISHAM2022102006}}]\label{def:quasipolynomial}
	A function $f$ of the prime powers $q$ is a \emph{quasi-polynomial} if there are an
	integer $N\ge1$ and polynomials $g_0,\dots,g_{N-1}\in\Q[t]$ such that $f(q)=g_i(q)$
	whenever $q\equiv i\pmod N$.
\end{definition}

\begin{theorem}[{cf.\ Katz \cite[Thm.~6.1.2(3)]{Katz}}]\label{thm:katz}
	Let $U$ be a complex variety defined over $\Z$. Suppose that there are a polynomial $P\in
	\Q[t]$ and integers $m,N\ge1$ such that
	\begin{equation}\label{eq:polycount}
		\#U(\F_q)=P(q)\qquad\text{for all prime powers $q =p^r$ with $\gcd(q,N)=1$ and $q
			\equiv1\pmod m$.}
	\end{equation}
	Then $E(U;x,y)=P(xy)$.
\end{theorem}

\begin{remark}\label{rem:katz}
	\zcref{thm:katz} is a slightly different formulation of Katz’s theorem. Katz assumes that $U$ is
	\emph{strongly polynomial count}: $U$ has a model $U_R$ over a finitely generated subring $R\subset\C$ such that $\#U_R(\F_q)=P(q)$ for every ring homomorphism $R\to\F_q$, and he
	concludes that $E(U;x,y)=P(xy)$ \cite[Thm.~6.1.2(3)]{Katz}. To deduce \zcref{thm:katz}, take
	$R=\Z[\zeta_m,1/(mN)]$ and for $U_R$ the base change of the model of $U$ over $\Z$. A ring
	homomorphism $R\to\F_q$ exists exactly when $\gcd(q,mN)=1$ and $q\equiv1\pmod m$, and both
	conditions pass from $q$ to all its powers. So \eqref{eq:polycount} says that $U_R$ is
	strongly polynomial count with counting polynomial $P$.
\end{remark}

\begin{corollary}\label{cor:pointcount}
	Let $U$ be an affine smooth irreducible complex variety of dimension $d$, defined over $\Z$.
	Suppose that there are a polynomial $P$ and integers $m,N\ge1$ such that \eqref{eq:polycount}
	holds. 

	Then $g(U)=0$ and $\cromb(U)=(-1)^dP(0)$. In particular, if $\#U(\F_q)$ is a
	quasi-polynomial in $q$, then $g(U)=0$.
\end{corollary}

\begin{proof}
	By \zcref{thm:katz}, $E(U;x,y)=P(xy)$. Every monomial of $P(xy)$ has equal exponents in $x$
	and $y$, so $e_{p,0}(U)=0$ for $p\ge1$ and $e_{0,0}(U)=P(0)$. By \zcref{prop:ep0} this gives
	$g(U)=0$ and $\cromb(U)=(-1)^dP(0)$. If $\#U(\F_q)$ is a quasi-polynomial with modulus $N$,
	then \eqref{eq:polycount} holds with $m=N$ and $P=g_1$, since every $q\equiv1\pmod N$ is
	prime to $N$.
\end{proof}

\subsection{Configuration spaces of points}

We now apply this method to configuration spaces of points. The first example is the moduli space of points on the line.

\begin{example}
Fixing three points at $\infty, 0, 1$ identifies $M_{0,n} = \cX_{2,n}$ with the open subscheme of $\A^{n-3}$ with coordinates $x_4,\ldots,x_n$ on which $x_i \notin \{0,1\}$ and $x_i \neq x_j$ for $i \neq j$. It is smooth and affine of dimension $n-3$ and
\[
\#M_{0,n}(\F_q)=(q-2)(q-3)\cdots(q-n+2).
\]
\zcref{cor:pointcount} gives $g(M_{0,n})=0$ and $\cromb(M_{0,n})=(-1)^{n-3}\prod_{j=2}^{n-2}(-j)=(n-2)!$. This recovers the well-known values of the genus and the combinatorial rank in \zcref{thm:guide-M0n}.
\end{example} 

We now turn to $k=3$. The configuration space $\cX_{3,n}$ is a smooth very affine variety defined over $\Z$. The point counts were computed by Glynn \cite{Gly88} and by Iampolskaia, Skorobogatov and Sorokin \cite{ISS95}. We use them in the form of Lam's appendix to \cite{Agostini:2021rze}. For $q = p^r$ for a prime $p>3$
\begin{align*}
	\#\cX_{3,6}(\F_q)&=(q-2)(q-3)(q^2-9q+21),\\
	\#\cX_{3,7}(\F_q)&=(q-3)(q-5)(q^4-20q^3+148q^2-468q+498),\\
	\#\cX_{3,8}(\F_q)&=Q_8(q)+840\,b(q),\\
	\#\cX_{3,9}(\F_q)&=Q_9(q)+840\,(9q^2-243q+1684)\,b(q)+30240\,\bigl(9\,d(q)+2\,e(q)\bigr),
\end{align*}
where
\begin{align*}
	Q_8(q)={}&(q-5)(q^7-43q^6+788q^5-7937q^4+47097q^3-162834q^2+299280q-222960),\\
	Q_9(q)={}&q^{10}-75q^9+2530q^8-50466q^7+657739q^6-5835825q^5+35563770q^4\\
	&-146288034q^3+386490120q^2-588513120q+389442480,
\end{align*}
and $b(q)$, $d(q)$, $e(q)$ are the numbers of roots in $\F_q$ of $x^2+x+1$, $x^2+x-1$ and $x^2+1$. The counts for six and seven points are polynomials, so
 \eqref{eq:polycount} holds with $m=1$ and $N=6$.
 Thus evaluating the point counts at $q =0$ gives $\cromb(\cX_{3,6}) = 126$ and
$\cromb(\cX_{3,7}) = 7470$.
 
 If $q\equiv1\pmod3$ then $b(q) = 2$ and $\cX_{3,8}$ satisfies   \eqref{eq:polycount} with $m=3, N = 6$ and  
\begin{equation}\label{eq:P8}
	\begin{aligned}
		\#\cX_{3,8}(\F_q)=P_8(q):={}&Q_8(q)+1680
	\end{aligned}
\end{equation}
has constant term $P_8(0) = 1\,116\,480$.

Similarly, for $n=9$, if $q\equiv1\pmod{60}$ then $\F_q^\times$ contains the third, fourth and fifth roots of unity, hence also $\sqrt{-3}$, $\sqrt{-1}$ and $\sqrt5$ (as $\sqrt5\in\Q(\zeta_5)$), so $b(q)=d(q)=e(q)=2$. Thus $\cX_{3,9}$ satisfies \eqref{eq:polycount} with $m=60$, $N=6$ and $P_9(q):=Q_9(q)+1680\,(9q^2-243q+1684)+665\,280$, whose constant term is $P_9(0)=392\,936\,880$.

\begin{proposition}\label{prop:ff-X3n}
	For $6\leq n\leq9$ the configuration space $\cX_{3,n}$ has genus zero, and
	\[
	\cromb(\cX_{3,n})=\dim\Omega^{2n-8}_{\log}(\cX_{3,n})=
	\begin{cases}
		126 & \text{for } n=6,\\
		7470 & \text{for } n=7,\\
		1\,116\,480 & \text{for } n=8,\\
		392\,936\,880 & \text{for } n=9.
	\end{cases}
	\]
\end{proposition}

\begin{proof}
	The dimension $d=2n-8$ is even, so \zcref{cor:pointcount} gives genus zero and identifies both numbers with the constant term of the counting polynomial.
\end{proof}

\begin{theorem}\label{thm:posarr-3n}
	For $6\leq n\leq9$ the pair $(\Gr_{3,n},Z_{3,n})$ is a positive arrangement.
\end{theorem}

\begin{proof}
	We verify the four conditions of \zcref{def:posarr}. 
	\emph{Condition 1.}  The Grassmannian $\Gr_{3,n}$ is smooth and defined over $\R$, and every configuration of $n$ real points in linearly general position in $\PP^2_\R$ gives a real point of $\Groo_{3,n}$.
	
	\emph{Condition 2.}  The irreducible components of $Z_{3,n}$ are the Schubert divisors $Z_I=\{p_I=0\}\cap\Gr_{3,n}$, $I\in\binom{[n]}{3}$, each defined over $\Q$ and irreducible.  They have a smooth real point: Choose a real configuration in which the three points indexed by $I$ are distinct and collinear and no other three points are collinear. 
	Such a configuration lies in the smooth locus of $Z_I$ and in no other component. 
	
	\emph{Condition 3.}  $\Gr_{3,n}$ is smooth, so $\Sing(\Gr_{3,n})=\varnothing\subseteq Z_{3,n}$.
	
	\emph{Condition 4.}  Genus zero is \zcref{prop:ff-X3n} since $g(\Gr_{3,n},Z_{3,n}) = g(\cX_{3,n})$ by \zcref{cor:genus_and_cr_descent}.
\end{proof}

\begin{remark}\label{rem:bergvall}
	For $n\le 7$ genus zero also follows from known purity results. Recall that a smooth complex variety $U$ is \emph{minimally pure} if $H^k(U)$ is pure of Hodge--Tate type $(k,k)$ for every $k$. Bergvall proves that $\cX_{3,7}$ is minimally pure \cite[Thm.~1.1]{Bergvall}. For $n\le 6$ this is due to Das and O'Connor \cite[Thm.~1.1]{DasOConnor}, and for $n=5$ independently to Bergvall and Gounelas \cite[Thm.~1.3]{BergvallGounelas}. Minimal purity gives genus zero at once: $H_d(X,Y)\cong H^d(\cX_{3,n})(d)$ is then of type $(0,0)$, so $h^{-p,0}(H_d(X,Y))=0$ for $p>0$, and $\cromb(\cX_{3,n})=\dim H^d(\cX_{3,n})$ is the top Betti number. For $n=7$ the Betti numbers can be read off from Bergvall’s point count: the row of the identity
	in \cite[Table~1]{Bergvall} is the polynomial $\#\cX_{3,7}(\F_q)$ above and we obtain $\dim H^6(\cX_{3,7})=7470$, its constant term, in
	agreement with \zcref{prop:ff-X3n}.
	
	For eight points in linearly general position minimal purity is not known. Bergvall's argument applies to eight points in general position, where in addition no six points lie on a conic, and he explains why it does not extend to $\cX_{3,8}$ \cite[\S 4]{Bergvall}.
\end{remark}

The term $840\,b(q)$ in the count of eight points comes from the M\"obius--Kantor configuration, which is realizable over $\C$ but not over $\R$, and in \zcref{ss:38-real} it reappears in the corners of a compactification of $\cX_{3,8}$ that no region reaches. 

\section{The real combinatorial rank}\label{sec:realrank}

By \zcref{thm:posarr-3n}, for $6\le n\le 9$ the pair $(\Gr_{3,n},Z_{3,n})$ is a positive arrangement, so its regions have canonical forms. In this section we compare the span of these forms with the space of all logarithmic top forms. We first reduce the question to the configuration space $\cX_{3,n}$ and show $\cromb_\R(\Gr_{3,n},Z_{3,n}) = \cromb_\R(\cX_{3,n})$. We then consider a tropical compactification of $\cX_{3,n}$ where a canonical form is determined by its iterated residues at the corners and a region only reaches the real corners. For seven points every corner is real by work of Corey \cite{Corey}. By the work of Antolini and Early on chirotropical Grassmannians \cite{Early:2024asu} the iterated residues of the regions are modulo $2$ incidences of chirotropical Dressians with cones of the Dressian. A rank computation over $\F_2$ with the fans of Antolini and Early \cite{AE24data} then shows $\cromb_\R=\cromb=7470$. For eight points, the $840$ initial degenerations with two components found by Corey and Luber \cite{CoreyLuber} give corners that are pairs of complex conjugate points, and these force $\cromb_\R<\cromb$.
 
\subsection{From the Grassmannian to the configuration space}\label{ss:descent}

Let $(X,Y)$ be a positive arrangement of dimension $d$ with $U=X\setminus Y$. The \emph{real combinatorial rank} $\cromb_\R(X,Y)$ is the dimension of the $\Q$-vector space spanned by the canonical forms $\omega(\sigma)$ of the regions $\sigma$ of $(X,Y)$ (\cite[Def.~6.10]{Sturmfels:2026nzc}). These forms lie in $\Omega^d_{\log}(U)$, which has dimension $\cromb(X,Y)$, so $\cromb_\R(X,Y)\le\cromb(X,Y)$.

The class of a region $\sigma$ does not depend on the compactification of $U$. This follows from the isomorphism of  $H_d(X,Y)$ with the locally finite homology $H^{\mathrm{lf}}_d(U)$ (see \cite[Eq. (8)]{Brown:2025jjg} and \cite[\S12]{Sturmfels:2026nzc} for this argument). For a smooth real variety $U$ of genus zero we therefore write $\cromb_\R(U)$ for the dimension of the span of the canonical forms of the closures of the orientable connected components in any compactification defined over $\R$. 

For six points the real combinatorial rank is known. The $372$ regions of $\cX_{3,6}(\R)$ form four orbits under $S_6$, one for each class of uniform oriented matroids of rank $3$ on six elements up to relabeling and reorientation. Their canonical forms span a space of dimension $126$, hence $\cromb_\R(\cX_{3,6})=\cromb(\cX_{3,6})=126$ (see  \cite[\S12]{Sturmfels:2026nzc},  \cite[Prop.~5.3]{HPP}).

We now pass from the Grassmannian to the configuration space on the level of forms and regions, using the trivialization $\Phi=(\pi,\mathbf t) \colon\Groo_{3,n}\longrightarrow\cX_{3,n}\times\Gm^{n-1}$ of \eqref{eq:trivialisation}. Since $\Phi$ is an isomorphism, the $(d+n-1)$ one-forms
\begin{equation}\label{eq:frame}
	dy_{15},\,dy_{25},\,\dots,\,dy_{1n},\,dy_{2n},\qquad \dlog t_i\ \ (i\in[n]\setminus\{3\})
\end{equation}
form a frame of $\Omega^1_{\Groo_{3,n}}$, where $d=2n-8$. Indeed, the $dy_{ja}$ trivialize $\Omega^1_{\cX_{3,n}}$, and $\dlog t_i=dt_i/t_i$ differs from $dt_i$ by a unit. Let
\[
\theta:=\bigwedge_{i\in[n]\setminus\{3\}}\dlog t_i\ \in\ \Omega^{n-1}(\Groo_{3,n}).
\]
It is logarithmic because the $t_i$ are invertible regular functions on $\Groo_{3,n}$. By \eqref{eq:t-equivariance}, an element $\tau$ of $T$ pulls $t_i$ back to $(\tau_i/\tau_3)\,t_i$, so $\tau^*\dlog t_i=\dlog t_i$ and $\theta$ is $T$-invariant. Since $\pi$ is $T$-invariant, so is $\pi^*\omega\wedge\theta$ for every form $\omega$ on $\cX_{3,n}$.

\begin{proposition}\label{prop:descent}
	Let $6\le n\le 9$. The map
	\[
	\Lambda\colon\Omega^{d}_{\log}\bigl(\cX_{3,n}\bigr)\longrightarrow\Omega^{d+n-1}_{\log}\bigl(\Groo_{3,n}\bigr),
	\qquad
	\omega \mapsto \pi^*\omega \wedge \theta,
	\]
	is an isomorphism of $\C$-vector spaces. In particular every logarithmic top form on $\Groo_{3,n}$ is $T$-invariant.
\end{proposition}

\begin{proof}
	Logarithmic forms pull back along morphisms \cite[\S1.5.1]{Brown:2025jjg} and are closed under wedge products, so $\Lambda$ is well defined. Write $\omega=f\,dy$, where $dy$ is the wedge product of the $d$ one-forms $dy_{ja}$. Then $\Lambda(\omega)=\pi^*f\cdot dy\wedge\theta$, and $dy\wedge\theta$ vanishes nowhere because \eqref{eq:frame} is a frame. So $\Lambda(\omega)=0$ forces $\pi^*f=0$, hence $f=0$ since $\pi$ is surjective, and $\Lambda$ is injective. Both spaces have dimension equal to the combinatorial rank plus the genus. The genera vanish and the combinatorial ranks agree by \zcref{cor:genus_and_cr_descent,prop:ff-X3n}, so $\Lambda$ is an isomorphism. The last statement holds since every form in the image of $\Lambda$ is $T$-invariant.
\end{proof}

The same product structure governs the real points. The isomorphism $\Phi$ identifies $\Groo_{3,n}(\R)$ with $\cX_{3,n}(\R) \times (\R^\times)^{n-1}$. The connected components on the right are the products $\sigma^\circ \times Q_\varepsilon$ of a connected component $\sigma^\circ$ of $\cX_{3,n}(\R)$ and an open orthant $Q_{\varepsilon} = \{\varepsilon_i t_i > 0 \}, \varepsilon \in \{\pm 1\}^{n-1}$. 
So the regions of $(\Gr_{3,n},Z_{3,n})$ are indexed by pairs $(\sigma,\varepsilon)$ of a
region $\sigma$ of $\cX_{3,n}$ and a sign vector: $\rho_{\sigma,\varepsilon}$ is the closure in $\Gr_{3,n}(\R)$ of
$\Phi^{-1}(\sigma^\circ\times Q_\varepsilon)$. In particular $(\Gr_{3,n}, Z_{3,n})$ has $2^{n-1}$ times as many regions as $\cX_{3,n}(\R)$.
Their canonical forms are obtained by the same wedge product
with $\theta$.

\begin{proposition}\label{prop:descent-real}
	Let $6\le n\le 9$. For every region $\sigma$ of $\cX_{3,n}$ and every $\varepsilon\in\{\pm1\}^{n-1}$,
	\[
	\omega(\rho_{\sigma,\varepsilon})=\pm\,\Lambda\bigl(\omega(\sigma)\bigr).
	\]
	Hence $\Lambda$ maps the span of the canonical forms of the regions of $\cX_{3,n}$ onto the span of the canonical forms of the regions of $(\Gr_{3,n},Z_{3,n})$, and $\cromb_\R(\Gr_{3,n},Z_{3,n})=\cromb_\R(\cX_{3,n})$.
\end{proposition}

\begin{proof}
	Let $\overline X$ be a compactification of $\cX_{3,n}$ defined over $\R$. We compute the canonical form of $\rho_{\sigma,\varepsilon}$ in the compactification $\overline X\times(\PP^1)^{n-1}$ of $\cX_{3,n}\times\Gm^{n-1}\cong\Groo_{3,n}$. The boundary is the union of $(\overline X\setminus\cX_{3,n})\times(\PP^1)^{n-1}$ and of the loci $\{t_i\in\{0,\infty\}\}$. The closure of $\sigma^\circ\times Q_\varepsilon$ is $\sigma\times\overline{Q_\varepsilon}$ with
	\[
	\overline{Q_\varepsilon}=\prod_{i\ne3}I_{\varepsilon_i}\ \subseteq\ (\PP^1_\R)^{n-1},\qquad
	I_+=[0,\infty],\quad I_-=[-\infty,0].
	\]
	The pair $(\overline X,\overline X\setminus\cX_{3,n})$ has genus zero by \zcref{prop:ff-X3n}, and $(\PP^1,\{0,\infty\})$ has genus zero \cite[\S5.2]{Brown:2025jjg}. Each closed half-line $I_\pm$ is a path from $0$ to $\infty$ in $\PP^1(\R)$, so its canonical form is $\pm\,dt/t$ \cite[\S5.2]{Brown:2025jjg}. The product formula \cite[Prop.~2.12]{Brown:2025jjg} gives $\omega(\overline{Q_\varepsilon})=\pm\,\theta$ in the coordinates $t_i$, and applied once more,
	\[
	\omega(\rho_{\sigma,\varepsilon})=\pm\,\mathrm{pr}_1^*\omega(\sigma)\wedge\mathrm{pr}_2^*\omega(\overline{Q_\varepsilon})
	=\pm\,\pi^*\omega(\sigma)\wedge\theta=\pm\,\Lambda\bigl(\omega(\sigma)\bigr).
	\]
	Every region of $(\Gr_{3,n},Z_{3,n})$ is some $\rho_{\sigma,\varepsilon}$ and $\Lambda$ is injective, so $\Lambda$ maps the span of the $\omega(\sigma)$ onto the span of the $\omega(\rho_{\sigma,\varepsilon})$, and the two real combinatorial ranks agree.
\end{proof}

\subsection{Tropical compactifications of $\cX_{3,7}$ and $\cX_{3,8}$}\label{ss:corners}

We now describe a tropical compactification of $\cX_{3,n}$ for $n=7,8$ with simple normal crossing boundary. Brown and Dupont determine a canonical form by its iterated residues at the corners of such a compactification. We describe
the corners by the maximal cones of the tropical variety, and we describe how a region of $\cX_{3,n}(\R)$ meets a corner.

Let $N=\Z^{\binom{[n]}{3}}$ with basis $(\mathbf e_I)_I$ and $\mathbf 1=\sum_I\mathbf e_I$, so that $N/\Z\mathbf 1$ is the cocharacter lattice of the Plücker torus $T_P = \G_m^{\binom{n}{3}} / \G_m$. Let $\varphi\colon\R^n\to N_\R$, $\varphi(a)_I=\sum_{i\in I}a_i$, and let $L=\varphi(\R^n)$. The map $\varphi$ is injective and $\mathbf 1=\varphi(\tfrac13,\dots,\tfrac13)\in L$. The torus $T = \G_m^n/ \G_m$ acts freely on $\Groo_{3,n}$, and its cocharacter lattice is $(N\cap L)/\Z\mathbf 1$.  The torus $T_X:=T_P/T$ of $\cX_{3,n}$ has cocharacter lattice $N_X=N/(N\cap L)$ and character lattice
\[
M_X=(N\cap L)^\perp=\Bigl\{m\in\Z^{\binom{[n]}{3}}:\ \textstyle\sum_{I\ni i}m_I=0\text{ for } i=1,\dots,n\Bigr\}.
\]
A character $m\in M_X$ is a Laurent monomial $p^m=\prod_Ip_I^{m_I}$ in which every index occurs as often in the numerator as in the denominator, such as the cross-ratios $y_{ja}$. The configuration space is a closed subvariety $\cX_{3,n}\subseteq T_X$. Let $\pi\colon N_\R\to(N_X)_\R=N_\R/L$ be the projection. The quotient map $\Groo_{3,n}\to\cX_{3,n}$ is the restriction of the quotient map of tori, so $\Trop\cX_{3,n}=\pi(\Trop\Groo_{3,n})$, a pure $d$-dimensional fan without lineality space. For the initial degenerations we have
\begin{equation}\label{eq:in-split}
	\inw\Groo_{3,n}\ \cong\ \operatorname{in}_{\pi(w)}\cX_{3,n}\times T\qquad(w\in\Trop\Groo_{3,n})
\end{equation}
(see \cite[\S 7]{CoreyLuber}). So results on initial degenerations of $\Groo_{3,n}$ are results on $\cX_{3,n}$. 

For $n=7$ every initial degeneration $\inw\cX_{3,7}$ is smooth and irreducible \cite[Thm.~1.2]{Corey}. For $n=8$ every initial degeneration $\inw\cX_{3,8}$ is smooth \cite[Thm.~1.1]{CoreyLuber}, and it is irreducible except over the $840$ Möbius--Kantor cones of \zcref{ss:38-real}. In particular $\cX_{3,n}$ is \emph{sch\"on} for $n=7,8$. 

Let $\Sigma$ be a unimodular fan in $(N_X)_\R$ with support $\Trop\cX_{3,n}$. Let $X$ be the closure of $\cX_{3,n}$ in the toric variety $X_{\Sigma}$, the \emph{tropical compactification} of $\cX_{3,n}$ with respect to $\Sigma$. Let $Y=X\setminus\cX_{3,n}$ denote the boundary divisor. For a cone $\tau\in \Sigma$ let $O_\tau$ be its torus orbit and $X_\tau=X\cap O_\tau$ its stratum. Since $|\Sigma|=\Trop\cX_{3,n}$, the variety $X$ is compact \cite[Prop.~2.3]{Tevelev}. Since $\cX_{3,n}$ is sch\"on and the fan $\Sigma$ is unimodular, the tropical compactification $X$ is smooth and $Y$ is a simple normal crossing divisor \cite[Thm.~1.4]{Tevelev}. Further, every stratum $X_\tau$ is nonempty and smooth of codimension $\dim\tau$ \cite[Thm.~2.4, Lemma~2.7]{Hacking}. For $w\in\relint\tau$ there is an isomorphism
\begin{equation}\label{eq:in-stratum}
	\inw\cX_{3,n}\ \cong\ X_\tau\times\Gm^{\dim\tau}
\end{equation}
\cite[Lemma~3.6]{helmkatz}. The corners of $(X,Y)$ in the sense of Brown and Dupont \cite[Def.~2.16]{Brown:2025jjg} are the points of $X$ on $d$ local branches of $Y$. 

\begin{lemma}\label{lem:corners}
	Let $n = 7,8$ and let $(X,Y)$ be the tropical compactification of $\cX_{3,n}$ with respect to $\Sigma$.
	\begin{enumerate}
		\item A point of $X_\tau$ lies on exactly $\dim\tau$ local branches of $Y$, one on each $Y_\rho:=X\cap\overline{O_\rho}$ for the rays $\rho$ of $\tau$. If $\inw \cX_{3,n}$ is irreducible for $w$ in the relative interior of a ray $\rho$, then $Y_\rho$ is an irreducible component of $Y$.
		\item The corners of $Y$ are the points of the zero-dimensional strata $X_C$, where $C$ is a maximal cone of $\Sigma$. For $w\in\relint C$, the number of corners in $X_C$ is the number of irreducible components of $\inw\cX_{3,n}$.
		\item Let $C \in \Sigma$ be a maximal cone. If $\inw\cX_{3,n}$ is irreducible for $w\in\relint C$, then $X_C$ is a single point, defined over $\Q$. In particular, this corner is real.
	\end{enumerate}
\end{lemma}
 
\begin{proof}
	Since $\tau\in\Sigma$ is unimodular, the affine chart of $\tau$ is $U_\tau\cong\A^k\times\Gm^{r-k}$, where $k=\dim\tau$, with coordinates $z_1,\dots,z_r$. Near a point of $X_\tau=X\cap O_\tau$ the divisor $Y$ is $X\cap\{z_1\cdots z_k=0\}$. Since $Y$ is simple normal crossing, the $k$ hypersurfaces $X\cap\{z_j=0\}$ are its local branches. By \eqref{eq:in-stratum}, the divisor $Y_\rho$ is irreducible if $\inw\cX_{3,n}$ is, which proves (1). A point of a stratum $X_C$ is a corner if and only if it lies on $d=\dim C$ branches of $Y$, that is, if $C$ is a maximal cone. For such a cone $C$ the stratum $X_C$ is smooth of dimension zero, and by \eqref{eq:in-stratum} its complex points correspond to the irreducible components of $\inw\cX_{3,n}$. This proves (2). If there is only one irreducible component, then $X_C$ is a single point defined over $\Q$, so it is real.
\end{proof}

Let $C$ be a maximal cone of $\Sigma$ and let $x \in X_C$ be a real corner. In the chart $U_C \cong \mathbb{A}^d \times \G_m^{r-d}$ of the proof, $x$ lies on the $d$ branches $\{z_j = 0\}$ for $j \leq d$ of the simple normal crossing divisor $Y$. So $z_1,\ldots,z_d$ are local coordinates of $X$ at $x$. Thus for a sufficiently small open neighborhood $V \subseteq X(\R)$ of $x$ we have that $V \setminus Y(\R)$ is the disjoint union of the $2^d$ open orthants
\begin{equation}\label{eq:orthants}
	V_\varepsilon\ =\ \{\,y\in V:\ \varepsilon_1z_1(y)>0,\ \dots,\ \varepsilon_dz_d(y)>0\,\},\qquad \varepsilon\in\{\pm1\}^d.
\end{equation}
For a real configuration $y$ we write $\chi(y)$ for its chirotope, the sign vector of the Plücker coordinates of a $3\times n$ matrix representing $y$, up to reorientation. So $S_\chi=\{y\in\cX_{3,n}(\R):\chi(y)=\chi\}$ is the oriented matroid stratum of $\chi$. The coordinate $z_j$ is a character in $M_X$, a Laurent monomial in the Plücker coordinates in which every index occurs equally often in the numerator and in the denominator, so its sign at $y$ only depends on the class $\chi(y)$. On $V_\varepsilon$ the coordinate $z_j$ has the sign $\varepsilon_j$ for $j\le d$. Hence the $2^d$ orthants at $x$ have pairwise different chirotopes, so they lie in $2^d$ different oriented matroid strata.

Now let $Y=Y_1\cup\dots\cup Y_N$ be the decomposition into irreducible components and let $\mathcal C$ be the set of corners. For $J\subseteq[N]$ let $Y_J=\bigcap_{j\in J}Y_j$, and for $x\in Y_J$ let $Y_J(x)$ be its connected component through $x$. A corner $c$ lies on $d$ components $Y_{i_1},\dots,Y_{i_d}$. Brown and Dupont define its iterated boundary $\partial_c\colon H_d(X,Y)\to H_0(\{c\})\cong\Q$ as the composition of the partial boundary maps along $Y_{i_d}$, $Y_{i_{d-1}}\cap Y_{i_d}$, \dots, and the iterated residue $\Res_c\colon\Omega^d_{\log}(\cX_{3,n})\to\C$ in the same way \cite[(10), (36), (37)]{Brown:2025jjg}. For $n=7,8$ the genus is zero by \zcref{prop:ff-X3n}, so the iterated residues give an injective map
\[
\Res=\textstyle\bigoplus_{c\in\mathcal C}\Res_c\colon\ \Omega^d_{\log}(\cX_{3,n})\hookrightarrow\C^{\mathcal C}.
\]
The canonical form $\omega(\sigma)$ of a class $\sigma\in H_d(X,Y)$ is the unique logarithmic form with $\Res_c\,\omega(\sigma)=\partial_c\sigma$ for all corners $c$ \cite[Prop.~2.17, Cor.~2.18]{Brown:2025jjg}. The image of $\Res$ is described by the dual complex $\Delta(Y)$, which has one $(k-1)$-simplex for every connected component of every $Y_J$ with $|J|=k$. It has no simplices of dimension $d$, and its $(d-1)$-simplices are the corners. Brown and Dupont identify $\Omega^d_{\log}(\cX_{3,n})$ with the top homology of $\Delta(Y)$, compatibly with the iterated boundaries \cite[Prop.~2.17, Prop.~4.4]{Brown:2025jjg}. Thus $\Res$ is an isomorphism of $\Omega^d_{\log}(\cX_{3,n})$ onto $W\otimes_\Q\C$, where
\begin{equation}\label{eq:W}
	W:=\widetilde H_{d-1}\bigl(\Delta(Y);\Q\bigr)=\ker\Bigl(\partial\colon\Q^{\mathcal C}\to\bigoplus_{|J|=d-1}H_0(Y_J)\Bigr),
	\qquad
	\partial(c)=\sum_{k=1}^d(-1)^{k-1}\bigl[Y_{I\setminus\{i_k\}}(c)\bigr],
\end{equation}
for a corner $c$ on the components indexed by $I=\{i_1<\dots<i_d\}$. In particular $\cromb(\cX_{3,n})=\dim_\Q W$. 

\begin{remark}\label{rem:tropical-rank-formula}\label{rem:no-fvector}
	Suppose that $\inw\cX_{3,n}$ is irreducible whenever $w$ lies in the relative interior of a cone of $\Sigma$ that is not maximal. By \eqref{eq:in-stratum} every positive-dimensional stratum $X_\tau$ is then connected, that is, $\cX_{3,n}$ is \emph{sehr sch\"on} in the sense of \cite{Koefler}. The intersection $Y_J$ is the closure of $X_\tau$ if the rays indexed by $J$ span a cone $\tau$ of dimension less than $d$ and it is $X_C$ if they span a maximal cone $C$. So $\Delta(Y)$ is the link of $\Sigma$ in which each maximal cone $C$ appears as many times as $X_C$ has points. Koefler's tropical formula for the combinatorial rank \cite[Prop.~2.2.12]{Koefler} expresses $\cromb$ through the face numbers of this complex. Since $\Trop\cX_{3,n}$ has no lineality space, it reads
	\[
	\cromb(\cX_{3,n})\ =\ \sum_{k=0}^{d}(-1)^{d-k}f_k ,
	\]
	where $f_0=1$, $f_k$ is the number of $k$-dimensional cones of $\Sigma$ for $0<k<d$, and $f_d$ counts the maximal cones with their numbers of corners, i.e. with their tropical multiplicities. For $n=7$ all initial degenerations are smooth and irreducible \cite[Thm.~1.2]{Corey}, so every maximal cone counts once, and the alternating sum is $(-1)^{d-1}$ times the reduced Euler characteristic of the link of $\Trop\cX_{3,7}$. The Gr\"obner fan structure on $\Trop\cX_{3,7}$ is simplicial with $f$-vector $(1,721,16800,124180,386155,522585,252000)$ \cite[Thm.~2.1]{HJJS}, and its alternating sum is $7470$. This number is also $\dim\widetilde
	H_5(\link\Trop\cX_{3,7})$. This is the last equality in \zcref{thm:main37}.
\end{remark}

Since $\Res$ is injective and the iterated boundaries are rational, the real combinatorial rank is a rank of an integer matrix:
\begin{equation}\label{eq:crR-matrix}
	\cromb_\R(\cX_{3,n})=\dim_\Q\operatorname{span}\bigl\{(\partial_c\sigma)_{c\in\mathcal C}\ :\ \sigma\text{ a region of }\cX_{3,n}\bigr\}.
\end{equation}

We call a corner \emph{real} if it lies in $X(\R)$. A region only sees the real corners, and at a real corner it is a union of closed orthants.
 
\begin{lemma}\label{lem:realcorners}
	Let $c$ be a corner of $Y$ and let $\sigma$ be a region of $\cX_{3,n}$.
	\begin{enumerate}
		\item If $c\notin\sigma$, then $\partial_c\sigma=0$. In particular $\partial_c\sigma=0$ for every corner $c$ that is not real.
		\item Let $c$ be real, and let $V$ and the orthants $V_\varepsilon$ be as in \eqref{eq:orthants}. If $\sigma\cap V$ is the closure in $V$ of a single orthant $V_\varepsilon$, then $\partial_c\sigma=\pm1$.
	\end{enumerate}
\end{lemma}
 
\begin{proof}    
  Let the corner $c$ lie on $Y_{i_1}\cap\dots\cap Y_{i_d}$. The
  first partial boundary map in $\partial_c$ takes the boundary of $\sigma$ and keeps its part
  in $Y_{i_d}$, the next one takes the boundary of that chain and keeps its part in
  $Y_{i_{d-1}}\cap Y_{i_d}$, and so on. So all chains that occur are supported in $\sigma$, and the last one is a 0-chain on $\sigma \cap Y_{i_1} \cap \dots \cap Y_{i_d}$ whose
    coefficient at $c$ is $\partial_c \sigma$. It vanishes if $c \notin \sigma$, and a corner that is not real does not lie in $\sigma\subseteq X(\R)$. This proves (1). In (2) the image of $[\sigma]$ is the image of the closed orthant $\{\varepsilon_1z_1\ge0,\dots,\varepsilon_dz_d\ge0\}$ in the coordinates $z_j$, and the iterated boundary of an orthant is $\pm1$ \cite[proof of Prop.~6.10]{Brown:2025jjg}.

\end{proof}

\subsection{Seven points}\label{ss:real37}

In this subsection $n=7$, so $d=\dim\cX_{3,7}=6$ and $\dim T_X=28$. Let $\mathcal K$ be the set of the reorientation classes of uniform chirotopes of rank $3$ on $[7]$. Up to relabeling there are $11$ of them, and $|\mathcal K|=27\,240$. \zcref{fig:chirotope_table} lists representatives, and \zcref{fig:chirotope_images} shows a configuration realizing each of them. So every class in $\mathcal K$ is realizable, as is every uniform oriented matroid of rank $3$ on at most eight elements \cite[Ch.~8]{BLVSWZ}, and every stratum $S_\chi$, $\chi\in\mathcal K$, is nonempty. Nonempty oriented matroid strata of $\cX_{3,n}$ are connected for $n\le9$ \cite{richter1988kombinatorische}. Hence the strata $S_\chi$, $\chi\in\mathcal K$, are exactly the connected components of $\cX_{3,7}(\R)$, and the regions of $\cX_{3,7}$ are their closures $R_\chi=\overline{S_\chi}$, $\chi\in\mathcal K$.

\begin{figure}
	\[
	\setlength{\arraycolsep}{0.5pt}
	\begin{array}{r|@{\hspace{4pt}}ccccccccccccccccccccccccccccccccccc}
		& \multicolumn{3}{@{\hspace{4pt}}l}{p_{123}} & \multicolumn{29}{c}{\cdots} & \multicolumn{3}{r}{p_{567}} \\ \hline
		\chi_{1} & + & + & + & + & + & + & + & + & + & + & + & + & + & + & + & + & + & + & + & - & + & + & + & + & + & + & + & + & - & - & - & - & - & + & - \\
		\chi_{2} & + & + & + & + & + & + & + & + & + & + & + & + & + & + & + & + & + & + & + & - & + & + & + & + & + & + & + & + & - & - & + & - & - & - & - \\
		\chi_{3} & + & + & + & + & + & + & + & + & + & + & + & + & + & + & + & + & + & + & + & - & + & + & + & + & + & + & + & + & + & - & + & - & - & - & - \\
		\chi_{4} & + & + & + & + & + & + & + & + & + & + & + & + & + & + & + & + & + & + & + & + & + & + & + & + & + & + & - & - & + & + & - & - & + & + & + \\
		\chi_{5} & + & + & + & + & + & + & + & + & + & + & + & + & + & + & + & + & + & + & + & + & + & + & + & + & + & + & - & + & + & + & - & - & + & + & + \\
		\chi_{6} & + & + & + & + & + & + & + & + & + & + & + & + & + & + & + & + & + & + & + & + & + & + & + & + & + & + & + & + & + & - & + & - & - & - & - \\
		\chi_{7} & + & + & + & + & + & + & + & + & + & + & + & + & + & + & + & + & + & + & + & + & + & + & + & + & + & + & + & + & + & - & + & + & - & - & - \\
		\chi_{8} & + & + & + & + & + & + & + & + & + & + & + & + & + & + & + & + & + & + & + & + & + & + & + & + & + & + & + & + & + & - & + & + & + & - & - \\
		\chi_{9} & + & + & + & + & + & + & + & + & + & + & + & + & + & + & + & + & + & + & + & + & + & + & + & + & + & + & + & + & + & + & + & + & + & - & - \\
		\chi_{10} & + & + & + & + & + & + & + & + & + & + & + & + & + & + & + & + & + & + & + & + & + & + & + & + & + & + & + & + & + & + & + & + & + & + & - \\
		\chi_{11} & + & + & + & + & + & + & + & + & + & + & + & + & + & + & + & + & + & + & + & + & + & + & + & + & + & + & + & + & + & + & + & + & + & + & +
	\end{array}
	\]
	\caption{The 11 chirotope representatives of equivalence classes of uniform oriented matroids of cardinality 7 and rank 3. The rows and columns are ordered following the conventions of Finschi's catalog \cite{FinschiCatalog,FinschiFukuda}, with the columns in colexicographic order ($p_{123},p_{124},p_{134},p_{234},\dots$).}\label{fig:chirotope_table}
\end{figure}

\begin{figure}
	\begingroup
\centering
\makebox[0.32\linewidth]{\begin{tikzpicture}[scale=0.7]
  \begin{scope}\clip (-3,-3) rectangle (3,3);
  \coordinate (p1) at (0,2.2);
  \coordinate (p2) at (-2.092,0.6798);
  \coordinate (p3) at (-1.293,-1.78);
  \coordinate (p4) at (1.293,-1.78);
  \coordinate (p5) at (1.041,-0.4266);
  \coordinate (p6) at (2.092,0.6798);
  \coordinate (p7) at (1.531,-0.381);
  \draw[black!45, line width=0.6pt] (p2) -- (p3) -- (p4) -- (p6) -- (p1) -- cycle;
  \draw[red!70!black, dashed, line width=0.5pt, shorten <=-12cm, shorten >=-12cm] (p1) -- (p4);
  \draw[blue!70!black, dashed, line width=0.5pt, shorten <=-12cm, shorten >=-12cm] (p2) -- (p5);
  \draw[blue!70!black, dashed, line width=0.5pt, shorten <=-12cm, shorten >=-12cm] (p3) -- (p5);
  \draw[red!70!black, dashed, line width=0.5pt, shorten <=-12cm, shorten >=-12cm] (p3) -- (p6);
  \fill[white] (p1) circle (3.3pt); \fill[black] (p1) circle (1.9pt);
  \node[font=\footnotesize, text=black, fill=white, fill opacity=0.8, text opacity=1, inner sep=0.6pt] at ([shift={(75.0:0.38)}]p1) {1};
  \fill[white] (p2) circle (3.3pt); \fill[black] (p2) circle (1.9pt);
  \node[font=\footnotesize, text=black, fill=white, fill opacity=0.8, text opacity=1, inner sep=0.6pt] at ([shift={(197.0:0.38)}]p2) {2};
  \fill[white] (p3) circle (3.3pt); \fill[black] (p3) circle (1.9pt);
  \node[font=\footnotesize, text=black, fill=white, fill opacity=0.8, text opacity=1, inner sep=0.6pt] at ([shift={(-111.0:0.38)}]p3) {3};
  \fill[white] (p4) circle (3.3pt); \fill[black] (p4) circle (1.9pt);
  \node[font=\footnotesize, text=black, fill=white, fill opacity=0.8, text opacity=1, inner sep=0.6pt] at ([shift={(-39.0:0.38)}]p4) {4};
  \fill[white] (p5) circle (3.3pt); \fill[red!70!black] (p5) circle (1.9pt);
  \node[font=\footnotesize, text=red!70!black, fill=white, fill opacity=0.8, text opacity=1, inner sep=0.6pt] at ([shift={(300:0.38)}]p5) {5};
  \fill[white] (p6) circle (3.3pt); \fill[black] (p6) circle (1.9pt);
  \node[font=\footnotesize, text=black, fill=white, fill opacity=0.8, text opacity=1, inner sep=0.6pt] at ([shift={(3.0:0.38)}]p6) {6};
  \fill[white] (p7) circle (3.3pt); \fill[blue!70!black] (p7) circle (1.9pt);
  \node[font=\footnotesize, text=blue!70!black, fill=white, fill opacity=0.8, text opacity=1, inner sep=0.6pt] at ([shift={(345:0.38)}]p7) {7};
  \end{scope}
  \draw[black!25, line width=0.3pt] (-3,-3) rectangle (3,3);
  \node[anchor=north] at (0,-3.1) {\small $\chi_{1}$};
\end{tikzpicture}}%
\makebox[0.32\linewidth]{\begin{tikzpicture}[scale=0.7]
  \begin{scope}\clip (-3,-3) rectangle (3,3);
  \coordinate (p1) at (0,2.2);
  \coordinate (p2) at (-2.092,0.6798);
  \coordinate (p3) at (-1.293,-1.78);
  \coordinate (p4) at (1.293,-1.78);
  \coordinate (p5) at (0.8594,0.385);
  \coordinate (p6) at (1.056,1.037);
  \coordinate (p7) at (2.092,0.6798);
  \draw[black!45, line width=0.6pt] (p2) -- (p3) -- (p4) -- (p7) -- (p1) -- cycle;
  \draw[red!70!black, dashed, line width=0.5pt, shorten <=-12cm, shorten >=-12cm] (p1) -- (p4);
  \draw[red!70!black, dashed, line width=0.5pt, shorten <=-12cm, shorten >=-12cm] (p2) -- (p7);
  \draw[blue!70!black, dashed, line width=0.5pt, shorten <=-12cm, shorten >=-12cm] (p3) -- (p5);
  \draw[red!70!black, dashed, line width=0.5pt, shorten <=-12cm, shorten >=-12cm] (p3) -- (p7);
  \draw[blue!70!black, dashed, line width=0.5pt, shorten <=-12cm, shorten >=-12cm] (p4) -- (p5);
  \fill[white] (p1) circle (3.3pt); \fill[black] (p1) circle (1.9pt);
  \node[font=\footnotesize, text=black, fill=white, fill opacity=0.8, text opacity=1, inner sep=0.6pt] at ([shift={(75.0:0.38)}]p1) {1};
  \fill[white] (p2) circle (3.3pt); \fill[black] (p2) circle (1.9pt);
  \node[font=\footnotesize, text=black, fill=white, fill opacity=0.8, text opacity=1, inner sep=0.6pt] at ([shift={(147.0:0.38)}]p2) {2};
  \fill[white] (p3) circle (3.3pt); \fill[black] (p3) circle (1.9pt);
  \node[font=\footnotesize, text=black, fill=white, fill opacity=0.8, text opacity=1, inner sep=0.6pt] at ([shift={(-101.0:0.38)}]p3) {3};
  \fill[white] (p4) circle (3.3pt); \fill[black] (p4) circle (1.9pt);
  \node[font=\footnotesize, text=black, fill=white, fill opacity=0.8, text opacity=1, inner sep=0.6pt] at ([shift={(-39.0:0.38)}]p4) {4};
  \fill[white] (p5) circle (3.3pt); \fill[red!70!black] (p5) circle (1.9pt);
  \node[font=\footnotesize, text=red!70!black, fill=white, fill opacity=0.8, text opacity=1, inner sep=0.6pt] at ([shift={(0:0.38)}]p5) {5};
  \fill[white] (p6) circle (3.3pt); \fill[blue!70!black] (p6) circle (1.9pt);
  \node[font=\footnotesize, text=blue!70!black, fill=white, fill opacity=0.8, text opacity=1, inner sep=0.6pt] at ([shift={(120:0.38)}]p6) {6};
  \fill[white] (p7) circle (3.3pt); \fill[black] (p7) circle (1.9pt);
  \node[font=\footnotesize, text=black, fill=white, fill opacity=0.8, text opacity=1, inner sep=0.6pt] at ([shift={(90:0.38)}]p7) {7};
  \end{scope}
  \draw[black!25, line width=0.3pt] (-3,-3) rectangle (3,3);
  \node[anchor=north] at (0,-3.1) {\small $\chi_{2}$};
\end{tikzpicture}}%
\makebox[0.32\linewidth]{\begin{tikzpicture}[scale=0.7]
  \begin{scope}\clip (-3,-3) rectangle (3,3);
  \coordinate (p1) at (0,2.2);
  \coordinate (p2) at (-2.092,0.6798);
  \coordinate (p3) at (-1.293,-1.78);
  \coordinate (p4) at (1.293,-1.78);
  \coordinate (p5) at (1.006,-0.3268);
  \coordinate (p6) at (1.078,1.01);
  \coordinate (p7) at (2.092,0.6798);
  \draw[black!45, line width=0.6pt] (p2) -- (p3) -- (p4) -- (p7) -- (p1) -- cycle;
  \draw[red!70!black, dashed, line width=0.5pt, shorten <=-12cm, shorten >=-12cm] (p1) -- (p4);
  \draw[blue!70!black, dashed, line width=0.5pt, shorten <=-12cm, shorten >=-12cm] (p2) -- (p7);
  \draw[red!70!black, dashed, line width=0.5pt, shorten <=-12cm, shorten >=-12cm] (p3) -- (p7);
  \draw[blue!70!black, dashed, line width=0.5pt, shorten <=-12cm, shorten >=-12cm] (p4) -- (p5);
  \fill[white] (p1) circle (3.3pt); \fill[black] (p1) circle (1.9pt);
  \node[font=\footnotesize, text=black, fill=white, fill opacity=0.8, text opacity=1, inner sep=0.6pt] at ([shift={(75.0:0.38)}]p1) {1};
  \fill[white] (p2) circle (3.3pt); \fill[black] (p2) circle (1.9pt);
  \node[font=\footnotesize, text=black, fill=white, fill opacity=0.8, text opacity=1, inner sep=0.6pt] at ([shift={(147.0:0.38)}]p2) {2};
  \fill[white] (p3) circle (3.3pt); \fill[black] (p3) circle (1.9pt);
  \node[font=\footnotesize, text=black, fill=white, fill opacity=0.8, text opacity=1, inner sep=0.6pt] at ([shift={(-111.0:0.38)}]p3) {3};
  \fill[white] (p4) circle (3.3pt); \fill[black] (p4) circle (1.9pt);
  \node[font=\footnotesize, text=black, fill=white, fill opacity=0.8, text opacity=1, inner sep=0.6pt] at ([shift={(-39.0:0.38)}]p4) {4};
  \fill[white] (p5) circle (3.3pt); \fill[red!70!black] (p5) circle (1.9pt);
  \node[font=\footnotesize, text=red!70!black, fill=white, fill opacity=0.8, text opacity=1, inner sep=0.6pt] at ([shift={(0:0.38)}]p5) {5};
  \fill[white] (p6) circle (3.3pt); \fill[blue!70!black] (p6) circle (1.9pt);
  \node[font=\footnotesize, text=blue!70!black, fill=white, fill opacity=0.8, text opacity=1, inner sep=0.6pt] at ([shift={(345:0.38)}]p6) {6};
  \fill[white] (p7) circle (3.3pt); \fill[black] (p7) circle (1.9pt);
  \node[font=\footnotesize, text=black, fill=white, fill opacity=0.8, text opacity=1, inner sep=0.6pt] at ([shift={(90:0.38)}]p7) {7};
  \end{scope}
  \draw[black!25, line width=0.3pt] (-3,-3) rectangle (3,3);
  \node[anchor=north] at (0,-3.1) {\small $\chi_{3}$};
\end{tikzpicture}}

\medskip
\makebox[0.32\linewidth]{\begin{tikzpicture}[scale=0.7]
  \begin{scope}\clip (-3,-3) rectangle (3,3);
  \coordinate (p1) at (-0.1408,2.195);
  \coordinate (p2) at (-1.831,1.22);
  \coordinate (p3) at (-1.831,-1.22);
  \coordinate (p4) at (-0.1408,-2.195);
  \coordinate (p5) at (1.972,-0.9758);
  \coordinate (p6) at (1.972,0.9758);
  \coordinate (p7) at (0,0);
  \draw[black!45, line width=0.6pt] (p3) -- (p4) -- (p5) -- (p6) -- (p1) -- (p2) -- cycle;
  \draw[red!70!black, dashed, line width=0.5pt, shorten <=-12cm, shorten >=-12cm] (p1) -- (p4);
  \draw[red!70!black, dashed, line width=0.5pt, shorten <=-12cm, shorten >=-12cm] (p2) -- (p5);
  \draw[red!70!black, dashed, line width=0.5pt, shorten <=-12cm, shorten >=-12cm] (p3) -- (p6);
  \fill[white] (p1) circle (3.3pt); \fill[black] (p1) circle (1.9pt);
  \node[font=\footnotesize, text=black, fill=white, fill opacity=0.8, text opacity=1, inner sep=0.6pt] at ([shift={(128.7:0.38)}]p1) {1};
  \fill[white] (p2) circle (3.3pt); \fill[black] (p2) circle (1.9pt);
  \node[font=\footnotesize, text=black, fill=white, fill opacity=0.8, text opacity=1, inner sep=0.6pt] at ([shift={(111.3:0.38)}]p2) {2};
  \fill[white] (p3) circle (3.3pt); \fill[black] (p3) circle (1.9pt);
  \node[font=\footnotesize, text=black, fill=white, fill opacity=0.8, text opacity=1, inner sep=0.6pt] at ([shift={(-111.3:0.38)}]p3) {3};
  \fill[white] (p4) circle (3.3pt); \fill[black] (p4) circle (1.9pt);
  \node[font=\footnotesize, text=black, fill=white, fill opacity=0.8, text opacity=1, inner sep=0.6pt] at ([shift={(-128.7:0.38)}]p4) {4};
  \fill[white] (p5) circle (3.3pt); \fill[black] (p5) circle (1.9pt);
  \node[font=\footnotesize, text=black, fill=white, fill opacity=0.8, text opacity=1, inner sep=0.6pt] at ([shift={(8.7:0.38)}]p5) {5};
  \fill[white] (p6) circle (3.3pt); \fill[black] (p6) circle (1.9pt);
  \node[font=\footnotesize, text=black, fill=white, fill opacity=0.8, text opacity=1, inner sep=0.6pt] at ([shift={(-8.7:0.38)}]p6) {6};
  \fill[white] (p7) circle (3.3pt); \fill[red!70!black] (p7) circle (1.9pt);
  \node[font=\footnotesize, text=red!70!black, fill=white, fill opacity=0.8, text opacity=1, inner sep=0.6pt] at ([shift={(180:0.38)}]p7) {7};
  \end{scope}
  \draw[black!25, line width=0.3pt] (-3,-3) rectangle (3,3);
  \node[anchor=north] at (0,-3.1) {\small $\chi_{4}$};
\end{tikzpicture}}%
\makebox[0.32\linewidth]{\begin{tikzpicture}[scale=0.7]
  \begin{scope}\clip (-3,-3) rectangle (3,3);
  \coordinate (p1) at (0,2.2);
  \coordinate (p2) at (-1.905,1.1);
  \coordinate (p3) at (-1.905,-1.1);
  \coordinate (p4) at (0,-2.2);
  \coordinate (p5) at (1.905,-1.1);
  \coordinate (p6) at (1.905,1.1);
  \coordinate (p7) at (0.4026,0.6974);
  \draw[black!45, line width=0.6pt] (p3) -- (p4) -- (p5) -- (p6) -- (p1) -- (p2) -- cycle;
  \draw[red!70!black, dashed, line width=0.5pt, shorten <=-12cm, shorten >=-12cm] (p1) -- (p4);
  \draw[red!70!black, dashed, line width=0.5pt, shorten <=-12cm, shorten >=-12cm] (p1) -- (p5);
  \draw[red!70!black, dashed, line width=0.5pt, shorten <=-12cm, shorten >=-12cm] (p2) -- (p6);
  \draw[red!70!black, dashed, line width=0.5pt, shorten <=-12cm, shorten >=-12cm] (p3) -- (p6);
  \fill[white] (p1) circle (3.3pt); \fill[black] (p1) circle (1.9pt);
  \node[font=\footnotesize, text=black, fill=white, fill opacity=0.8, text opacity=1, inner sep=0.6pt] at ([shift={(55.0:0.38)}]p1) {1};
  \fill[white] (p2) circle (3.3pt); \fill[black] (p2) circle (1.9pt);
  \node[font=\footnotesize, text=black, fill=white, fill opacity=0.8, text opacity=1, inner sep=0.6pt] at ([shift={(142.0:0.38)}]p2) {2};
  \fill[white] (p3) circle (3.3pt); \fill[black] (p3) circle (1.9pt);
  \node[font=\footnotesize, text=black, fill=white, fill opacity=0.8, text opacity=1, inner sep=0.6pt] at ([shift={(-115.0:0.38)}]p3) {3};
  \fill[white] (p4) circle (3.3pt); \fill[black] (p4) circle (1.9pt);
  \node[font=\footnotesize, text=black, fill=white, fill opacity=0.8, text opacity=1, inner sep=0.6pt] at ([shift={(-55.0:0.38)}]p4) {4};
  \fill[white] (p5) circle (3.3pt); \fill[black] (p5) circle (1.9pt);
  \node[font=\footnotesize, text=black, fill=white, fill opacity=0.8, text opacity=1, inner sep=0.6pt] at ([shift={(-22.0:0.38)}]p5) {5};
  \fill[white] (p6) circle (3.3pt); \fill[black] (p6) circle (1.9pt);
  \node[font=\footnotesize, text=black, fill=white, fill opacity=0.8, text opacity=1, inner sep=0.6pt] at ([shift={(65.0:0.38)}]p6) {6};
  \fill[white] (p7) circle (3.3pt); \fill[red!70!black] (p7) circle (1.9pt);
  \node[font=\footnotesize, text=red!70!black, fill=white, fill opacity=0.8, text opacity=1, inner sep=0.6pt] at ([shift={(240:0.38)}]p7) {7};
  \end{scope}
  \draw[black!25, line width=0.3pt] (-3,-3) rectangle (3,3);
  \node[anchor=north] at (0,-3.1) {\small $\chi_{5}$};
\end{tikzpicture}}%
\makebox[0.32\linewidth]{\begin{tikzpicture}[scale=0.7]
  \begin{scope}\clip (-3,-3) rectangle (3,3);
  \coordinate (p1) at (0,2.2);
  \coordinate (p2) at (-2.092,0.6798);
  \coordinate (p3) at (-1.293,-1.78);
  \coordinate (p4) at (1.293,-1.78);
  \coordinate (p5) at (1.516,-0.1138);
  \coordinate (p6) at (1.16,0.9828);
  \coordinate (p7) at (2.092,0.6798);
  \draw[black!45, line width=0.6pt] (p2) -- (p3) -- (p4) -- (p7) -- (p1) -- cycle;
  \draw[red!70!black, dashed, line width=0.5pt, shorten <=-12cm, shorten >=-12cm] (p1) -- (p4);
  \draw[blue!70!black, dashed, line width=0.5pt, shorten <=-12cm, shorten >=-12cm] (p1) -- (p5);
  \draw[blue!70!black, dashed, line width=0.5pt, shorten <=-12cm, shorten >=-12cm] (p2) -- (p7);
  \draw[red!70!black, dashed, line width=0.5pt, shorten <=-12cm, shorten >=-12cm] (p3) -- (p7);
  \draw[blue!70!black, dashed, line width=0.5pt, shorten <=-12cm, shorten >=-12cm] (p4) -- (p5);
  \fill[white] (p1) circle (3.3pt); \fill[black] (p1) circle (1.9pt);
  \node[font=\footnotesize, text=black, fill=white, fill opacity=0.8, text opacity=1, inner sep=0.6pt] at ([shift={(75.0:0.38)}]p1) {1};
  \fill[white] (p2) circle (3.3pt); \fill[black] (p2) circle (1.9pt);
  \node[font=\footnotesize, text=black, fill=white, fill opacity=0.8, text opacity=1, inner sep=0.6pt] at ([shift={(147.0:0.38)}]p2) {2};
  \fill[white] (p3) circle (3.3pt); \fill[black] (p3) circle (1.9pt);
  \node[font=\footnotesize, text=black, fill=white, fill opacity=0.8, text opacity=1, inner sep=0.6pt] at ([shift={(-111.0:0.38)}]p3) {3};
  \fill[white] (p4) circle (3.3pt); \fill[black] (p4) circle (1.9pt);
  \node[font=\footnotesize, text=black, fill=white, fill opacity=0.8, text opacity=1, inner sep=0.6pt] at ([shift={(-39.0:0.38)}]p4) {4};
  \fill[white] (p5) circle (3.3pt); \fill[red!70!black] (p5) circle (1.9pt);
  \node[font=\footnotesize, text=red!70!black, fill=white, fill opacity=0.8, text opacity=1, inner sep=0.6pt] at ([shift={(210:0.38)}]p5) {5};
  \fill[white] (p6) circle (3.3pt); \fill[blue!70!black] (p6) circle (1.9pt);
  \node[font=\footnotesize, text=blue!70!black, fill=white, fill opacity=0.8, text opacity=1, inner sep=0.6pt] at ([shift={(135:0.38)}]p6) {6};
  \fill[white] (p7) circle (3.3pt); \fill[black] (p7) circle (1.9pt);
  \node[font=\footnotesize, text=black, fill=white, fill opacity=0.8, text opacity=1, inner sep=0.6pt] at ([shift={(90:0.38)}]p7) {7};
  \end{scope}
  \draw[black!25, line width=0.3pt] (-3,-3) rectangle (3,3);
  \node[anchor=north] at (0,-3.1) {\small $\chi_{6}$};
\end{tikzpicture}}

\medskip
\makebox[0.32\linewidth]{\begin{tikzpicture}[scale=0.7]
  \begin{scope}\clip (-3,-3) rectangle (3,3);
  \coordinate (p1) at (0,2.2);
  \coordinate (p2) at (-2.092,0.6798);
  \coordinate (p3) at (-1.293,-1.78);
  \coordinate (p4) at (1.293,-1.78);
  \coordinate (p5) at (1.627,-0.0664);
  \coordinate (p6) at (1.479,0.4806);
  \coordinate (p7) at (2.092,0.6798);
  \draw[black!45, line width=0.6pt] (p2) -- (p3) -- (p4) -- (p7) -- (p1) -- cycle;
  \draw[red!70!black, dashed, line width=0.5pt, shorten <=-12cm, shorten >=-12cm] (p1) -- (p4);
  \draw[blue!70!black, dashed, line width=0.5pt, shorten <=-12cm, shorten >=-12cm] (p1) -- (p5);
  \draw[blue!70!black, dashed, line width=0.5pt, shorten <=-12cm, shorten >=-12cm] (p2) -- (p7);
  \draw[red!70!black, dashed, line width=0.5pt, shorten <=-12cm, shorten >=-12cm] (p3) -- (p7);
  \draw[blue!70!black, dashed, line width=0.5pt, shorten <=-12cm, shorten >=-12cm] (p4) -- (p5);
  \fill[white] (p1) circle (3.3pt); \fill[black] (p1) circle (1.9pt);
  \node[font=\footnotesize, text=black, fill=white, fill opacity=0.8, text opacity=1, inner sep=0.6pt] at ([shift={(75.0:0.38)}]p1) {1};
  \fill[white] (p2) circle (3.3pt); \fill[black] (p2) circle (1.9pt);
  \node[font=\footnotesize, text=black, fill=white, fill opacity=0.8, text opacity=1, inner sep=0.6pt] at ([shift={(147.0:0.38)}]p2) {2};
  \fill[white] (p3) circle (3.3pt); \fill[black] (p3) circle (1.9pt);
  \node[font=\footnotesize, text=black, fill=white, fill opacity=0.8, text opacity=1, inner sep=0.6pt] at ([shift={(-111.0:0.38)}]p3) {3};
  \fill[white] (p4) circle (3.3pt); \fill[black] (p4) circle (1.9pt);
  \node[font=\footnotesize, text=black, fill=white, fill opacity=0.8, text opacity=1, inner sep=0.6pt] at ([shift={(-39.0:0.38)}]p4) {4};
  \fill[white] (p5) circle (3.3pt); \fill[red!70!black] (p5) circle (1.9pt);
  \node[font=\footnotesize, text=red!70!black, fill=white, fill opacity=0.8, text opacity=1, inner sep=0.6pt] at ([shift={(210:0.38)}]p5) {5};
  \fill[white] (p6) circle (3.3pt); \fill[blue!70!black] (p6) circle (1.9pt);
  \node[font=\footnotesize, text=blue!70!black, fill=white, fill opacity=0.8, text opacity=1, inner sep=0.6pt] at ([shift={(90:0.38)}]p6) {6};
  \fill[white] (p7) circle (3.3pt); \fill[black] (p7) circle (1.9pt);
  \node[font=\footnotesize, text=black, fill=white, fill opacity=0.8, text opacity=1, inner sep=0.6pt] at ([shift={(300:0.38)}]p7) {7};
  \end{scope}
  \draw[black!25, line width=0.3pt] (-3,-3) rectangle (3,3);
  \node[anchor=north] at (0,-3.1) {\small $\chi_{7}$};
\end{tikzpicture}}%
\makebox[0.32\linewidth]{\begin{tikzpicture}[scale=0.7]
  \begin{scope}\clip (-3,-3) rectangle (3,3);
  \coordinate (p1) at (0,2.2);
  \coordinate (p2) at (-2.092,0.6798);
  \coordinate (p3) at (-1.293,-1.78);
  \coordinate (p4) at (1.293,-1.78);
  \coordinate (p5) at (1.588,-0.3942);
  \coordinate (p6) at (1.512,0.051);
  \coordinate (p7) at (2.092,0.6798);
  \draw[black!45, line width=0.6pt] (p2) -- (p3) -- (p4) -- (p7) -- (p1) -- cycle;
  \draw[red!70!black, dashed, line width=0.5pt, shorten <=-12cm, shorten >=-12cm] (p1) -- (p4);
  \draw[blue!70!black, dashed, line width=0.5pt, shorten <=-12cm, shorten >=-12cm] (p1) -- (p5);
  \draw[red!70!black, dashed, line width=0.5pt, shorten <=-12cm, shorten >=-12cm] (p3) -- (p7);
  \draw[blue!70!black, dashed, line width=0.5pt, shorten <=-12cm, shorten >=-12cm] (p4) -- (p5);
  \fill[white] (p1) circle (3.3pt); \fill[black] (p1) circle (1.9pt);
  \node[font=\footnotesize, text=black, fill=white, fill opacity=0.8, text opacity=1, inner sep=0.6pt] at ([shift={(75.0:0.38)}]p1) {1};
  \fill[white] (p2) circle (3.3pt); \fill[black] (p2) circle (1.9pt);
  \node[font=\footnotesize, text=black, fill=white, fill opacity=0.8, text opacity=1, inner sep=0.6pt] at ([shift={(162.0:0.38)}]p2) {2};
  \fill[white] (p3) circle (3.3pt); \fill[black] (p3) circle (1.9pt);
  \node[font=\footnotesize, text=black, fill=white, fill opacity=0.8, text opacity=1, inner sep=0.6pt] at ([shift={(-111.0:0.38)}]p3) {3};
  \fill[white] (p4) circle (3.3pt); \fill[black] (p4) circle (1.9pt);
  \node[font=\footnotesize, text=black, fill=white, fill opacity=0.8, text opacity=1, inner sep=0.6pt] at ([shift={(-39.0:0.38)}]p4) {4};
  \fill[white] (p5) circle (3.3pt); \fill[red!70!black] (p5) circle (1.9pt);
  \node[font=\footnotesize, text=red!70!black, fill=white, fill opacity=0.8, text opacity=1, inner sep=0.6pt] at ([shift={(195:0.38)}]p5) {5};
  \fill[white] (p6) circle (3.3pt); \fill[blue!70!black] (p6) circle (1.9pt);
  \node[font=\footnotesize, text=blue!70!black, fill=white, fill opacity=0.8, text opacity=1, inner sep=0.6pt] at ([shift={(105:0.38)}]p6) {6};
  \fill[white] (p7) circle (3.3pt); \fill[black] (p7) circle (1.9pt);
  \node[font=\footnotesize, text=black, fill=white, fill opacity=0.8, text opacity=1, inner sep=0.6pt] at ([shift={(3.0:0.38)}]p7) {7};
  \end{scope}
  \draw[black!25, line width=0.3pt] (-3,-3) rectangle (3,3);
  \node[anchor=north] at (0,-3.1) {\small $\chi_{8}$};
\end{tikzpicture}}%
\makebox[0.32\linewidth]{\begin{tikzpicture}[scale=0.7]
  \begin{scope}\clip (-3,-3) rectangle (3,3);
  \coordinate (p1) at (0,2.2);
  \coordinate (p2) at (-1.905,1.1);
  \coordinate (p3) at (-1.905,-1.1);
  \coordinate (p4) at (0,-2.2);
  \coordinate (p5) at (1.905,-1.1);
  \coordinate (p6) at (1.27,0.465);
  \coordinate (p7) at (1.905,1.1);
  \draw[black!45, line width=0.6pt] (p3) -- (p4) -- (p5) -- (p7) -- (p1) -- (p2) -- cycle;
  \draw[red!70!black, dashed, line width=0.5pt, shorten <=-12cm, shorten >=-12cm] (p1) -- (p5);
  \draw[red!70!black, dashed, line width=0.5pt, shorten <=-12cm, shorten >=-12cm] (p3) -- (p7);
  \draw[red!70!black, dashed, line width=0.5pt, shorten <=-12cm, shorten >=-12cm] (p4) -- (p7);
  \fill[white] (p1) circle (3.3pt); \fill[black] (p1) circle (1.9pt);
  \node[font=\footnotesize, text=black, fill=white, fill opacity=0.8, text opacity=1, inner sep=0.6pt] at ([shift={(82.0:0.38)}]p1) {1};
  \fill[white] (p2) circle (3.3pt); \fill[black] (p2) circle (1.9pt);
  \node[font=\footnotesize, text=black, fill=white, fill opacity=0.8, text opacity=1, inner sep=0.6pt] at ([shift={(150.0:0.38)}]p2) {2};
  \fill[white] (p3) circle (3.3pt); \fill[black] (p3) circle (1.9pt);
  \node[font=\footnotesize, text=black, fill=white, fill opacity=0.8, text opacity=1, inner sep=0.6pt] at ([shift={(-115.0:0.38)}]p3) {3};
  \fill[white] (p4) circle (3.3pt); \fill[black] (p4) circle (1.9pt);
  \node[font=\footnotesize, text=black, fill=white, fill opacity=0.8, text opacity=1, inner sep=0.6pt] at ([shift={(-82.0:0.38)}]p4) {4};
  \fill[white] (p5) circle (3.3pt); \fill[black] (p5) circle (1.9pt);
  \node[font=\footnotesize, text=black, fill=white, fill opacity=0.8, text opacity=1, inner sep=0.6pt] at ([shift={(-22.0:0.38)}]p5) {5};
  \fill[white] (p6) circle (3.3pt); \fill[red!70!black] (p6) circle (1.9pt);
  \node[font=\footnotesize, text=red!70!black, fill=white, fill opacity=0.8, text opacity=1, inner sep=0.6pt] at ([shift={(120:0.38)}]p6) {6};
  \fill[white] (p7) circle (3.3pt); \fill[black] (p7) circle (1.9pt);
  \node[font=\footnotesize, text=black, fill=white, fill opacity=0.8, text opacity=1, inner sep=0.6pt] at ([shift={(-5.0:0.38)}]p7) {7};
  \end{scope}
  \draw[black!25, line width=0.3pt] (-3,-3) rectangle (3,3);
  \node[anchor=north] at (0,-3.1) {\small $\chi_{9}$};
\end{tikzpicture}}

\medskip
\makebox[0.32\linewidth]{\begin{tikzpicture}[scale=0.7]
  \begin{scope}\clip (-3,-3) rectangle (3,3);
  \coordinate (p1) at (0,2.2);
  \coordinate (p2) at (-1.905,1.1);
  \coordinate (p3) at (-1.905,-1.1);
  \coordinate (p4) at (0,-2.2);
  \coordinate (p5) at (1.905,-1.1);
  \coordinate (p6) at (1.611,0);
  \coordinate (p7) at (1.905,1.1);
  \draw[black!45, line width=0.6pt] (p3) -- (p4) -- (p5) -- (p7) -- (p1) -- (p2) -- cycle;
  \draw[red!70!black, dashed, line width=0.5pt, shorten <=-12cm, shorten >=-12cm] (p1) -- (p5);
  \draw[red!70!black, dashed, line width=0.5pt, shorten <=-12cm, shorten >=-12cm] (p4) -- (p7);
  \fill[white] (p1) circle (3.3pt); \fill[black] (p1) circle (1.9pt);
  \node[font=\footnotesize, text=black, fill=white, fill opacity=0.8, text opacity=1, inner sep=0.6pt] at ([shift={(82.0:0.38)}]p1) {1};
  \fill[white] (p2) circle (3.3pt); \fill[black] (p2) circle (1.9pt);
  \node[font=\footnotesize, text=black, fill=white, fill opacity=0.8, text opacity=1, inner sep=0.6pt] at ([shift={(150.0:0.38)}]p2) {2};
  \fill[white] (p3) circle (3.3pt); \fill[black] (p3) circle (1.9pt);
  \node[font=\footnotesize, text=black, fill=white, fill opacity=0.8, text opacity=1, inner sep=0.6pt] at ([shift={(-150.0:0.38)}]p3) {3};
  \fill[white] (p4) circle (3.3pt); \fill[black] (p4) circle (1.9pt);
  \node[font=\footnotesize, text=black, fill=white, fill opacity=0.8, text opacity=1, inner sep=0.6pt] at ([shift={(-82.0:0.38)}]p4) {4};
  \fill[white] (p5) circle (3.3pt); \fill[black] (p5) circle (1.9pt);
  \node[font=\footnotesize, text=black, fill=white, fill opacity=0.8, text opacity=1, inner sep=0.6pt] at ([shift={(-22.0:0.38)}]p5) {5};
  \fill[white] (p6) circle (3.3pt); \fill[red!70!black] (p6) circle (1.9pt);
  \node[font=\footnotesize, text=red!70!black, fill=white, fill opacity=0.8, text opacity=1, inner sep=0.6pt] at ([shift={(285:0.38)}]p6) {6};
  \fill[white] (p7) circle (3.3pt); \fill[black] (p7) circle (1.9pt);
  \node[font=\footnotesize, text=black, fill=white, fill opacity=0.8, text opacity=1, inner sep=0.6pt] at ([shift={(22.0:0.38)}]p7) {7};
  \end{scope}
  \draw[black!25, line width=0.3pt] (-3,-3) rectangle (3,3);
  \node[anchor=north] at (0,-3.1) {\small $\chi_{10}$};
\end{tikzpicture}}%
\makebox[0.32\linewidth]{\begin{tikzpicture}[scale=0.7]
  \begin{scope}\clip (-3,-3) rectangle (3,3);
  \coordinate (p1) at (0,2.2);
  \coordinate (p2) at (-1.72,1.372);
  \coordinate (p3) at (-2.145,-0.4896);
  \coordinate (p4) at (-0.9546,-1.982);
  \coordinate (p5) at (0.9546,-1.982);
  \coordinate (p6) at (2.145,-0.4896);
  \coordinate (p7) at (1.72,1.372);
  \draw[black!45, line width=0.6pt] (p3) -- (p4) -- (p5) -- (p6) -- (p7) -- (p1) -- (p2) -- cycle;
  \fill[white] (p1) circle (3.3pt); \fill[black] (p1) circle (1.9pt);
  \node[font=\footnotesize, text=black, fill=white, fill opacity=0.8, text opacity=1, inner sep=0.6pt] at ([shift={(90.0:0.38)}]p1) {1};
  \fill[white] (p2) circle (3.3pt); \fill[black] (p2) circle (1.9pt);
  \node[font=\footnotesize, text=black, fill=white, fill opacity=0.8, text opacity=1, inner sep=0.6pt] at ([shift={(141.4:0.38)}]p2) {2};
  \fill[white] (p3) circle (3.3pt); \fill[black] (p3) circle (1.9pt);
  \node[font=\footnotesize, text=black, fill=white, fill opacity=0.8, text opacity=1, inner sep=0.6pt] at ([shift={(-167.1:0.38)}]p3) {3};
  \fill[white] (p4) circle (3.3pt); \fill[black] (p4) circle (1.9pt);
  \node[font=\footnotesize, text=black, fill=white, fill opacity=0.8, text opacity=1, inner sep=0.6pt] at ([shift={(-115.7:0.38)}]p4) {4};
  \fill[white] (p5) circle (3.3pt); \fill[black] (p5) circle (1.9pt);
  \node[font=\footnotesize, text=black, fill=white, fill opacity=0.8, text opacity=1, inner sep=0.6pt] at ([shift={(-64.3:0.38)}]p5) {5};
  \fill[white] (p6) circle (3.3pt); \fill[black] (p6) circle (1.9pt);
  \node[font=\footnotesize, text=black, fill=white, fill opacity=0.8, text opacity=1, inner sep=0.6pt] at ([shift={(-12.9:0.38)}]p6) {6};
  \fill[white] (p7) circle (3.3pt); \fill[black] (p7) circle (1.9pt);
  \node[font=\footnotesize, text=black, fill=white, fill opacity=0.8, text opacity=1, inner sep=0.6pt] at ([shift={(38.6:0.38)}]p7) {7};
  \end{scope}
  \draw[black!25, line width=0.3pt] (-3,-3) rectangle (3,3);
  \node[anchor=north] at (0,-3.1) {\small $\chi_{11}$};
\end{tikzpicture}}%
\par\endgroup
	\caption{Representative point configurations in $\cX_{3,7}$ for each of the 11 chirotopes $\chi_i$ from \zcref{fig:chirotope_table}. In each frame, the black points form a convex polygon. The red point is in the chamber bounded by edges of the polygon and dashed red lines through pairs of vertices. The blue point is in the chamber bounded by edges of the polygon and dashed red/blue lines through the other points.} \label{fig:chirotope_images}
\end{figure}
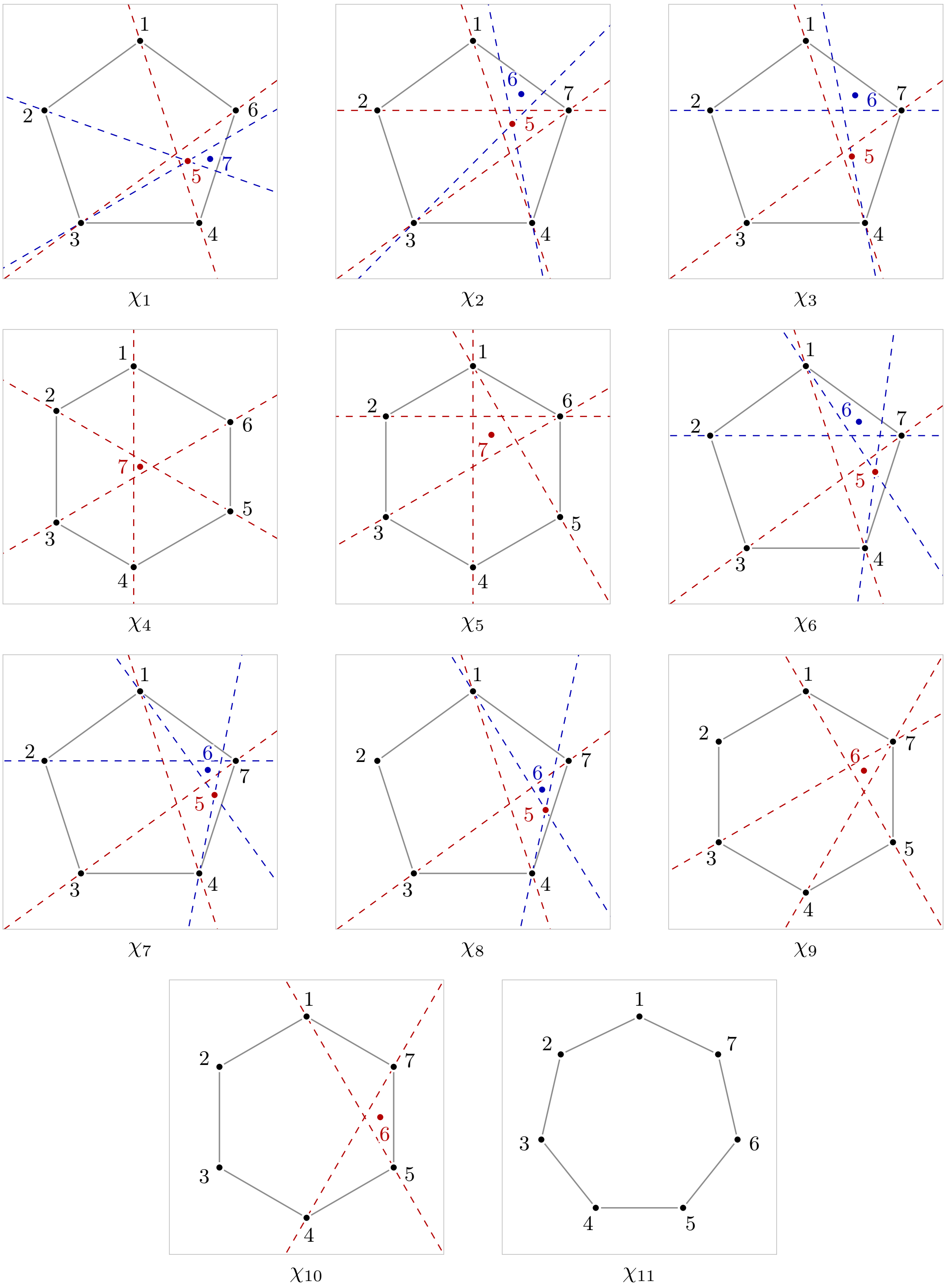

\begin{theorem}\label{thm:real37}
$\cromb_\R(\Gr_{3,7},Z_{3,7})=\cromb_\R(\cX_{3,7})=\cromb(\cX_{3,7})=7470$.
\end{theorem}

The first equality is \zcref{prop:descent-real}, and $\cromb_\R(\cX_{3,7})\le\cromb(\cX_{3,7})=7470$ by \zcref{prop:ff-X3n}. For the lower bound we find $7470$ regions and $7470$ corners of a tropical compactification whose matrix of iterated boundaries is invertible modulo $2$. We first outline the argument, which the rest of this section carries out in detail. At a corner the complement of the boundary consists of $64$ orthants. They lie in $64$ different oriented matroid strata, so a region contains at most one of them. Which strata occur can be read off from the initial degeneration at the corner: it is a coset of a six-dimensional torus whose real points carry exactly the chirotopes of the $64$ orthants. By the theorem of Antolini and Early \cite{Early:2024asu} these are the chirotopes whose chirotropical Dressian contains the cone of the corner. So modulo $2$ the matrix of iterated boundaries is the incidence matrix of the fans computed by Antolini and Early and we compute its rank over $\F_2$.

\subsubsection{Chirotropical Dressians} 

We  recall the chirotropical objects \cite[\S 2]{Early:2024asu}. Let $\chi \in \{\pm 1\}^{\binom{[n]}{3}}$ be a uniform chirotope and let $g \in \Q[p_I]$ be a polynomial in the Plücker coordinates. We write $p^v = \prod_I p_I^{v_I}$ and $\chi^v = \prod_I \chi_I^{v_I}$.  A term $cp^v$ of $g$ has the sign $\sign(c)\chi^v$ at $\chi$ and the weight $\sum_I v_I w_I$ at $w \in \R^{\binom{[n]}{3}}.$ The initial form $\inw(g)$ is the sum of the terms of $g$ of minimal weight at $w$ and $\Trop^\chi V(g)$ is the set of all $w$ for which the terms of $\inw(g)$ do not all have the same sign at $\chi$. 

The \emph{chirotropical Dressian} $\Dr^\chi(3,n)$ is the intersection of the $\Trop^\chi V(f)$ over the three-term Plücker relations
\begin{equation}\label{eq:three-term}
	p_{lab}\,p_{lce}\ -\ p_{lac}\,p_{lbe}\ +\ p_{lae}\,p_{lbc}\ =\ 0\qquad(a<b<c<e,\ l\notin\{a,b,c,e\}),
\end{equation}
where $p_{xyz}$ is the determinant of the columns $x,y,z$ in this order. The \emph{chirotropical Grassmannian} $\Trop^\chi\Gr(3,n)$ is the intersection  of the sets $\Trop^\chi V(g)$ over all polynomials $g$ in the Plücker ideal $I_{3,n} \subseteq \Q[p_I]$. Every index occurs equally often in the three terms of \eqref{eq:three-term}. Therefore $\Dr^\chi(3,n)$ only depends on the reorientation class of $\chi$, and it is invariant under translation by the lineality space $L$. Antolini and Early prove
\begin{equation}\label{eq:AE}
	\Dr^\chi(3,n)\ =\ \Trop^\chi\Gr(3,n)\qquad(n=6,7,8)
\end{equation}
for every uniform chirotope $\chi$ \cite[Thm.~1.1]{Early:2024asu}.

\subsubsection{Plücker cones}
Two vectors of the Dressian $\Dr(3,7)$ lie in the same \emph{Plücker cone} if they have the same terms of minimal weight in every relation \eqref{eq:three-term}. These polyhedral cones form the Plücker fan structure on the Dressian $\Dr(3,7)$ \cite[\S2]{HJJS}, and each $\Dr^\chi(3,7)$ is a union of Plücker cones. Plücker cones and the sets $\Dr^\chi(3,7)$ are invariant under translation by $L$, so from now on we regard them as subsets of $(N_X)_\R = N_\R / L$. Let $\mathcal P$ be the set of the $211\,365$ Plücker cones of dimension $6$. Their union is $\Trop\cX_{3,7}$ \cite[Thm.~2.2]{HJJS}. We note that the Dressian $\Dr(3,7)$ also has $30$ further maximal cones of dimension $7$ modulo $L$, coming from the Fano matroid, whose relative interiors do not meet $\Trop \cX_{3,7}$.

\subsubsection{The compactification and its corner}
 
Let $\Sigma$ be a unimodular refinement of the Gr\"obner fan structure on $\Trop\cX_{3,7}$. Let $(X,Y)$ be the tropical compactification of $\cX_{3,7}$ with respect to $\Sigma$, as in \zcref{ss:corners}. From now on $R_\chi$ denotes the closure of $S_\chi$ in $X(\R)$. Every initial degeneration of $\cX_{3,7}$ is smooth and irreducible \cite[Thm.~1.2]{Corey}. Hence, by \zcref{lem:corners}(3) the corners of $Y$ are the points $x_C$, one over every maximal cone $C$ of $\Sigma$, and they are defined over $\R$. 

The Gröbner fan structure on $\Trop \cX_{3,7}$ refines its Plücker fan structure. So for a maximal cone $C$ of $\Sigma$ the relative interior of $C$ lies in a single Plücker cone. This cone has dimension at least $6$ and meets $\Trop \cX_{3,7}$, so it is not the relative interior of one of the $30$ cones of dimension $7$.  Hence it is the relative interior of a cone $P(C) \in \mathcal P$, i.e. $\relint C \subseteq \relint P(C)$.  Conversely, every $P \in \mathcal P$ has dimension $6$ and lies in the support of $\Sigma$. So $\relint P$ meets the relative interior of some maximal cone $C$ of $\Sigma$, and then  $P =P(C)$.

For a maximal cone $C$ of $\Sigma$ we use the chart $U_C \cong \A^6 \times \G_m^{22}$ of the proof of \zcref{lem:corners}. Since $C$ is unimodular, its primitive ray generators $u_1,\ldots,u_6$ extend to a $\Z$-basis $u_1,\ldots,u_{28}$ of $N_X$.  The coordinates $z_1,\ldots,z_{28}$ of $U_C$ are the dual basis of the character lattice $M_X$. So $C = \RR^6_{\geq 0} \times \{0\}$ in the basis $u_1,\ldots,u_{28}$ of $N_X$. The orbit $O_C$ is $\{z_1 = \dots = z_6 = 0\}$  and $x_C = (0,a)$ with $a \in (\Q^\times)^{22}$. We choose the neighborhood $V$ of $x_C$ in $X(\R)$ as in \eqref{eq:orthants} so that $(z_1,\ldots,z_6)$ maps it homeomorphically to an open box around $0$ in $\R^6$, and that $\sign z_k = \sign a_k$ on $V$ for $7 \leq k \leq 28$. Then $x_C$ lies in the closure of each of the $64$ orthants $V_\varepsilon,\ \varepsilon \in \{\pm 1\}^6$.

The chirotope of a real configuration is recorded by the signs of the torus coordinates. Each $z_k$ is a Laurent monomial in the Plücker coordinates in which every index occurs equally often in the numerator and in the denominator. So its sign at a real point of $T_P$ only depends on the reorientation class of the sign vector at that point. In particular, the signs of $z_1,\ldots,z_{28}$ are the same at all real points with a given chirotope $\chi$. We write $\sign_z(\chi) \in \{\pm 1\}^{28}$ for this sign vector, which determines $\chi$. We can view $S_\chi$ as the intersection of $\cX_{3,7}(\R)$ with the orthant of $T_X(\R)$ with sign vector $\sign_z(\chi)$.

At a corner the initial degeneration is a torus coset, and its real points have the same sign vectors as the orthants.

\begin{lemma}\label{lem:corner-initial}
	Let $C$ be a maximal cone of $\Sigma$ and $w\in\relint C$. Then
	\[
	\inw\cX_{3,7}\ =\ \{\,z\in T_X:\ z_k=a_k\ \text{ for }7\le k\le 28\,\}.
	\]
	In particular, the sign vectors $(\sign z_j)_{j\le28}$ of the real points of
	$\inw\cX_{3,7}$ are exactly the $64$ vectors $(\varepsilon,\sign a)$,
	$\varepsilon\in\{\pm1\}^6$, and these are the sign vectors of the $64$ orthants
	$V_\varepsilon$ at $x_C$.
\end{lemma}

\begin{proof}
	The first statement follows from the isomorphism $\inw \cX_{3,7} \cong X_C \times \G_m^6$ where we use $X_C = \{x_C\}$ and $(z_7,\ldots,z_{28})(x_C) = a$.
	
	Its real points are the points $(s,a)$ with $s \in (\R^\times)^6$; their sign vectors are $(\sign s, \sign a)$, and every $\varepsilon \in \{\pm 1\}^6$ occurs. On $V_\varepsilon$ the coordinate $z_j$ has the sign $\varepsilon_j$ for $j \leq 6$ and $z_k$ the sign of $a_k$ for $k \geq 7$.
\end{proof}

The next proposition identifies the regions through a corner.

\begin{proposition}\label{prop:orthants-dressian}
	Let $C$ be a maximal cone of $\Sigma$, let $w\in\relint C$, and let
	$\chi\in\mathcal K$. The following are equivalent.
	\begin{enumerate}[(1)]
		\item The corner $x_C$ lies in the region $R_\chi$.
		\item $S_\chi\cap V=V_\varepsilon$ for some $\varepsilon\in\{\pm1\}^6$.
		\item $\sign_z(\chi)$ is the sign vector of a real point of $\inw\cX_{3,7}$, i.e. $\sign_z(\chi) = (\varepsilon, \sign a)$ for some $\varepsilon \in \{\pm 1\}^6$.
		\item $P(C)\subseteq\Dr^\chi(3,7)$.
	\end{enumerate}
	In this case $\varepsilon$ is unique and $R_\chi\cap V$ is the closure of 	$V_\varepsilon$ in $V$. Exactly $64$ classes $\chi\in\mathcal K$ satisfy these
	conditions.
\end{proposition}

\begin{proof}
By the choice of $V$, the coordinates $z_1,\ldots,z_{28}$ have the signs $(\varepsilon,\sign a)$ on $V_\varepsilon$. The orthants are nonempty subsets of $\cX_{3,7}(\R)$, and $\sign_z$ determines the chirotope. So $V_\varepsilon \subseteq S_\chi$ if $\sign_z(\chi) = (\varepsilon,\sign a)$, and $V_\varepsilon \cap S_\chi = \emptyset$ otherwise. As $V \setminus Y(\R)$ is the disjoint union of the orthants, $S_\chi \cap V$ is $V_\varepsilon$ if (3) holds for $\varepsilon$, and it is empty if (3) fails. This proves (2)$\Leftrightarrow$(3) and the uniqueness of $\varepsilon$. Since $V$ is an open neighborhood of $x_C$ in $X(\R)$, the set $R_\chi \cap V$ is the closure of $S_\chi \cap V$ in $V$. It contains $x_C$ if $S_\chi \cap V = V_\varepsilon$, and it is empty otherwise. This proves (1)$\Leftrightarrow$(2) and the description of $R_\chi \cap V$.

For the remaining two implications recall from \eqref{eq:in-split} that $\inw\Groo_{3,7}$ is the preimage of $\inw\cX_{3,7}$ in the Plücker torus $T_P$.

(3)$\Rightarrow$(4). Let $\varepsilon$ be as in (3), and let $\hat y \in T_P(\R)$ be the Plücker vector of a real matrix representing a configuration in $V_\varepsilon$. Its sign vector lies in $\chi$, and its image in $T_X$ lies in the orthant with sign vector $(\varepsilon, \sign a)$. Multiplying the coordinates of $\hat y$ by positive reals changes neither of these facts, and it moves the image to any point of that orthant, since the positive part of $T_P(\R)$ maps onto the positive part of $T_X(\R)$. So we may assume that the image is the point $(\varepsilon, a)$ of $\inw\cX_{3,7}$ from \zcref{lem:corner-initial}. Then $\hat y$ is a real point of $\inw\Groo_{3,7}$. Let $f$ be a three-term relation \eqref{eq:three-term}. Its initial form $\inw(f)$ lies in $\inw(I_{3,7})$, so it vanishes at $\hat y$, while each of its terms is nonzero at $\hat y$ and has there the sign it has at the sign vector of $\hat y$. So the terms of $\inw(f)$ do not all have the same sign at $\chi$, that is, $w \in \Trop^\chi V(f)$. Hence $w \in \Dr^\chi(3,7)$. Since $\Dr^\chi(3,7)$ is a union of Plücker cones and $w \in \relint P(C)$, we get $P(C) \subseteq \Dr^\chi(3,7)$.

(4)$\Rightarrow$(3). By \eqref{eq:AE}, applied to a representative of $\chi$, we have $w\in\Trop^\chi\Gr(3,7)$. Suppose that (3) fails. Then there is $k \ge 7$ with $\sign_z(\chi)_k = -\sign a_k$. Write $z_k = p^{m^+}/p^{m^-}$ with $m^\pm \ge 0$ and $m^+ - m^- \in M_X$. By \zcref{lem:corner-initial} the function $z_k - a_k$ vanishes on $\inw\cX_{3,7}$, so $p^{m^+} - a_k\,p^{m^-}$ vanishes on $\inw\Groo_{3,7}$. This scheme is cut out of $T_P$ by $\inw(I_{3,7})$, and it is reduced, being smooth \cite[Thm.~1.2]{Corey}. Hence there is a monomial $p^u$ with
\[
h\ :=\ p^u\bigl(p^{m^+}-a_k\,p^{m^-}\bigr)\ \in\ \inw(I_{3,7}).
\]
The two terms of $h$ have the same weight at $w$, because $z_k$ has weight zero on $C$. So $h$ is the initial form of an element of $I_{3,7}$: write $h = \sum_i c_i\,p^{v_i} \inw(g_i)$ with $g_i \in I_{3,7}$, keep the summands of the same weight as $h$, and put $g = \sum_i c_i\,p^{v_i} g_i$; then $\inw(g) = h$. At $\chi$ the two terms of $h$ have the signs $\chi^{u+m^+}$ and $-\sign(a_k)\,\chi^{u+m^-}$, whose quotient is $-\sign(a_k)\,\sign_z(\chi)_k = 1$. So the terms of $\inw(g)$ all have the same sign at $\chi$, that is, $w \notin \Trop^\chi V(g)$, a contradiction.

Finally, the $64$ orthants $V_\varepsilon$ are nonempty subsets of $\cX_{3,7}(\R)$ with pairwise different sign vectors, so they lie in $64$ different strata $S_{\chi_\varepsilon}$, and $\chi_\varepsilon$ satisfies (3) with this $\varepsilon$. Conversely, a class satisfying (3) with $\varepsilon$ equals $\chi_\varepsilon$, as $\sign_z$ determines the class.
\end{proof}

Let $B\in\F_2^{\mathcal K\times\mathcal P}$ be the incidence matrix with
$B_{\chi,P}=1$ if and only if $P\subseteq\Dr^\chi(3,7)$.

\begin{corollary}\label{cor:B-columns}
	\begin{enumerate}[(1)]
		\item Every column of $B$ has exactly $64$ nonzero entries.
		\item For every maximal cone $C$ of $\Sigma$ and every $\chi\in\mathcal K$,
		\begin{equation}\label{eq:boundary-mod2}
			\partial_{x_C}R_\chi\ \equiv\ B_{\chi,P(C)}\pmod 2 .
		\end{equation}
	\end{enumerate}
\end{corollary}

\begin{proof}
	Every $P\in\mathcal P$ is $P(C)$ for some maximal cone $C$ of $\Sigma$, so (1)
	follows from \zcref{prop:orthants-dressian}. For (2), if
	$B_{\chi,P(C)}=1$, then $R_\chi\cap V$ is the closure of a single orthant by
	\zcref{prop:orthants-dressian}, and $\partial_{x_C}R_\chi=\pm1$ by
	\zcref{lem:realcorners}(2). Otherwise $x_C\notin R_\chi$, and
	$\partial_{x_C}R_\chi=0$ by \zcref{lem:realcorners}(1).
\end{proof}

We build $B$ from the data of Antolini and Early \cite{Early:2024asu}. Let $\psi_0,\dots,\psi_{10}$ be their eleven chirotopes, one in each reorientation class up to relabeling. For each $\psi_j$ they list the maximal cones of the fan $\Dr^{\psi_j}(3,7)/L$ by their rays \cite[Thm.~5.3, \S7]{Early:2024asu}, and a relabeling $g\in S_7$ maps $\Dr^{\psi_j}(3,7)$ onto $\Dr^{g\cdot\psi_j}(3,7)$. For a maximal cone $F$ of the fan of $\psi_j$ let $w_F$ be the sum of its rays. Since $\Dr^{g\cdot\psi_j}(3,7)$ is a union of Plücker cones, it contains the Plücker cone of $g\cdot w_F$, which we read off from the terms of minimal weight of the relations \eqref{eq:three-term} at $g\cdot w_F$. A computer search over all $g$, $j$ and $F$ yields $211\,365$ distinct Plücker cones, all of dimension $6$, and $13\,527\,360 = 64 \cdot 211\,365$ distinct pairs of a class $g\cdot\psi_j$ and such a cone, exactly $64$ for each cone. Each such pair is a nonzero entry of $B$. Since $|\mathcal P| = 211\,365$ and every column of $B$ has exactly $64$ nonzero entries by \zcref{cor:B-columns}, these pairs are all the nonzero entries of $B$. 

The matrix $B$ and the code are available on the Zenodo page \cite{KL26data}. The data includes:
\begin{itemize}
	\item the matrix $B$, with each class given by a sign vector and as $g\cdot\psi_j$, and each cone $P\in\mathcal P$ by an integer point in its relative interior,
	\item the file \texttt{build\_B.sage}, which builds $B$ from the data of \cite{Early:2024asu} as described above,
	\item the file \texttt{rank\_B.sage}, which computes the rank of $B$ over $\F_2$.
\end{itemize}

\begin{proposition}\label{prop:certificate}
	The matrix $B$ has $\operatorname{rank}_{\F_2}B=7470$.
\end{proposition}

\begin{proof}
	The matrix $B$ is a $27\,240\times211\,365$ matrix over $\F_2$ with $13\,527\,360$ ones. We computed its rank with \texttt{rank\_B.sage}, see \cite{KL26data}.
\end{proof}

\begin{proof}[Proof of \zcref{thm:real37}]
	It remains to show $\cromb_\R(\cX_{3,7})\ge7470$. Choose for every $P\in\mathcal P$ a maximal cone $C_P$ of $\Sigma$ with $P(C_P)=P$. By \eqref{eq:boundary-mod2} the integer matrix $(\partial_{x_{C_P}}R_\chi)_{\chi\in\mathcal K,\,P\in\mathcal P}$ is congruent to $B$ modulo $2$. The rank over $\Q$ of an integer matrix is at least the rank over $\F_2$ of its reduction modulo $2$, so by \zcref{prop:certificate} this matrix has rank at least $7470$. Hence the vectors $(\partial_cR_\chi)_{c\in\mathcal C}$ of $7470$ regions $R_\chi$ are linearly independent, and $\cromb_\R(\cX_{3,7})\ge7470$ by \eqref{eq:crR-matrix}.
\end{proof}

\begin{remark}\label{remark:cegm_integrands} 
Cachazo, Early and Zhang construct the integrands of \cite{Cachazo:2019ngv} on $\cX_{3,n}$ region by region \cite[\S3, Algorithm~I]{Cachazo:2023ltw}. The integrand of a region has simple poles along the Pl\"ucker divisors bounding it,\footnote{One way to find these boundaries is to study the mutations of uniform oriented matroids \cite{BLVSWZ}.} and its numerator is a polynomial of known degree. Starting from the integrands whose numerator is fixed by the torus weight, each further numerator is determined by requiring that its residue along every wall $Z_{\{I\}}$ agrees, up to sign, with that of the region across the wall. For $n\le8$ this determines every integrand up to sign \cite[\S7]{Cachazo:2023ltw}.

    For $6\le n\le8$ this construction gives the canonical forms. Let $I\in\binom{[n]}{3}$ and $Z_{\mathrm{rest}}=Z_{\binom{[n]}{3}\setminus\{I\}}$. Since $\Gr_{3,n}\setminus Z_{\mathrm{rest}}$ is smooth and affine, \cite[Cor.~3.12]{Brown:2025jjg} gives
    \[
    0=g(\Gr_{3,n},Z_{3,n})=g(\Gr_{3,n},Z_{\mathrm{rest}})+g\bigl(Z_{\{I\}},Z_{\{I\}}\cap Z_{\mathrm{rest}}\bigr),
    \]
    and both terms vanish as they are nonnegative. So \cite[Prop.~2.15]{Brown:2025jjg} applies: $\Res_{Z_{\{I\}}}\omega(\sigma)=\omega(\partial_{Z_{\{I\}}}\sigma)$ for every region $\sigma$. Hence $\omega(\sigma)$ has poles only along the walls of $\sigma$, and two regions sharing a wall have canonical forms whose residues there agree up to sign. Thus the canonical forms satisfy every constraint of the construction, and since the constraints have a unique solution, induction along the construction shows that each integrand is a nonzero multiple of $\omega(\sigma)$. The residues along a common wall are nonzero, so neighboring regions carry the same multiple up to sign, and since any two regions are connected through walls, the integrands are $\pm\lambda\,\omega(\sigma)$ for a single constant $\lambda$ set by the normalization of the measure. So $\cromb_\R(\cX_{3,n})$ is the dimension of the span of the integrands. For $n=6,7$ this also gives a way to compute the real combinatorial rank. For $n=8$ the same holds in principle, but computing this span is out of reach.
\end{remark}

\subsection{Eight points}\label{ss:38-real}

For eight points not every corner is real. Let $Q_{\mathrm{mk}}$ be the matroid of rank $3$ on $[8]$ whose nonbases are
\begin{equation}\label{eq:Qmk-nonbases}
	126,\quad 145,\quad 178,\quad 235,\quad 248,\quad 347,\quad 368,\quad 567.
\end{equation}
They are the eight lines of the M\"obius--Kantor configuration, and Corey and Luber call $Q_{\mathrm{mk}}$ the M\"obius--Kantor or MacLane matroid \cite[\S1]{CoreyLuber}. Up to relabeling it is the sparse paving matroid of \cite[Rem.~4.5]{BBRS}. It is realizable over $\Q(\sqrt{-3})$ but not over $\R$, since its realizations come from the two roots of $a^2-a+1$ \cite[Prop.~4.2]{CoreyLuber}, as recalled below. Let $\Sigma_8$ be the fan structure on $\Trop\cX_{3,8}$ given by matroid subdivisions of the hypersimplex $\Delta(3,8)$ \cite[\S2.3]{CoreyLuber}. The indicator vector $w_{\mathrm{mk}}\in N = \Z^{\binom{8}{3}}$ of the nonbases \eqref{eq:Qmk-nonbases} lies in the relative interior of the cone
\begin{equation}\label{eq:Cmk}
	C_{\mathrm{mk}}\ =\ \R_{\ge0}\bigl\langle\,
	\overline{\mathbf e}_{126},\ \overline{\mathbf e}_{145},\ \overline{\mathbf e}_{178},\ \overline{\mathbf e}_{235},\
	\overline{\mathbf e}_{248},\ \overline{\mathbf e}_{347},\ \overline{\mathbf e}_{368},\ \overline{\mathbf e}_{567}\,\bigr\rangle
\end{equation}
of $\Sigma_8$ (see \cite[\S5]{CoreyLuber}). It is a simplicial maximal cone, and its orbit $\mathcal O=S_8\cdot C_{\mathrm{mk}}$ consists of $840$ maximal cones \cite[Rem.~4.5]{BBRS}. So the stabilizer of $C_{\mathrm{mk}}$ in $S_8$, the group of permutations of $[8]$ that preserve the eight lines \eqref{eq:Qmk-nonbases}, has order $8!/840=48$. Corey and Luber show that $\Groo_{3,8}$, and hence $\cX_{3,8}$ by \eqref{eq:in-split}, is sch\"on \cite[Thm.~1.1]{CoreyLuber}. Also $\inw\cX_{3,8}$ is irreducible unless $w$ lies in the relative interior of a cone of $\mathcal O$ \cite[Thm.~1.4]{CoreyLuber}, and it has two connected components for $w$ in the relative interior of a cone of $\mathcal O$ \cite[Thm.~1.3]{CoreyLuber}.

\begin{theorem}\label{thm:38-realrank}
	The real combinatorial rank of eight points is strictly smaller than the combinatorial rank:
	\[
	\cromb_{\R}\bigl(\Gr_{3,8},Z_{3,8}\bigr)=\cromb_{\R}\bigl(\cX_{3,8}\bigr)
	\ \le\ \cromb\bigl(\cX_{3,8}\bigr)-840
	\ <\ \cromb\bigl(\cX_{3,8}\bigr)=\cromb\bigl(\Gr_{3,8},Z_{3,8}\bigr).
	\]
\end{theorem}

We first show that the cones of $\mathcal O$ are unimodular, so that the tropical compactification contains their strata. Over each of them lie two corners, which are complex conjugate and which no region reaches, while their difference is a cycle of the dual complex.

\begin{proposition}\label{prop:Cunimodular}
	Every cone $C\in\mathcal O$ is unimodular, that is, its primitive ray generators are part of a basis of $N_X$.
\end{proposition}

\begin{proof}
	By symmetry it suffices to consider $C_{\mathrm{mk}}$. Let $I_1,\dots,I_8$ be the nonbases \eqref{eq:Qmk-nonbases}, let $A\in\Z^{16\times56}$ be the matrix with rows
	\[
	\mathbf e_{I_1},\ \dots,\ \mathbf e_{I_8},\ \varphi(\mathbf e_1),\ \dots,\ \varphi(\mathbf e_7),\ \mathbf 1,
	\]
	and let $\Gamma\subseteq N$ be the lattice spanned by its rows. The maximal minor of $A$ with the columns $I_1,\dots,I_8$ and $123,124,125,127,128,134,136,234$ equals $\pm1$. Hence $A$ has rank $16$ and $\Gamma$ is saturated in $N$, since the index of the row lattice of an integer matrix of full row rank in its saturation is the greatest common divisor of its maximal minors. Since $\varphi(\mathbf e_8)=3\cdot\mathbf 1-\sum_{i\le7}\varphi(\mathbf e_i)$, the last eight rows span $L$ over $\R$, so $N\cap L\subseteq(\Gamma\otimes\R)\cap N=\Gamma$ and $\Gamma=\sum_j\Z\,\mathbf e_{I_j}+(N\cap L)$. Since $N_X/(\Gamma/(N\cap L))\cong N/\Gamma$ is torsion free, the lattice $\Gamma/(N\cap L)$ is saturated in $N_X$. It has rank $8$ and is spanned by $\overline{\mathbf e}_{I_1},\dots,\overline{\mathbf e}_{I_8}$, so these eight vectors are a basis of a direct summand of $N_X$. In particular they are the primitive ray generators of $C_{\mathrm{mk}}$.
\end{proof}

By \zcref{prop:Cunimodular} and \cite[Thm.~11.1.9]{CLS} there is a unimodular refinement $\widetilde\Sigma$ of $\Sigma_8$ that contains every cone of $\mathcal O$. Let $(X,Y)$ be the tropical compactification of $\cX_{3,8}$ with respect to $\widetilde\Sigma$. By \zcref{lem:corners} and \cite[Thm.~1.4]{CoreyLuber} there is exactly one corner over every maximal cone of $\widetilde\Sigma$ that is not in $\mathcal O$, and it is real. By \cite[Thm.~1.3]{CoreyLuber} there are two corners over every cone $C\in\mathcal O$, as Corey and Luber also note for their compactification \cite[proof of Thm.~7.2]{CoreyLuber}.

The two corners over a cone $C\in\mathcal O$ are complex conjugate, and neither of them is real. By symmetry we may take $C=C_{\mathrm{mk}}$ and $w=w_{\mathrm{mk}}$. Let $\Gr(Q_{\mathrm{mk}})\subseteq\Gr_{3,8}$ be the thin Schubert cell of the points whose nonzero Pl\"ucker coordinates are those of the bases of $Q_{\mathrm{mk}}$. Every $3\times8$ matrix representing $Q_{\mathrm{mk}}$ is, up to the actions of $\mathrm{GL}_3$ and of the torus, the matrix $A'(a)$ of \cite[(11)]{CoreyLuber}, where $a$ is a root of $a^2-a+1$. Both roots are non-real, and $\Gr(Q_{\mathrm{mk}})$ has two connected components, one for each root \cite[Prop.~4.2]{CoreyLuber}. Complex conjugation maps $A'(a)$ to $A'(\bar a)$, so it exchanges these two components. By \cite[Prop.~5.1, Thm.~5.2 and their proofs]{CoreyLuber} the two connected components of $\operatorname{in}_{w_{\mathrm{mk}}}\Groo_{3,8}$ correspond to the two components of $\Gr(Q_{\mathrm{mk}})$, through morphisms defined over $\Q$, so complex conjugation exchanges them as well. The isomorphisms \eqref{eq:in-split} and \eqref{eq:in-stratum} are also defined over $\Q$, and they identify these two components with the two points of $X_C$. Hence complex conjugation exchanges the two corners over $C$, and neither of them is real.

\begin{lemma}\label{lem:38-cycle}
	Let $C\in\mathcal O$ and let $c,c'$ be the two corners in $X_C$. Then $\partial c=\partial c'$ for the boundary map in \eqref{eq:W}, so $z_C:=c-c'$ is a nonzero element of $W$.
\end{lemma}

\begin{proof}
	Let $\rho_1,\dots,\rho_8$ be the rays of $C$, so that $c$ and $c'$ both lie exactly on the components $Y_{\rho_1},\dots,Y_{\rho_8}$. For $j\in[8]$ let $\tau_j$ be the facet of $C$ spanned by the rays other than $\rho_j$, and let $Y^{(j)}=\bigcap_{i\ne j}Y_{\rho_i}$. This is a smooth curve, the union of the stratum $X_{\tau_j}$ and of the finitely many corners over the maximal cones containing $\tau_j$. Every connected component of $Y^{(j)}$ is a curve, so it meets the open subset $X_{\tau_j}$. We have that $\inw\cX_{3,8}$ is irreducible for $w\in\relint\tau_j$ \cite[Thm.~1.4]{CoreyLuber}, and $X_{\tau_j}$ is irreducible by \eqref{eq:in-stratum}. Hence $Y^{(j)}$ is connected, and $c$ and $c'$ lie on the same connected component of $Y^{(j)}$ for every $j$. Both corners lie on the same components, so the terms of $\partial c$ and $\partial c'$ in \eqref{eq:W} coincide. So $z_C\in\ker\partial=W$, and $z_C\ne0$ because $c\ne c'$.
\end{proof}

\begin{proof}[Proof of \zcref{thm:38-realrank}]
	We compute both ranks using the tropical compactification $(X,Y)$. By \eqref{eq:W}, $\cromb(\cX_{3,8})=\dim_\Q W$. Let $\mathcal C_{\mathcal O}$ be the set of the $2\cdot840$ corners over the cones of $\mathcal O$, and let $r\colon W\to\Q^{\mathcal C_{\mathcal O}}$ be the restriction of a vector in $\Q^{\mathcal C}$ to these corners. By \zcref{lem:38-cycle}, $W$ contains the $840$ vectors $z_C$, $C\in\mathcal O$. Distinct cones of $\mathcal O$ have disjoint sets of corners, so the vectors $r(z_C)$ have pairwise disjoint supports and are linearly independent. Hence $r$ has rank at least $840$ and $\dim_\Q\ker r\le\cromb(\cX_{3,8})-840$. For a region $\sigma$ the vector $(\partial_c\sigma)_{c\in\mathcal C}$ lies in $W$, and it vanishes on $\mathcal C_{\mathcal O}$ by \zcref{lem:realcorners}(1), since none of these corners is real, as shown before \zcref{lem:38-cycle}. By \eqref{eq:crR-matrix},
	\[
	\cromb_{\R}\bigl(\cX_{3,8}\bigr)\ \le\ \dim_\Q\ker r\ \le\ \cromb\bigl(\cX_{3,8}\bigr)-840.
	\]
\end{proof}

Together with \zcref{thm:posarr-3n}, this proves \zcref{thm:main38}.

\begin{remark}\label{rem:38-missing}
	The proof exhibits the forms that the regions miss. By \eqref{eq:W}, each $z_C$ is the residue vector of a unique $\omega_C\in\Omega^8_{\log}(\cX_{3,8})$, with iterated residues $\pm1$ at the two corners over $C$ and $0$ at all other corners, and no nonzero linear combination of the $840$ forms $\omega_C$ lies in the span of the canonical forms of the regions.
\end{remark}

For six and seven points the homology of the link of the tropical variety is known. The link of
$\Trop\cX_{3,6}$ has reduced integral homology only in degree $3$, free of rank $126$
\cite[Thm.~5.4]{SpeyerSturmfels}, and the link of $\Trop\cX_{3,7}$ has reduced integral homology
only in degree $5$, free of rank $7470$ \cite[Thm.~2.1]{HJJS}. In both cases the rank is the
combinatorial rank (\zcref{thm:main-posarr}), as explained by \cite[Thm.~2.5, Cor.~3.3]{Hacking},
since all strata of a tropical compactification are smooth and connected \cite[Ex.~4.2]{Hacking},
\cite[Thm.~1.2]{Corey}. For eight points the M\"obius--Kantor cones make the two numbers differ.

\begin{corollary}\label{cor:link38}
	The link of $\Trop\cX_{3,8}$ has the rational homology of a wedge of $1\,115\,640$ spheres of
	dimension~$7$.
\end{corollary}

\begin{proof}
	Let $(X,Y)$ be the tropical compactification with respect to $\widetilde\Sigma$ from
	\zcref{ss:38-real}. Its strata are smooth, and those of positive dimension are
	connected by \cite[Thm.~1.1, Thm.~1.4]{CoreyLuber} and \eqref{eq:in-stratum}, since the cones
	of $\mathcal O$ are maximal. So by \cite[Thm.~2.5]{Hacking} the link $L$ of
	$\Trop\cX_{3,8}$ has reduced rational homology only in degree~$7$. By the proof of
	\cite[Thm.~2.5]{Hacking}, the dual complex $\Delta(Y)$ arises from $L$ by replacing each
	maximal cell $C$ by $|X_C|$ copies glued along their boundaries. By \zcref{lem:corners} and
	\cite[Thm.~1.3, Thm.~1.4]{CoreyLuber}, exactly the $840$ cones of $\mathcal O$ carry two
	corners, so $\dim\widetilde H_7(\Delta(Y);\Q)=\dim\widetilde H_7(L;\Q)+840$. The left side is
	$\dim_\Q W=\cromb(\cX_{3,8})=1\,116\,480$ by \eqref{eq:W} and \zcref{prop:ff-X3n}. Hence
	$\dim\widetilde H_7(L;\Q)=1\,115\,640$.
\end{proof}

\section{Positive genus for ten and more points}\label{sec:main310}

In this section we show that $g(\cX_{3,10})\ge1$, and we deduce that the genus of $\cX_{k,n}$
is positive in a large range of $k$ and $n$. The point counts of \zcref{ss:finitefield} cannot
show this, since they only detect genus zero. Instead we find the genus in the boundary: one
stratum in the boundary of $\cX_{3,10}$ is an open subset of a K3 surface, and we show that its
genus survives in $\cX_{3,10}$ (\zcref{ss:ten-points}). Adding points and the duality of
Grassmannians then give the general statement (\zcref{ss:genus-range}). As a consequence of
\zcref{cor:pointcount}, the number of points of $\cX_{3,n}$ over $\F_q$ is not a
quasi-polynomial in $q$ for $n\ge10$.

\subsection{Ten points and a K3 surface}\label{ss:ten-points}

We label the ten points by $1,\dots,9,0$ as in \cite{Sturmfels1991,Sink-K3surface} and consider the ten triples
\[
\mathcal L_V:=\{124,138,179,237,259,350,456,480,678,690\}.
\]
Every point lies in three of them. Let $\mathcal Y\supset\cX_{3,10}$ be the space of configurations of ten points in $\PP^2$, up to $\PGL_3$, in which no triple outside $\mathcal L_V$ is collinear. Let $\mathcal R_{\mathcal L_V}\subset\mathcal Y$ be the realization space of $\mathcal L_V$, in which exactly the triples of $\mathcal L_V$ are collinear. By \cite[Thm.~1.2]{Sink-K3surface}, $\mathcal R_{\mathcal L_V}$ is isomorphic to a Zariski open subset of a K3 surface $\tilde S$, hence  $g(\mathcal R_{\mathcal L_V}) \geq g(\tilde S)=1$. It is a stratum of codimension ten in $\mathcal Y$. Note that a boundary stratum of positive genus does not force the genus of a pair to be positive \cite[\S2.4.3, \S5.9]{Brown:2025jjg}. Here it does, because $\mathcal R_{\mathcal L_V}$ is reached from $\cX_{3,10}$ by imposing the ten collinearities one at a time, and at each step the stratum together with the next one is a smooth affine variety.

We use the affine chart of $\mathcal Y$ by the matrices \cite{Sink-K3surface,Sturmfels1991}
\[
    \begin{pmatrix}
        1 & 0 & 0 & 1 & 1 & a & x_7 & x_8 & x_9 & x_0 \\ 0 & 1 & 0 & y_4 & 1 & b & u & y_8 & y_9 & y_0 \\ 0 & 0 & 1 & z_4 & 1 & 1 & 1 & 1 & 1 & 1
    \end{pmatrix},
\]
whose columns are the projective coordinates of the points $P_1, \dots, P_9, P_0$. In $\mathcal Y$ the points $P_1,P_2,P_3,P_5$ form a projective frame, and the chart only requires  $p_{234}\ne0$ and $p_{12q}\ne0$ for $q\ne4$. None of these triples lies in $\mathcal L_V$. Hence this chart identifies $\mathcal Y$ with the principal open subset of $\mathbb A^{12}$ on which the $110$ Pl\"ucker coordinates with labels outside $\mathcal L_V$ are nonzero.

\begin{theorem}\label{thm:X310-positive-genus}
    We have $g(\Gr_{3,10},Z_{3,10})=g(\cX_{3,10})\ge1$. In particular, $(\Gr_{3,10},Z_{3,10})$ is not a positive arrangement.
\end{theorem}
\begin{proof}
    We order the triples in $\mathcal L_V$ matching the construction sequence in \cite[\S 2]{Sturmfels1991}
    \[
        I_1,\dots,I_{10} =237, 124, 456, 138, 678, 179, 259, 350, 480, 690.
    \]
    For $0 \leq k \leq 10$, we define
    \[
        R_k := \bigl\{p_{I} = 0 ~ \forall ~ I \in \{I_1 , \dots , I_k\}\bigr\} \cap \left\{p_{I} \neq 0 ~\forall ~ I \in \binom{[10]}{3} \setminus \{I_1 , \dots , I_k\}\right\} \subset \mathcal Y,
    \]
    so $R_0 = \cX_{3,10}$ and $R_{10} = \mathcal R_{\mathcal L_V}$. For $0 \leq k \leq 9$, let $W_k := R_k \cup R_{k+1}$, where $p_{I_{k+1}}$ is allowed to vanish.
    
    Each $W_k$ is isomorphic to a nonempty principal open subset of $\mathbb A^{12-k}$. In particular, the $W_k$ are smooth, affine, and irreducible. Indeed, $W_0$ is the principal open subset of $\mathcal Y\subset\mathbb A^{12}$ where $p_I\neq0$ for $I\in\mathcal L_V\setminus\{237\}$. Next, $W_1$ is the principal open subset of $\mathcal Y \cap \{x_7 = 0\} \subset \mathbb A^{11}$ where $p_I \neq 0$ for $I \in \mathcal L_V \setminus \{237,124\}$, etc. With each iteration, imposing $p_{I_{k+1}} = 0$ solves for one coordinate as a regular function of the others. After nine steps one reaches $W_9$, a principal open subset of $\mathbb A^3_{u,a,b}$ whose points are the configurations
    \[
        \begin{pmatrix}
            1 & 0 & 0 & 1 & 1 & a & 0 & \frac{ua}{u-b} & 1 & x_0 \\
            0 & 1 & 0 & \frac{b-1}{a-1} & 1 & b & u & 0 & u & x_0 \\
            0 & 0 & 1 & 0 & 1 & 1 & 1 & 1 & 1 & 1
        \end{pmatrix},
        \qquad x_0 = -\frac{ua(b-1)}{(a-b)(u-b)},
    \]
    where the denominators $a-1=p_{256}$, $b-u=p_{167}$ and $b-a=p_{356}$ are brackets outside $\mathcal L_V$, hence nonzero on $\mathcal Y$.
	On $W_9$ we have
    \[
        p_{690} = -\frac{F}{(u-b)(b-a)}, \qquad F = u^2a^2 - ua^2 - b^3 + ub^2 + ab^2 - uab + ua - u^2a,
    \]
    and $F$ is the equation of the surface of \cite[(2.1)]{Sink-K3surface} in this chart. So we obtain $R_{10} = W_9 \cap \{F = 0\} = \mathcal R_{\mathcal L_V}$.

     Let $X_k$ be a compactification of $W_k$, and let $Y_k = X_k \setminus R_k$ and $Z_k = X_k \setminus W_k$. Then $Z_k\subset Y_k$, and $Y_k$ has codimension one in $X_k$ because $Y_k\setminus Z_k=R_{k+1}$. As $X_k \setminus Z_k = W_k$ is smooth and affine, \cite[Prop.~3.10]{Brown:2025jjg} gives
    \[
    g(R_k) = g(X_k,Y_k) = g(X_k,Z_k) + g(Y_k,Z_k) = g(W_k) + g(R_{k+1}) \geq g(R_{k+1}).
    \]

    Hence, $g(\cX_{3,10}) = g(R_0) \geq g(R_{10}) = g(\mathcal R_{\mathcal L_V})$. Finally, $\tilde S\setminus\mathcal R_{\mathcal L_V}$ is a curve together with finitely many points, so
    \[
        g(\mathcal R_{\mathcal L_V}) = g(\tilde S,\tilde S\setminus\mathcal R_{\mathcal L_V}) \geq g(\tilde S) = h^0(K_{\tilde S}) = 1
    \]
    by \cite[Prop.~3.9, Cor.~3.11, Prop.~3.2]{Brown:2025jjg}. 
\end{proof}

\begin{remark}\label{rem:X310-form}
	In the chart above, the form
	\[
	\omega=\frac{da\wedge db\wedge du\wedge dx_7\wedge dy_4\wedge dz_4\wedge dx_8\wedge
		dy_8\wedge dx_9\wedge dy_9\wedge dx_0\wedge dy_0}{\prod_{I\in\mathcal L_V}p_I}
	\]
	has simple poles along the ten divisors $\{p_I=0\}$, $I\in\mathcal L_V$, and its
	iterated residue along them is the holomorphic $2$-form $\pm\,du\wedge da/\partial_bF$ of
	the K3 surface $\tilde S$. We verify the residue by computer, through a telescoping factorization
	of $\omega$ along the construction sequence and a Jacobian computation on $W_9$ \cite{KL26data}. 
	We expect $\omega$ to be logarithmic on $\cX_{3,10}$ and to lie in the kernel of the map $R$ in \eqref{eqn:Rmap}; it would then be a nonzero element of $\ker R$, whose dimension is $g(\cX_{3,10})$. 
\end{remark}

\begin{remark}
    In \cite[Theorem 1.2]{Sink-K3surface}, it is shown that there are three more analogous line configurations that also result in realization spaces that are open subsets of K3 surfaces. Alongside the relabelings by the action of $S_{10}$, there is a large number of strata of positive genus in the boundary of $\cX_{3,10}$. We expect therefore that $g(\cX_{3,10})$ is much larger than one, though the strategy proving \zcref{thm:X310-positive-genus} treats one stratum at a time and gives no better bound.
\end{remark}

\begin{remark}\label{rem:positive-genus-regions}
	Positive genus alone does not rule out that a region of $\cX_{3,10}$ is a positive
	geometry in the sense of \cite{Arkani-Hamed:2017tmz}: there are positive geometries whose
	pair has positive genus \cite[\S6]{Koefler:2026vkb}. What \zcref{thm:X310-positive-genus}
	shows is that the framework of Brown and Dupont does not single out canonical forms for the
	regions: the logarithmic forms that the map $R$ of \eqref{eqn:Rmap} sends to the class of a
	region differ by elements of $\ker R$, which has dimension $g(\cX_{3,10})\ge1$. For a single
	region $\sigma$ one may instead take its algebraic boundary $\partial_a\sigma$ as the
	boundary. This can only lower the genus \cite[Prop.~3.10]{Brown:2025jjg}, and if the genus
	becomes zero, $\sigma$ has a canonical form. This happens for the positive Grassmannian
	$\Gr^{\ge0}_{3,10}$, whose algebraic boundary consists of the ten cyclic Pl\"ucker divisors
	(\zcref{cor:positive-grassmannian-rank}). 
\end{remark}

\subsection{More points and higher rank}\label{ss:genus-range}

The genus cannot drop when a point is added, and it is invariant under the duality of Grassmannians. Together with \zcref{thm:X310-positive-genus} this gives positive genus in a large range.

\begin{proposition}\label{prop:genus-monotone}
	Let $2\le k\le n-3$. Then $\Groo_{k,n}$ is isomorphic to the complement of a divisor in $\Groo_{k,n-1}\times\Gm^k$, and
	\[
	g(\cX_{k,n})\ \ge\ g(\cX_{k,n-1}),\qquad \cromb(\cX_{k,n})\ \ge\ \cromb(\cX_{k,n-1}).
	\]
\end{proposition}

\begin{proof}
	Every point of $\Groo_{k,n}$ has a unique representative of the form $(E\,|\,A\,|\,z)$, where $E$ is the $k\times k$ identity matrix, $A$ is a $k\times(n-1-k)$ matrix and $z\in\C^k$. The matrix $(E\,|\,A)$ represents a point of $\Groo_{k,n-1}$, and $z_r=\pm p_{([k]\setminus r)\cup n}\ne0$ for all $r\in[k]$. So $(E\,|\,A\,|\,z)\mapsto\bigl((E\,|\,A),z\bigr)$ is an open immersion of $\Groo_{k,n}$ into $V:=\Groo_{k,n-1}\times\Gm^k$. Its image is the set of pairs for which all maximal minors of $(E\,|\,A\,|\,z)$ that involve the last column are nonzero. So its complement is a divisor.
	
	Let $X$ be a compactification of $V$, let $Z=X\setminus V$, and let $Y\supseteq Z$ be the complement of the image of $\Groo_{k,n}$ in $X$. Then $Y$ has codimension one, so $g(X,Z)\le g(X,Y)$ and $\cromb(X,Z)\le\cromb(X,Y)$ by \cite[Prop.~3.10, Prop.~4.8]{Brown:2025jjg}. Here $g(X,Y)=g(\Groo_{k,n})$ and $g(X,Z)=g(V)=g(\Groo_{k,n-1})$, since a torus factor does not change the genus \cite[Prop.~2.3.10]{Koefler}, and likewise for the combinatorial rank \cite[Prop.~2.2.19]{Koefler}. The claim follows from \zcref{cor:genus_and_cr_descent}.
\end{proof}

The number of ordered $n$-arcs in $\PP^2(\F_q)$, that is, of $n$-tuples of points no three of which are collinear, is $|\PGL_3(\F_q)|\cdot\#\cX_{3,n}(\F_q)$ for $n\ge4$. For $n\le9$ it is a quasi-polynomial in $q$ \cite{Gly88,ISS95}.

\begin{corollary}\label{cor:genus-range}
	\begin{enumerate}
		\item We have $g(\cX_{3,n})\ge1$ for all $n\ge10$. For $n\ge10$ neither the number of points of $\cX_{3,n}$ over $\F_q$ nor the number of $n$-arcs in $\PP^2(\F_q)$ is a quasi-polynomial in $q$.
		\item For $2\le k\le n-2$ we have $g(\cX_{k,n})=g(\cX_{n-k,n})$ and $\cromb(\cX_{k,n})=\cromb(\cX_{n-k,n})$.
		\item We have $g(\cX_{k,n})\ge1$ whenever $3\le k\le n-7$ or $7\le k\le n-3$. In these cases $(\Gr_{k,n},Z_{k,n})$ is not a positive arrangement.
	\end{enumerate}
\end{corollary}

\begin{proof}
	(1) The inequality follows from \zcref{thm:X310-positive-genus,prop:genus-monotone} by induction on $n$. If one of the two counts were a quasi-polynomial, then so would the other, because they
	differ by the polynomial factor $|\PGL_3(\F_q)|$, and $\#\cX_{3,n}(\F_q)$ would be a 
	polynomial in $q$ for all prime powers $q\equiv1\pmod N$. Then the genus would vanish by \zcref{cor:pointcount}.
	
	(2) The isomorphism $\Gr_{k,n}\to\Gr_{n-k,n}$ that sends a subspace to its orthogonal complement maps $p_I$ to $\pm p_{[n]\setminus I}$, so it maps $Z_{k,n}$ onto $Z_{n-k,n}$. Now we apply \zcref{cor:genus_and_cr_descent}.
	
	(3) Let $P$ be the set of pairs $(k,n)$ with $g(\cX_{k,n})\ge1$. By (1) it contains $(3,n)$ for $n\ge10$, hence by (2) also $(k,k+3)$ for $k\ge7$, hence by \zcref{prop:genus-monotone} all $(k,n)$ with $7\le k\le n-3$, and by (2) once more all $(k,n)$ with $3\le k\le n-7$.
\end{proof}

Part (1), together with \zcref{cor:genus_and_cr_descent}, proves \zcref{thm:main310}. For ten points the statement on arcs was announced in \cite[\S1.1]{ISHAM2022102006}, where the count is related to K3 surfaces and elliptic curves. For more than ten points it was conjectured in \cite[\S1.1]{ISHAM2022102006}.

\begin{remark}\label{rem:genus-open}
Together with \zcref{thm:guide-M0n,prop:ff-X3n} and duality, \zcref{cor:genus-range} decides for
every pair $(k,n)$ with $2\le k\le n-2$ whether $(\Gr_{k,n},Z_{k,n})$ is a positive arrangement,
except for the nine pairs with $4\le k\le6$ and $4\le n-k\le6$. By \zcref{prop:genus-monotone} and duality, $g(\cX_{4,8})\le g(\cX_{k,n})\le g(\cX_{6,12})$ for all of them. We do not know whether $g(\cX_{4,8})=0$.
\end{remark}

\section{Recursions for $\Gr_{k,n}$ with smaller boundaries}\label{sec:smaller-boundaries}

By \zcref{cor:genus-range}, the pair $(\Gr_{k,n},Z_{k,n})$ has positive genus whenever $3\le k\le n-7$ or $7\le k\le n-3$, so for most $k$ and $n$ the full Pl\"ucker boundary does not give a positive arrangement. In this section we show that the genus and the combinatorial rank do not change when one adds a point together with a controlled set of Pl\"ucker divisors. The mechanism is an isomorphism
\[
(\Gr_{k,n-1}\setminus Z_S)\times(\C^*)^k\;\cong\;\Gr_{k,n}\setminus Z_{S\cup D(\Gamma)},
\]
which exists whenever $S$ contains a set $N(\Gamma)$ determined by a combinatorial datum $\Gamma$ that we call a nesting (\zcref{lem:nesting-isomorphism}). Since genus and combinatorial rank depend only on the complement and are unchanged by a torus factor, this carries both from $\Gr_{k,n-1}$ to $\Gr_{k,n}$. A second step removes $Z_{N(\Gamma)}$ from the boundary: the pushforward of a logarithmic form has no poles along $Z_{N(\Gamma)}$ (\zcref{sec:form_recursion}), and this gives \zcref{thm:recursion}. As applications we recover the positive Grassmannian (\zcref{cor:positive-grassmannian-rank}) and treat stacked polytopes in \zcref{sec:polytopes}.

\subsection{Relating $\Gr_{k,n-1}$ to $\Gr_{k,n}$}\label{sec:nestings}

A nesting chooses a canonical matrix representative for the first
$n-1$ columns, such that adjoining the $n$-th column becomes a product with a torus. We first define
nestings and their numerator and denominator, then construct the representatives, and then the
isomorphism.

\begin{definition}\label{def:nesting}
	A $k$\emph{-nesting} on a finite set $E$ is a tuple $\Gamma = (I_1, \ldots , I_{k})$ of distinct $(k-1)$-element subsets of $E$ such that for each $r \in [k-1]$ there is a \emph{nesting index} $j_r$ satisfying $j_r \not\in I_r$ and $j_r \in \bigcap_{s > r} I_s$.
\end{definition}

\begin{lemma}\label{lem:nesting-indices}
	For a $k$-nesting $\Gamma = (I_1, \ldots , I_{k})$, the nesting indices $j_1,\dots,j_{k-1}$ are unique. Moreover, there is a unique index $j_k$ with $j_k \not\in I_k$ and $j_k \in I_{k-1}$.
\end{lemma}

\begin{proof}
	The nesting indices are distinct, as for $s > r$ one has $j_r \in I_s$ and $j_s \not\in I_s$. Since the set of nesting indices is fixed as $I_k = \{j_1,\ldots,j_{k-1}\}$, a different choice of nesting indices would be related to the original by a permutation $(j'_1 , \dots , j'_{k-1}) = (j_{\sigma_1} , \dots , j_{\sigma_{k-1}})$. If $\sigma$ is not the identity, then there is an $r$ with $\sigma_r < r$, and $j_r' = j_{\sigma_r}\not\in I_r$ contradicts $j_{\sigma_r} \in \bigcap_{s > \sigma_r} I_s$.
	
	For the second statement, note that $I_k = \{j_1,\dots,j_{k-1}\}$ and that $I_{k-1} \supset \{j_1, \dots, j_{k-2}\}$ does not contain $j_{k-1}$. Consequently, $I_{k-1} \setminus I_k$ consists of a single element $j_k$.
\end{proof}

\begin{example}\label{exa:nesting}
	The tuple $\Gamma = (234,167,159,158)$ is a $4$-nesting on $[9]$ with nesting indices  $(j_1,j_2,j_3,j_4) = (1,5,8,9)$. 
	
	The tuple $\Gamma = (45,15,12)$ on $[5]$ is a $3$-nesting. Here $j_1 = 1$ lies in $15 \cap 12$ but not in $45$, $j_2 = 2$ lies in $12$ but not in $15$, and $j_3 = 5$ is the element of $15$ not in $12$.
\end{example}

To compare $\Gr_{k,n-1}$ with $\Gr_{k,n}$ we attach two sets of Plücker coordinates to a nesting.

\begin{definition}
	Let $\Gamma$ be a $k$-nesting on a ground set $E$. Its \emph{nesting numerator} is the set $N(\Gamma) = \{I_r \cup j_r : r \in [k-1]\}$. For $e \not\in E$ its \emph{nesting denominator} is  $D_e(\Gamma) = \{I_r \cup e : r \in [k]\}$.
\end{definition}

The set $J := \{j_1,\ldots,j_k\}$ equals $I_{k-1} \cup j_{k-1} = I_k \cup j_k$, so $J \in N(\Gamma)$. Typically, the ground set will be $[n-1]$, and the nesting denominator is $D(\Gamma) := D_n(\Gamma)$. Then $N(\Gamma)$ and $D(\Gamma)$ are both subsets of $\binom{[n]}{k}$, and define divisors $Z_{N(\Gamma)}$ and $Z_{D(\Gamma)}$ on $\Gr_{k,n}$. For the $4$-nesting of \zcref{exa:nesting} and $n=10$ we get $N(\Gamma) = \{1234,1567,1589\}$. For the $3$-nesting $\Gamma$ and $n = 6$ we have $N(\Gamma) = \{145,125\}$ and $D(\Gamma) = \{456,156,126\}$.

\begin{remark}
    In the language of Grassmann\textendash Cayley algebras \cite[Section 3.3]{InvariantTheoryAlgorithms_Sturmfels}, the conditions of \zcref{def:nesting} are sufficient for the meet $I_1 \wedge \cdots \wedge I_k$ to be a monomial in a bracket expansion. This monomial is $\prod_{r \in [k-1]} p_{I_r \cup j_r}$, so the nesting numerator lists the appearing sets. We leave the exploration of this relationship to Grassmann\textendash Cayley algebras to future work.
\end{remark}

It is clear from the definition that $D(\Gamma)$ does not depend on the ordering of the sets in $\Gamma$. We show the same is true for $N(\Gamma)$, provided the ordering still leaves a $k$-nesting. Consequently, the results that follow which depend only on $N(\Gamma)$ and $D(\Gamma)$ will be independent of the choice of ordering.

\begin{lemma}\label{lem:reorder-first}
    Let $\Gamma$ be a $k$-nesting such that for some $r \in [k]$ one has $j_r \not\in I_r$ and $j_r \in \bigcap_{s \neq r} I_s$. Then, $\Gamma' := (I_r , I_1 , \cdots , \hat I_r , \cdots , I_k)$ is a $k$-nesting with $N(\Gamma') = N(\Gamma)$.
\end{lemma}
\begin{proof}
    One can verify that the nesting indices
    \[(j'_1,j'_2, \cdots , j'_{k-1}) := (j_r,j_1 , \cdots, \hat j_r , \cdots , j_{k-1})\]
    satisfy the conditions of \zcref{def:nesting}. Furthermore,
    \[
        N(\Gamma') = \{I_r \cup j_r , I_1 \cup j_1, \cdots \} = N(\Gamma).
    \]
\end{proof}

\begin{lemma}\label{lem:reorder}
	Let $\Gamma = (I_1 , \dots, I_k)$ be a $k$-nesting on a ground set $E$, and let $\sigma \in S_k$ be a permutation such that $\sigma \Gamma = (I_{\sigma_1} , \dots , I_{\sigma_k})$ is again a $k$-nesting. Then $N(\sigma \Gamma) = N(\Gamma)$.
\end{lemma}

\begin{proof}
    We prove this by induction on $k$, starting with the base case $k=1$, where both nesting numerators are empty. Suppose the lemma is true for $k' < k$:
    
    Since $\sigma\Gamma$ is a $k$-nesting, $j_{\sigma_1} \not\in I_{\sigma_1}$ and $j_{\sigma_1} \in \bigcap_{s \neq \sigma_1} I_s$, so we may apply the preceding lemma to find that $\Gamma' = (I_{\sigma_1} , I_1 , \cdots , \hat I_{\sigma_1} , \cdots , I_k) $ is a $k$-nesting with $N(\Gamma') = N(\Gamma)$. Then, stripping $I_{\sigma_1}$ and $j_{\sigma_1}$ from $\Gamma'$ and $\sigma \Gamma$, one finds $(k-1)$-nestings on the ground set $E \setminus j_{\sigma_1}$,
    \[
        \begin{aligned}
            \widetilde{(\Gamma')} &:= (I_1 \setminus j_{\sigma_1} , \cdots , \widehat{I_{\sigma_1}} , \cdots , I_{k} \setminus j_{\sigma_1}), \\
            \widetilde{(\sigma \Gamma)} &:= (I_{\sigma_2} \setminus j_{\sigma_1} , \cdots , I_{\sigma_{k}} \setminus j_{\sigma_1}).
        \end{aligned}
    \]
    The original nesting numerators can be recovered by
    \[
        \begin{aligned}
            N(\Gamma') &= \bigcup_{I \in N\bigl(\widetilde{(\Gamma')}\bigr)} \{I \cup j_{\sigma_1} \} \cup \{I_{\sigma_1} \cup j_{\sigma_1}\}, \\
            N(\sigma \Gamma) &= \bigcup_{I \in N\bigl(\widetilde{(\sigma \Gamma)}\bigr)} \{I \cup j_{\sigma_1} \} \cup \{I_{\sigma_1} \cup j_{\sigma_1}\}.
        \end{aligned}
    \]
    Since $\widetilde{(\Gamma')}$ and $\widetilde{(\sigma \Gamma)}$ are $(k-1)$-nestings, the inductive hypothesis applies, and
    \[
        N\bigl(\widetilde{(\Gamma')}\bigr) = N\bigl(\widetilde{(\sigma \Gamma)}\bigr) \implies N(\Gamma') = N(\sigma \Gamma).
    \]
    Recalling that $N(\Gamma) = N(\Gamma')$, one also has $N(\Gamma) = N(\sigma \Gamma)$ as desired.
\end{proof}

Now, we proceed to use these nestings to construct explicit matrix representatives.

\begin{definition}
	Let $\Gamma = (I_1 , \dots , I_{k})$ be a $k$-nesting on a subset of $[m]$. A \emph{$\Gamma$-representative} of $C \in \Gr_{k,m} \setminus Z_{N(\Gamma)}$ is a matrix $C_\Gamma$ representing $C$ such that for each $r \in [k]$ one has
	\[
	(C_\Gamma)_{r\,j_r}=1\qquad\text{and}\qquad(C_\Gamma)_{r\,i}=0\quad\text{for all } i\in I_r .
	\]
\end{definition}

\begin{lemma}
    Let $\Gamma$ be a $k$-nesting on a subset of $[m]$. Every $C \in \Gr_{k,m} \setminus Z_{N(\Gamma)}$ has a unique $\Gamma$-representative $C_\Gamma$. Its entries are
    \[
    (C_\Gamma)_{r\,i}\ =\ \pm\,\frac{p_{I_r\cup i}(C)}{p_{I_r\cup j_r}(C)},
    \]
    where the sign depends on the position of $i$ and $j_r$ relative to $I_r$. In particular they are regular functions on $\Gr_{k,m}\setminus Z_{N(\Gamma)}$. 
\end{lemma}
\begin{proof}
    Fix some matrix $C$ representing the point. Every other representative is of the form $MC$ for $M \in GL_k$. We seek to construct an invertible matrix $M$ such that $MC$ is a $\Gamma$-representative. To do so, for $r \in [k]$ we will make use of the $k \times k$ submatrices $C_r$ with columns $I_r \cup j_r$. Since $C \in \Gr_{k,m} \setminus Z_{N(\Gamma)}$, the determinants $\det C_r = \pm p_{I_r \cup j_r}(C)$ are nonzero, since $I_r\cup j_r\in N(\Gamma)$ for $r<k$ and $I_k\cup j_k=J\in N(\Gamma)$.
    
    The $r$th row of $MC$ is $M_r C$ where $M_r$ is the $r$th row of $M$. Requiring that $MC$ is a $\Gamma$-representative amounts to requiring that $M_r C_r = e_{j_r}$ for each $r \in [k]$, where $e_{j_r}$ is the unit vector in $\C^k$ with $1$ in the position of $j_r$ within $I_r\cup j_r$. This has the unique solution $M_r = e_{j_r} (C_r)^{-1}$ for each $r$.

    Checking that $M$ is invertible is equivalent to verifying that $MC$ is full-rank. The submatrix of $MC$ with columns $I_k \cup j_k = J$, ordered as $j_1,\dots,j_k$, is upper unitriangular, so indeed $MC$ has rank $k$ and $M$ is invertible.
\end{proof}

Geometrically, the columns of $C$ indexed by $I_r$ span a hyperplane $H_r\subset\C^k$, since $p_{I_r\cup j_r}(C)\neq0$. Row $r$ of $C_\Gamma$ lists the values on the columns of $C$ of the linear form $\lambda_r=M_r$ that vanishes on $H_r$ and takes the value $1$ on the column $j_r$. The forms $\lambda_1,\dots,\lambda_k$ are linearly independent, so the hyperplanes $H_1,\dots,H_k$ bound a simplex, and $C_\Gamma$ is $C$ written in the coordinates $\lambda_1,\dots,\lambda_k$ of this simplex.

\begin{example}
	For $\Gamma = (45,15,12)$ a point of $\Gr_{3,5} \setminus Z_{\{145, 125\}}$ has the $\Gamma$-representative
	\[
	C_{\Gamma}\ =\ \begin{pmatrix} 1 & * & * & 0 & 0\\ 0 & 1 & * & * & 0\\ 0 & 0 & * & * & 1\end{pmatrix}.
	\]
\end{example}

\begin{example}\label{eg:internal-nesting}\label{eg:consecutive-nesting}\label{eg:stacked-nesting}
	For $\Gamma=(234,167,159,158)$ the first nine columns of the $\Gamma$-representative are
	\[
	\begin{pmatrix}
		1 & 0 & 0 & 0 & * & * & * & * & * \\ 0 & * & * & * & 1 & 0 & 0 & * & * \\ 0 & * & * & * & 0 & *  & * & 1 & 0 \\ 0 & * & * & * & 0 & * & * & 0 & 1
	\end{pmatrix}.
	\]
	For $\Gamma=(234,134,124,123)$ the $\Gamma$-representative has the identity matrix in the first four columns. More generally, every ordering of $\binom{[k]}{k-1}$ is a $k$-nesting on $[k]$, and its $\Gamma$-representative has a permutation matrix in the first $k$ columns. For $\Gamma=(123,126,156,456)$ and for $\Gamma=(234,135,126,123)$ the first six columns are
	\[
	\begin{pmatrix}
		0 & 0 & 0 & * & * & 1 \\ 0 & 0 & * & * & 1 & 0 \\ 0 & * & * &1 & 0 & 0 \\ 1 & * & * & 0 & 0 & 0
	\end{pmatrix}
	\qquad\text{and}\qquad
	\begin{pmatrix}
		1 & 0 & 0 & 0 & * & * \\ 0 & 1 & 0 & * & 0 & * \\ 0 & 0 & 1 & * & * & 0 \\ 0 & 0 & 0 & * & * & 1
	\end{pmatrix}.
	\]
\end{example}

Now we add a column. We write $(M\,|\,z)$ for the matrix $M$ with the column $z$ appended.

\begin{lemma}\label{lem:nesting-isomorphism}
    Let $\Gamma$ be a $k$-nesting on $[n-1]$. The map
    \[
        \begin{aligned}
            f_\Gamma : (\Gr_{k,n-1} \setminus Z_{N(\Gamma)}) \times (\mathbb C^*)^k &\rightarrow \Gr_{k,n} \setminus Z_{N(\Gamma) \cup D(\Gamma)}, \\
            (C,z) &\mapsto (C_\Gamma | z),
        \end{aligned}
    \]
    is an isomorphism of varieties. This map restricts to an isomorphism
    \[
        f_\Gamma : (\Gr_{k,n-1} \setminus Z_{S}) \times (\mathbb C^*)^k \rightarrow \Gr_{k,n} \setminus Z_{S \cup D(\Gamma)},
    \]
    for any $S \subset \binom{[n-1]}{k}$ containing $N(\Gamma)$. On its image,
    \begin{equation}
    	z_r\ =\ \pm\,\frac{p_{I_r\cup n}}{p_{I_r\cup j_r}}\qquad\text{for } r\in[k].
    \end{equation}
\end{lemma}

\begin{proof}
    We use $\Delta_{A,B}(M)$ to denote the minor of the matrix $M$ with rows $A$ and columns $B$.

    Since the entries of $C_\Gamma$ are regular functions on $\Gr_{k,n-1} \setminus Z_{N(\Gamma)}$, with denominators $p_{I}$ for $I \in N(\Gamma)$, the map $f_\Gamma$ is regular away from $Z_{N(\Gamma)}$. Since $C_\Gamma$ is full-rank, so is $(C_\Gamma | z)$. Since $\Delta_{[k],I}(C_\Gamma | z) = \Delta_{[k],I}(C_\Gamma)$ for $I \in \binom{[n-1]}{k}$, the nonvanishing of minors in the domain implies the nonvanishing of the corresponding minors in the image. Consequently, the image avoids $Z_{S}$. Since
    \begin{equation}\label{eqn:nesting-isomorphism_proof}
        \Delta_{[k],I_r \cup n}(C_\Gamma | z) = \epsilon \Delta_{[k] \setminus r , I_r}(C_\Gamma) \cdot z_r = \epsilon \Delta_{[k],I_r \cup j_r}(C_\Gamma) \cdot z_r,
    \end{equation}
    for some suitable sign $\epsilon$, the image also avoids $Z_{D(\Gamma)}$.

	Conversely,  for $C' \in \Gr_{k,n} \setminus Z_{S \cup D(\Gamma)}$, there is a $\Gamma$-representative $C'_\Gamma$ whose entries are regular functions of Pl\"ucker coordinates away from $Z_{N(\Gamma)}$. The first $n-1$ columns of $C'_\Gamma$ are a full-rank matrix, and an analogous computation to equation \eqref{eqn:nesting-isomorphism_proof} shows that the last column has nonzero entries. So, taking the first $n-1$ columns to be an element of $\Gr_{k,n-1} \setminus Z_{S}$ and the last column to be a vector in $(\mathbb C^*)^k$ is a well-defined regular map $g_\Gamma$.
	The maps $f_\Gamma$ and $g_\Gamma$ are inverse to each other, as both pass through the unique $\Gamma$-representatives.
\end{proof}

\begin{figure}
	\[
	\begin{tikzpicture}[scale=1.5]
		\def\r{0.05}
		\clip (-0.8,-1.1) rectangle (3,1.9);
		\coordinate (F) at (0.2,0.9);
		\coordinate (B) at (2,1.5);
		\coordinate (RTL) at (0,1.3);
		\coordinate (RTR) at (1.8,1.9);
		\coordinate (RBL) at (1.2,-1.1);
		\coordinate (RBR) at (3,-0.5);
		\coordinate (BTL) at (0.4,1.3);
		\coordinate (BTR) at (2.2,1.9);
		\coordinate (BBL) at (-0.8,-1.1);
		\coordinate (BBR) at (1.0,-0.5);
		\fill[red!50,opacity=0.3] (RTL) -- (RTR) -- (B) -- (F) -- cycle;
		\fill[white,opacity=0.7] (BTL) -- (BTR) -- (B) -- (F) -- cycle;
		\fill[blue!50,opacity=0.3] (BTL) -- (BTR) -- (B) -- (F) -- cycle;
		\fill[blue!50,opacity=0.3] (BBL) -- (BBR) -- (B) -- (F) -- cycle;
		\fill[white,opacity=0.7] (RBL) -- (RBR) -- (B) -- (F) -- cycle;
		\fill[red!50,opacity=0.3] (RBL) -- (RBR) -- (B) -- (F) -- cycle;
		\draw[dotted] (F) -- (B);
		\coordinate (1) at (0.2,0);
		\coordinate (2) at (0,-0.5);
		\coordinate (3) at (0.5,-0.4);
		\coordinate (4) at (1.5,-0.3);
		\coordinate (5) at (2.2,-0.6);
		\coordinate (6) at (2,0);
		\path[name path=L] ($(F)!-10cm!(B)$) -- ($(B)!-10cm!(F)$);
		\path[name path=L12] ($(1)!-10cm!(2)$) -- ($(2)!-10cm!(1)$);
		\path [name intersections={of=L12 and L, by=FF}];
		\draw[dotted] (1) -- (2);
		\draw (1) -- (FF);
		\path[name path=L56] ($(5)!-10cm!(6)$) -- ($(6)!-10cm!(5)$);
		\path [name intersections={of=L56 and L, by=BB}];
		\draw[dotted] (5) -- (6);
		\draw (6) -- (BB);
		\draw (FF) -- (BB);
		\draw (1) -- (6);
		\draw (FF) -- (6);
		\draw[dashed] (1) -- (BB);
		\draw[fill] (FF) circle (\r) node[above left] {$L_{12} \cap H_{456}$};
		\draw[fill] (BB) circle (\r) node[above left] {$L_{56} \cap H_{123}$};
		\draw[fill] (1) circle (\r) node[left] {1};
		\draw[fill] (2) circle (\r) node[left] {2};
		\draw[fill] (3) circle (\r) node[below] {3};
		\draw[fill] (4) circle (\r) node[below] {4};
		\draw[fill] (5) circle (\r) node[below left] {5};
		\draw[fill] (6) circle (\r) node[right] {6};
	\end{tikzpicture}
	\hspace{+3em}
	\begin{tikzpicture}[scale=1.5]
		\def\r{0.05}
		\coordinate (1) at (0,0);
		\coordinate (2) at (2,0);
		\coordinate (3) at (1.2,-0.8);
		\coordinate (4) at (2.1,-1.1);
		\coordinate (5) at (0.3,-1.1);
		\coordinate (6) at (0.9,-0.4);
		\coordinate (A) at (0.9,1.1);
		\draw[fill] (6) circle (\r) node[above] {6};
		\fill[gray!20,opacity=0.8] (1) -- (3) -- (2) -- (1);
		\draw[dashed] (1) -- (2);
		\draw (1) -- (3) -- (2);
		\draw (1) -- (A);
		\draw (2) -- (A);
		\draw (3) -- (A);
		\draw[fill] (1) circle (\r) node[left] {1};
		\draw[fill] (2) circle (\r) node[right] {2};
		\draw[fill] (3) circle (\r) node[below] {3};
		\draw[fill] (A) circle (\r) node[above] {$H_{234} \cap H_{135} \cap H_{126}$};
		\draw[fill] (4) circle (\r) node[below right] {4};
		\draw[fill] (5) circle (\r) node[below left] {5};
	\end{tikzpicture}
	\]
	\caption{Configurations of six points in $\PP^3$. Here $L_{ab}$ is the line through $a,b$ and $H_{abc}$ the plane through $a,b,c$. Left: the planes $H_{123},H_{126},H_{156},H_{456}$ of the $4$-nesting $\Gamma_L=(123,126,156,456)$ bound the tetrahedron with vertices $1$, $6$, $L_{12}\cap H_{456}$ and $L_{56}\cap H_{123}$; the shaded planes are $H_{123}$ and $H_{456}$. Right: for $\Gamma_R=(234,135,126,123)$ the vertices are $1,2,3$ and $H_{234}\cap H_{135}\cap H_{126}$. A seventh point lies outside $Z_{D(\Gamma)}$ if and only if it avoids the four facet planes, and $z$ records its coordinates relative to the tetrahedron.}\label{fig:k4-extensions}
\end{figure}
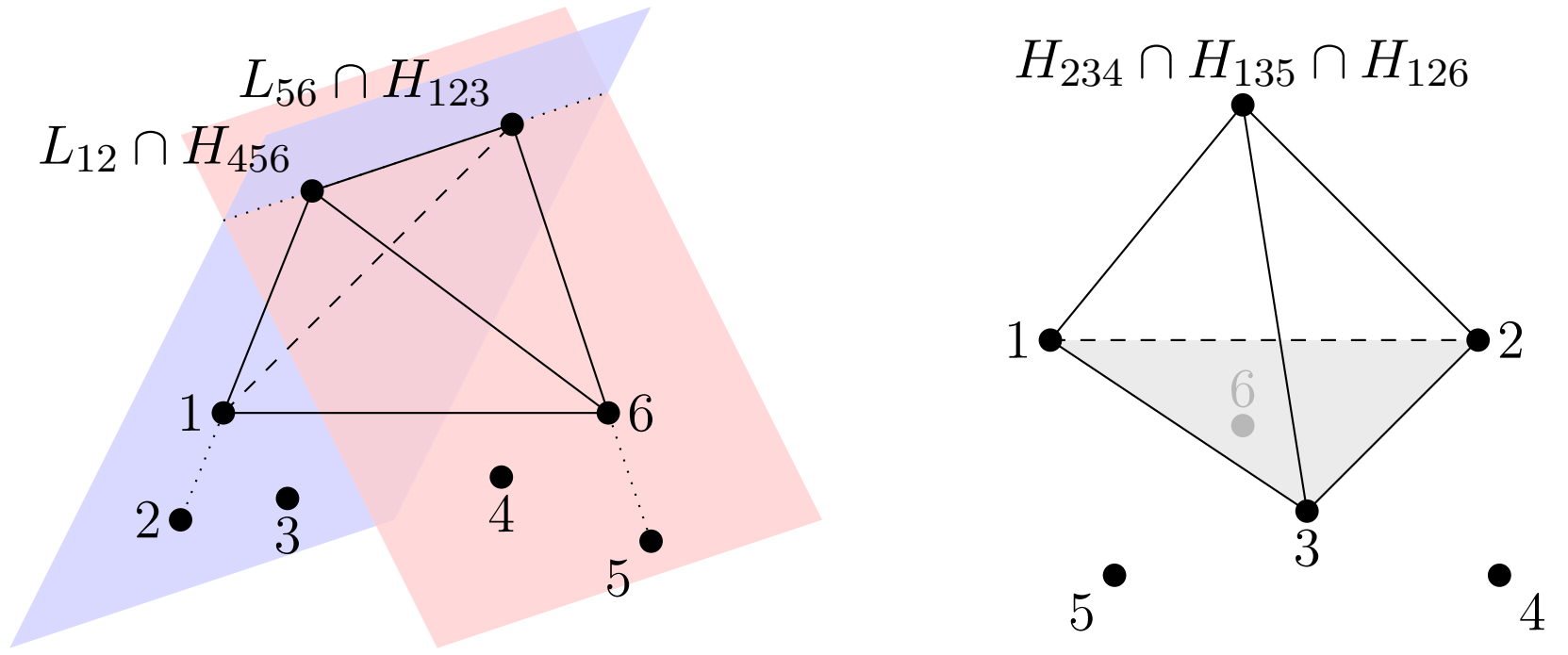

In these terms, the last column $z$ consists of the values $\lambda_r$ at the new point $n$ in the coordinates of the simplex bounded by $H_1,\dots,H_k$. It can be chosen freely, and its entries are nonzero exactly when the point avoids the hyperplanes $H_1,\dots,H_k$, which is the condition $D(\Gamma)$. \zcref{fig:k3-extensions,fig:k4-extensions} show this simplex for $k=3$ and $k=4$.

\begin{example}\label{eg:running-isomorphism}
	For $\Gamma_R = (45,15,12)$ the isomorphism of \zcref{lem:nesting-isomorphism} with $S=N(\Gamma_R)$ reads
	\[
	\bigl(\Gr_{3,5}\setminus Z_{\{145,125\}}\bigr)\times\Gm^3\ \xrightarrow{\ \sim\ }\ \Gr_{3,6}\setminus Z_{\{145,125,456,156,126\}},\qquad(C,z)\mapsto\begin{pmatrix} 1 & * & * & 0 & 0 & z_1\\ 0 & 1 & * & * & 0 & z_2\\ 0 & 0 & * & * & 1 & z_3\end{pmatrix}.
	\]
	The dimensions are $6+3=9$. A direct computation with this matrix gives the inverse:
	\[
	z_1=\frac{p_{456}}{p_{145}},\qquad z_2=-\frac{p_{156}}{p_{125}},\qquad z_3=\frac{p_{126}}{p_{125}} .
	\]
	So the sixth point may be anywhere off the three lines through $4,5$, through $1,5$ and through $1,2$, and $z$ records where it is.
\end{example}

\begin{example}\label{eg:k3-figure}
	\zcref{fig:k3-extensions} pictures two $3$-nestings. The left image shows $\Gamma_L=(24,25,45)$, with nesting indices $(5,4,2)$, $N(\Gamma_L)=\{245\}$ and $D(\Gamma_L)=\{246,256,456\}$; the right image shows $\Gamma_R$. The three blue lines form a triangle, that is, they are in general position, as long as the five points lie outside $Z_{N(\Gamma)}$. For $\Gamma_L$ this just says that $2,4,5$ are not collinear. For $\Gamma_R$ it says that the lines through $1,2$ and through $1,5$ are distinct and that the line through $4,5$ misses their intersection point $1$. A sixth point gives a configuration outside $Z_{D(\Gamma)}$ if and only if it avoids the blue lines. \zcref{fig:k4-extensions} shows the analogous picture for $k=4$, where the hyperplanes span the facets of a tetrahedron.
\end{example}

\begin{figure}
	\[
	\begin{tikzpicture}
		\def\r{0.05}
		\coordinate (1) at (-90+2*72:0.8);
		\coordinate (2) at (-90+72:1);
		\coordinate (3) at (-90:1);
		\coordinate (4) at (-90-72:1);
		\coordinate (5) at (-90-72*2:0.8);
		\draw[dashed] (-2,-2) rectangle (2,2);
		\clip (-2,-2) rectangle (2,2);
		\draw[fill] (1) circle (\r) node[above right] {1};
		\draw[fill] (2) circle (\r) node[above right] {2};
		\draw[fill] (3) circle (\r) node[below] {3};
		\draw[fill] (4) circle (\r) node[above left] {4};
		\draw[fill] (5) circle (\r) node[shift={(-0.1,0.3)}] {5};
		\draw[blue, thick] ($(2)!-10cm!(5)$) -- ($(5)!-10cm!(2)$);
		\draw[blue, thick] ($(2)!-10cm!(4)$) -- ($(4)!-10cm!(2)$);
		\draw[blue, thick] ($(4)!-10cm!(5)$) -- ($(5)!-10cm!(4)$);
	\end{tikzpicture}
	\hspace{+3em}
	\begin{tikzpicture}
		\def\r{0.05}
		\coordinate (1) at (-90+2*72:0.8);
		\coordinate (2) at (-90+72:1);
		\coordinate (3) at (-90:1);
		\coordinate (4) at (-90-72:1);
		\coordinate (5) at (-90-72*2:0.8);
		\draw[dashed] (-2,-2) rectangle (2,2);
		\clip (-2,-2) rectangle (2,2);
		\draw[fill] (1) circle (\r) node[above right] {1};
		\draw[fill] (2) circle (\r) node[right] {2};
		\draw[fill] (3) circle (\r) node[below] {3};
		\draw[fill] (4) circle (\r) node[left] {4};
		\draw[fill] (5) circle (\r) node[above left] {5};
		\draw[blue, thick] ($(1)!-10cm!(2)$) -- ($(2)!-10cm!(1)$);
		\draw[blue, thick] ($(4)!-10cm!(5)$) -- ($(5)!-10cm!(4)$);
		\draw[blue, thick] ($(1)!-10cm!(5)$) -- ($(5)!-10cm!(1)$);
	\end{tikzpicture}
	\]
	\caption{Configurations in $\cX_{3,5}$, to be imagined as shadows of points in $\Gr_{3,5}$. The blue lines are given by the $3$-nestings $\Gamma_L=(24,25,45)$ (left) and $\Gamma_R=(45,15,12)$ (right).}\label{fig:k3-extensions}
\end{figure}
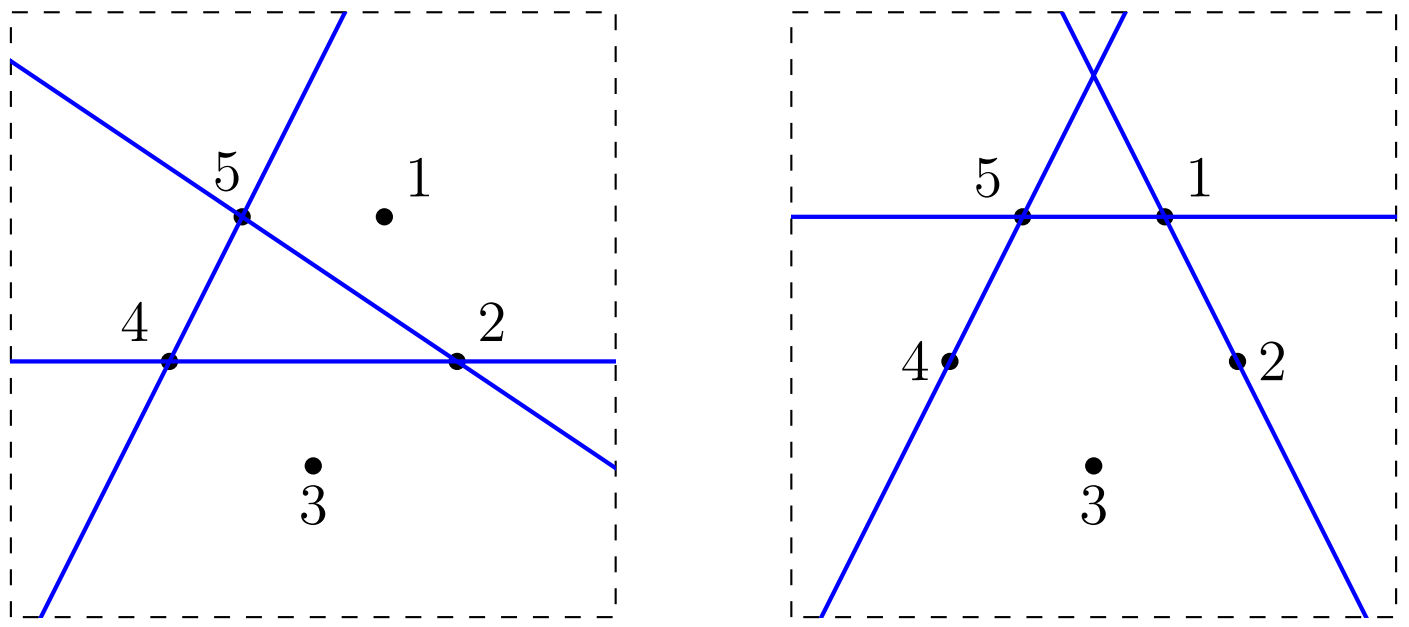

We now specialize to a nesting relevant to the positive Grassmannian.

\begin{example}\label{eg:positive_grassmannian_isomorphism}
    Let $S^+_{k,n}\subset\binom{[n]}{k}$ be the set of cyclic intervals $\{i,i+1,\dots,i+k-1\}$, $i\in[n]$, taken modulo $n$, so that $Z_{S^+_{k,n}}=\partial_a\Gr^{\ge0}_{k,n}$. For $n-1\ge k$, consider the $k$-nesting $\Gamma=(I_1,\dots,I_k)$ on $[n-1]$ with
    \[
        I_r=\{n-k+r,\dots,n-1,1,\dots,r-1\},\qquad r\in[k],
    \]
    the $k$ cyclic intervals of $[n-1]$ of size $k-1$ that become adjacent to $n$ in $[n]$. Its nesting indices are $j_r=r$ for $r<k$ and $j_k=n-1$. The denominator $D(\Gamma)=S^+_{k,n}\setminus S^+_{k,n-1}$ consists of the $k$ cyclic intervals of $[n]$ containing $n$. For $n\ge k+2$ the numerator $N(\Gamma)=S^+_{k,n-1}\setminus S^+_{k,n}$ consists of the $k-1$ cyclic intervals of $[n-1]$ containing both $n-1$ and $1$, which stop being cyclic intervals once $n$ is inserted; for $n=k+1$ all $I_r\cup j_r$ equal $[k]$. For $k=3$ and $n=6$ this is the nesting $\Gamma_R=(45,15,12)$ of \zcref{fig:k3-extensions}. In both cases $S^+_{k,n-1}\cup D(\Gamma)=S^+_{k,n}\cup N(\Gamma)$, so \zcref{lem:nesting-isomorphism} with $S=S^+_{k,n-1}$ gives an isomorphism
    \[
        f_\Gamma\colon(\Gr_{k,n-1}\setminus Z_{S^+_{k,n-1}})\times(\mathbb C^*)^k\ \longrightarrow\ \Gr_{k,n}\setminus Z_{S^+_{k,n}\cup N(\Gamma)} .
    \]
\end{example}

\begin{remark}
	Gorsky and Scroggin construct isomorphisms $\Phi_a\colon U_a\to\Pi^\circ_{k,n-a+1}\times\Pi^\circ_{k,a+k-1}$, where $U_a$ is an open subset of the open positroid variety $\Pi^\circ_{k,n}=\Gr_{k,n}\setminus Z_{S^+_{k,n}}$ \cite[Thm.~1.1]{GorskyScroggin}. For $n\ge k+2$, the isomorphism $f_\Gamma$ of \zcref{eg:positive_grassmannian_isomorphism} is the inverse of $\Phi_2$, with $\Pi^\circ_{k,k+1}\cong(\C^*)^k$, after the cyclic relabeling $i\mapsto i+1$ of $[n]$, which preserves $S^+_{k,n}$ and moves the column $n$ to the first position \cite[proof of Lemma~3.9]{GorskyScroggin}.
\end{remark}

\subsection{Recursion for differential forms}\label{sec:form_recursion}

Pushing forward along $f_\Gamma$ turns logarithmic forms on $\Gr_{k,n-1}$ into logarithmic forms on $\Gr_{k,n}$, with an explicit formula (\zcref{lem:form-recursion}). We first recall the standard notation $h\,dC$ for top forms on a Grassmannian \cite[\S 7.4 and App.~C.2, eq.~(C.9)]{Arkani-Hamed:2017tmz}, see also \cite[\S 13]{Lam:2014omj}, since we need it in the charts given by $\Gamma$-representatives.

For $I \in \binom{[n]}{k}$ let $U_I := \Gr_{k,n} \setminus Z_I$. For $C \in U_I$, let $C_I := (C_{[I]})^{-1} C$ be the matrix representative with the identity matrix in the columns $I$. The $k(n-k)$ entries $c_{ab}$ of $C_I$ outside of columns $I$ are coordinates on $U_I$, and we write $dC_I := \bigwedge_{a \in [k], b \in [n] \setminus I} dc_{ab}$ in a fixed order. On $U_I \cap U_J$, one can check that
\[
    dC_J = \pm \left(\frac{p_I}{p_J}\right)^n dC_I.
\]
This reflects $K_{\Gr_{k,n}} \cong \mathcal O(-n)$. For $h$ a rational function in the Plücker coordinates, homogeneous of degree $-n$, the function $h p_I^n$ is a well-defined rational function on $\Gr_{k,n}$. The forms $h p_I^n dC_I$ agree up to sign for all $I$, and we denote this form by $h(C) dC$. For $k = 1$, the differential $dC$ is the form $\sum_i (-1)^i x_i dx_1 \wedge \dots \widehat{dx_i} \dots \wedge dx_n$ on $\mathbb P^{n-1}$; in the physics literature $dC$ is often denoted by $d^{k\times n} C / \mathrm{vol}\, GL(k)$.

Let $\Gamma$ be a $k$-nesting on $[n-1]$. Since $J = \{j_1,\dots,j_k\} = I_{k-1} \cup j_{k-1} \in N(\Gamma)$, the domain and target of $f_\Gamma$ are contained in the charts $U_J \subset \Gr_{k,n-1}$ and $U_J' \subset \Gr_{k,n}$. Let $M := (C_\Gamma)_{[J]}$ be the submatrix with columns $J$, so that $(C_\Gamma)_J = M^{-1} C_\Gamma$ and
\[
    (f_\Gamma(C, z))_J = (C_J | M^{-1} z).
\]
In the affine coordinates of the two charts, the map $f_\Gamma$ is the map $(C_J, z) \mapsto (C_J | M^{-1}z) =: C'_J$. Since $M$ depends only on $C$ and is upper unitriangular,
the Jacobian matrix of this map is block triangular, with the identity on the $C_J$-block and $M^{-1}$ on the $\vec z$-block, so its determinant is $\det M^{-1} = 1$ and one finds
$(f_\Gamma)_*(dC_J \wedge dz_1 \wedge \cdots \wedge dz_k) = \pm dC'_J$. Consequently, analyzing the pushforward of a top-form on $(\Gr_{k,n-1} \setminus Z_{N(\Gamma)}) \times (\mathbb C^*)^k$ along $f_\Gamma$ simplifies as
\[
    (f_\Gamma)_* \left( h(C) dC \wedge \frac{dz_1 \wedge \cdots \wedge dz_k}{z_1 \cdots z_k}\right) = \pm p_J^{-1} (f_\Gamma)_* \left( \frac{h(C)}{z_1 \cdots z_k} \right) dC'.
\]
The factor of $p_J^{-1}$ comes from the different weights of $dC$ and $dC'$.

\begin{lemma}\label{lem:form-recursion}
    Let $\Gamma$ be a $k$-nesting on $[n-1]$ and $S \subseteq \binom{[n-1]}{k}$, and let $f_\Gamma$ be the corresponding isomorphism as in \zcref{lem:nesting-isomorphism}. Let $h(C) dC$ be a form in $\Omega^{k(n-1-k)}_{\log}(\Gr_{k,n-1} \setminus Z_{N(\Gamma) \cup S})$. Then
    \[
        \omega = h(C) dC \wedge \frac{dz_1 \wedge \cdots \wedge dz_k}{z_1 \cdots z_k}
    \]
    is a logarithmic form on $(\Gr_{k,n-1} \setminus Z_{N(\Gamma) \cup S}) \times (\mathbb C^*)^k$, and its pushforward along $f_\Gamma$ is the logarithmic form
    \[
        (f_\Gamma)_* \omega = \pm \frac{\prod_{r \in [k-1]} p_{I_r \cup j_r}}{\prod_{I \in D(\Gamma)} p_I} h(C') dC'.
    \]
\end{lemma}
\begin{proof}
    Since $h(C)$ depends only on Pl\"ucker coordinates with indices in $[n-1]$, $(f_\Gamma)_* h(C) = h(C')$. By \eqref{eqn:nesting-isomorphism_proof}, $z_r = \pm p_{I_r \cup n}/p_{I_r \cup j_r}$ for all
    $r \in [k]$, as functions of $C'$. Since $I_k \cup j_k = J$,
    \[
        \frac{1}{p_J\, z_1 \cdots z_k}
        = \pm \frac{\prod_{r \in [k]} p_{I_r \cup j_r}}{p_J \prod_{I \in D(\Gamma)} p_I}
        = \pm \frac{\prod_{r \in [k-1]} p_{I_r \cup j_r}}{\prod_{I \in D(\Gamma)} p_I},
    \]
    and the claim follows from the display above. The form is logarithmic because $\omega$ is,
    and $f_\Gamma$ is an isomorphism.
\end{proof}

The numerator and denominator in the expression for $(f_\Gamma)_*\omega$ are the reason for the names nesting numerator\footnote{Note that $\prod_{r \in [k-1]} p_{I_r \cup j_r} \neq \prod_{I \in N(\Gamma)} p_I$ in general as some terms may need to appear with higher multiplicity. Choosing $\Gamma=(234,134,124,123)$ from \zcref{eg:internal-nesting}, one has $\prod_{r \in [k-1]} p_{I_r \cup j_r} = (p_{1234})^3$.} and nesting denominator.

\begin{corollary}\label{cor:no-poles-NG}
    Let $n-1 > k$ and $N(\Gamma)\subseteq S$. Then $(f_\Gamma)_*\omega$ extends to a logarithmic form on $\Gr_{k,n} \setminus Z_{(S \setminus N(\Gamma)) \cup D(\Gamma)}$, and it restricts to a logarithmic form on $\Gr_{k,n}\setminus Z_{S'}$ for every $S'\supseteq(S\setminus N(\Gamma))\cup D(\Gamma)$.
\end{corollary}
\begin{proof}
    Put $V:=\Gr_{k,n}\setminus Z_{(S\setminus N(\Gamma))\cup D(\Gamma)}$. By \zcref{lem:form-recursion}, $(f_\Gamma)_*\omega$ is logarithmic on $V\setminus Z_{N(\Gamma)}=\Gr_{k,n}\setminus Z_{S\cup D(\Gamma)}$. Since $n-1>k$, each $Z_{\{I\}}\subset\Gr_{k,n-1}$ is a hypersurface, so $h\,dC$ has at most simple poles along $Z_{\{I\}}$ for $I\in S$ \cite[Prop.~1.13]{Brown:2025jjg}. The numerator $\prod_{r\in[k-1]}p_{I_r\cup j_r}$ contains $p_I$ for every $I\in N(\Gamma)$, so $(f_\Gamma)_*\omega$ is holomorphic on $V$.

    Let $\bar V$ be a compactification of $V$ with normal crossing boundary $E$. By \cite[Prop.~1.13]{Brown:2025jjg}, $(f_\Gamma)_*\omega$ has at most simple poles along $E\cup\overline{Z_{N(\Gamma)}\cap V}$. Since it is holomorphic on $V$, Hartogs' theorem shows that it has at most simple poles along $E$. As $E$ is normal crossing, the form is therefore logarithmic on $V$ \cite[Prop.~1.13]{Brown:2025jjg}.

    The restriction to $\Gr_{k,n}\setminus Z_{S'}\subseteq V$ is logarithmic by functoriality of logarithmic forms \cite[\S1.5.1]{Brown:2025jjg}, applied to a log resolution of $\bar V\setminus(\Gr_{k,n}\setminus Z_{S'})$.
\end{proof}

The next example shows that a region can have $Z_I$ with $I \in N(\Gamma)$ in its boundary while its canonical form has no pole there.

\begin{example}
Let $S = \{12,23,34,41\}$ as for $\partial_a\Gr^{\geq 0}_{2,4}$, and $\Gamma = (\{1\},\{4\})$. Then, $N(\Gamma) = \{14\}$ and $D(\Gamma) = \{15,45\}$. The region $\sigma = \{p_{12}, p_{13},p_{14},p_{15} \geq 0, ~ p_{42},p_{43},p_{45} \geq 0, ~ p_{32} \geq 0\}$ is a chamber in $\Gr_{2,5} \setminus Z_{S \cup N(\Gamma) \cup D(\Gamma)}$, with $Z_{\{14\}}$ as part of its boundary. One may decompose
    \[
        \sigma = (\sigma \cap \{p_{25} \geq 0 , p_{35} \geq 0\}) \cup (\sigma \cap \{p_{25} \leq 0 , p_{35} \geq 0\}) \cup (\sigma \cap \{p_{25} \leq 0 , p_{35} \leq 0\}),
    \]
each of which is a single oriented matroid stratum in $\Gr_{2,5} \setminus Z_{2,5}$. Computing the canonical forms, one finds
    \[
        \begin{aligned}
            \omega_\sigma &= \frac{dC}{p_{14} p_{43} p_{32} p_{25} p_{51}} + \frac{dC}{p_{14} p_{43} p_{35} p_{52} p_{21}} + \frac{dC}{p_{14} p_{45} p_{53} p_{32} p_{21}} \\
            &= -\frac{dC}{p_{15} p_{54} p_{43} p_{32} p_{21}},
        \end{aligned}
    \]
    which does not have a pole along $Z_{\{14\}}$, although $Z_{\{14\}}$ is part of the boundary of $\sigma$.
\end{example}

\begin{example}\label{eg:positive_grassmannian_form}
    The canonical form on the positive Grassmannian can be found by induction. Using the $k$-nesting $\Gamma$ of \zcref{eg:positive_grassmannian_isomorphism}, one finds
    \[
\begin{aligned}
	\omega_{\Gr_{k,n}^{\geq 0}} &= (f_\Gamma)_* \bigl( \omega_{\Gr_{k,n-1}^{\geq 0}} \wedge \omega_{(\mathbb C^*)^k} \bigr)\\
	&= \pm \frac{p_{(n-k+1) \dots (n-1) 1} \cdots p_{(n-1) 1 \dots (k-1)}}{p_{(n-k+1) \dots n} \cdots p_{n 1 \dots (k-1)}} \cdot \frac{1}{p_{1 \dots k} \cdots p_{(n-1) 1 \dots (k-1)}}\, dC \\
	&= \pm \frac{dC}{p_{1 \dots k}\, p_{2 \dots (k+1)} \cdots p_{n 1 \dots (k-1)}},
\end{aligned}
    \]
    where the numerator of the ratio cancels the spurious poles of the original form.
\end{example}

\subsection{Genus and combinatorial rank by induction}\label{sec:induction}

We collect the results of this section in a single statement.

\begin{theorem}\label{thm:recursion}
	Let $2\le k<n-1$, let $\Gamma$ be a $k$-nesting on $[n-1]$, and let $S\subseteq\binom{[n-1]}{k}$ with $N(\Gamma)\subseteq S$.
	\begin{enumerate}
		\item The genus and combinatorial rank satisfy
		\[
			g(\Gr_{k,n},Z_{S\cup D(\Gamma)})=g(\Gr_{k,n-1},Z_S),\qquad
			\cromb(\Gr_{k,n},Z_{S\cup D(\Gamma)})=\cromb(\Gr_{k,n-1},Z_S).
		\]
		In particular, $(\Gr_{k,n},Z_{S\cup D(\Gamma)})$ is a positive arrangement if and only if $(\Gr_{k,n-1},Z_S)$ is.
		\item The linear map
		\[
			\Omega^{k(n-1-k)}_{\log}\bigl(\Gr_{k,n-1}\setminus Z_S\bigr)\longrightarrow\Omega^{k(n-k)}_{\log}\bigl(\Gr_{k,n}\setminus Z_{S\cup D(\Gamma)}\bigr),\qquad
			\omega\longmapsto(f_\Gamma)_*\bigl(\omega\wedge\frac{dz_1 \wedge \dots \wedge dz_k}{z_1 \cdots z_k}\bigr),
		\]
		given explicitly by \zcref{lem:form-recursion}, is an isomorphism.
		\item Suppose $(\Gr_{k,n-1},Z_S)$ is a positive arrangement, and let $S'$ satisfy $(S\setminus N(\Gamma))\cup D(\Gamma)\subseteq S'\subseteq S\cup D(\Gamma)$. Then $(\Gr_{k,n},Z_{S'})$ is a positive arrangement with
		\[
			\cromb(\Gr_{k,n},Z_{S'})=\cromb(\Gr_{k,n-1},Z_S),
		\]
		and restriction identifies the logarithmic forms on $\Gr_{k,n}\setminus Z_{S'}$ with the forms in (2).
	\end{enumerate}
\end{theorem}

Part (3) of \zcref{thm:recursion} proves \zcref{thm:mainrecursion}. The case $S'=(S\setminus N(\Gamma))\cup D(\Gamma)$ is what we use below for the positive Grassmannian and, in \zcref{sec:polytopes}, for stacked polytopes.

\begin{proof}
	Throughout, $f_\Gamma\colon(\Gr_{k,n-1}\setminus Z_S)\times(\C^*)^k\to\Gr_{k,n}\setminus Z_{S\cup D(\Gamma)}$ is the isomorphism of \zcref{lem:nesting-isomorphism}, in the form for a set $S\supseteq N(\Gamma)$. The genus and the combinatorial rank of a pair $(X,Y)$, with $X$ compact and $X\setminus Y$ smooth, depend only on $X\setminus Y$ \cite[(16)]{Brown:2025jjg}.

	(1) The variety $\Gr_{k,n-1}\setminus Z_S$ is smooth and quasi-projective, so the torus factor changes neither the genus nor the combinatorial rank \cite[Prop.~2.3.10, Prop.~2.2.19]{Koefler}. For the last claim, $\Gr_{k,n}$ is smooth and defined over $\R$ with real points. Each $Z_{\{I\}}$ is, after permuting columns, a Schubert divisor, so it is irreducible and defined over $\R$. Its singular locus is where the columns in $I$ have rank at most $k-2$, so any real matrix whose columns in $I$ have rank exactly $k-1$ gives a smooth real point. Hence a pair $(\Gr_{k,n},Z_{S''})$ with $S''\neq\emptyset$ is a positive arrangement if and only if its genus vanishes.

	(2) By \zcref{lem:form-recursion} the map takes values in logarithmic forms. It is injective, since wedging with $\dlog z_1\wedge\cdots\wedge\dlog z_k$ is injective on top forms pulled back from the first factor, and $f_\Gamma$ is an isomorphism. For $X$ smooth and compact of dimension $m$, $\dim\Omega^m_{\log}(X\setminus Y)=g(X,Y)+\cromb(X,Y)$ by \eqref{eqn:Rmap}, see \cite[(16), (21), (24), Lemma~2.3]{Brown:2025jjg}. By (1) the source and the target therefore have the same dimension.

	(3) Put $W:=\Gr_{k,n}\setminus Z_{S\cup D(\Gamma)}$ and $W':=\Gr_{k,n}\setminus Z_{S'}$. Since $Z_{S'}\subseteq Z_{S\cup D(\Gamma)}$, \cite[Prop.~3.10]{Brown:2025jjg} and (1) give $g(\Gr_{k,n},Z_{S'})\le g(\Gr_{k,n},Z_{S\cup D(\Gamma)})=0$. As $S'\supseteq D(\Gamma)$ is nonempty, $(\Gr_{k,n},Z_{S'})$ is a positive arrangement by (1). Restriction $\Omega_{\log}(W')\to\Omega_{\log}(W)$ is well defined by functoriality of logarithmic forms \cite[\S1.5.1]{Brown:2025jjg}. It is injective because $W$ is dense in $W'$, and surjective by (2) and \zcref{cor:no-poles-NG}. Both pairs have genus zero, so the combinatorial rank is the dimension of the space of logarithmic forms, and
	\[
		\cromb(\Gr_{k,n},Z_{S'})=\dim\Omega_{\log}(W')=\dim\Omega_{\log}(W)=\cromb(\Gr_{k,n-1},Z_S). \qedhere
	\]
\end{proof}

\begin{corollary}[Positive Grassmannian revisited]\label{cor:positive-grassmannian-rank}
    Let $S_{k,n}^+=\bigl\{(i,i+1,\cdots,i+k-1):i\in[n]\bigr\}$ be the set of cyclically consecutive $k$-subsets, so that $Z_{S^+_{k,n}}=\partial_a\Gr_{k,n}^{\ge0}$ is the algebraic boundary of the
    nonnegative Grassmannian. Then $(\Gr_{k,n},Z_{S_{k,n}^+})$ is a positive arrangement for all $k<n$, in particular
    $g\bigl(\Gr_{k,n},Z_{S_{k,n}^+}\bigr) = 0$, and
    \[
        \cromb_\R\bigl(\Gr_{k,n},Z_{S_{k,n}^+}\bigr)=\cromb\bigl(\Gr_{k,n},Z_{S_{k,n}^+}\bigr)=1 .
    \]
\end{corollary}

\begin{proof}
    For $n=k+1$ the Pl\"ucker coordinates identify $(\Gr_{k,k+1},Z_{S^+_{k,k+1}})$ with $\PP^k$ and its coordinate hyperplanes, a toric pair with $\cromb_\R=\cromb=1$ \cite[\S5.3]{Brown:2025jjg}, and likewise $(\Gr_{1,n},Z_{S^+_{1,n}})$ with $\PP^{n-1}$.

    Now let $k\ge2$ and $n\ge k+2$, and suppose the claim holds for $\Gr_{k,n-1}$. Take the $k$-nesting $\Gamma$ of \zcref{eg:positive_grassmannian_isomorphism}, $S=S^+_{k,n-1}$, and $S'=S^+_{k,n}$. Then $N(\Gamma)=S^+_{k,n-1}\setminus S^+_{k,n}\subseteq S$ and $D(\Gamma)=S^+_{k,n}\setminus S^+_{k,n-1}$, so $S'=(S\setminus N(\Gamma))\cup D(\Gamma)$. Part (3) of \zcref{thm:recursion} shows that $(\Gr_{k,n},Z_{S^+_{k,n}})$ is a positive arrangement with $\cromb=\cromb(\Gr_{k,n-1},Z_{S^+_{k,n-1}})=1$. Finally, $\cromb_\R = 1$ because the region $\Gr^{\ge0}_{k,n}$ has the nonzero canonical form of \zcref{eg:positive_grassmannian_form}.
\end{proof}

\begin{remark}
    These isomorphisms allow one to construct positive arrangements in $\Gr_{3,n}$ for $n \geq 10$, such as $(\Gr_{3,10} , Z_{\binom{[9]}{3}} \cup Z_{\{230,150,140\}})$ by starting with $(\Gr_{3,9},Z_{3,9})$ and using the isomorphism associated to $\Gamma =(23,15,14)$. We note that the argument proving \zcref{thm:X310-positive-genus} fails for such arrangements, as having the K3 surface in a boundary stratum would require $D(\Gamma) = \{0a_1a_2,0a_3a_4,0a_5a_6\}$ for distinct $a_i$, which could not come from a nesting.
\end{remark}

\section{Realization spaces of simplicial polytopes}\label{sec:polytopes}

So far our main examples of $Z_S$ were the full
Plücker boundary and the boundary of the positive Grassmannian. In this section our choices of $S$ come from convex geometry. To a polytope $K$ with $n$ labeled vertices one attaches a semialgebraic subset $\Gr_K^{>0}$ of the real Grassmannian, the Grassmannian model of its realization space \cite{GMW-Polytope}. It is cut out by sign conditions on Plücker coordinates. If $K$ is an $n$-gon in the plane with vertices labeled in cyclic order, $\Gr_K^{>0}$ is the positive Grassmannian $\Gr_{3,n}^{>0}$. We show that for a simplicial $3$-polytope $K$ the realization space $\Gr_K^{>0}$ is a connected component of $\Gr_{4,n}(\R) \setminus Z_K(\R)$, where the \emph{polytopal divisor} $Z_K$ consists of the Plücker divisors along which two adjacent facets become coplanar (\zcref{lem:polytope-boundary}). Once the pair $(\Gr_{4,n},Z_K)$ has genus zero, it is a positive arrangement, the closure of $\Gr_K^{>0}$ is one of its regions, and this region has a unique canonical form $\omega_K$. We show that for small polytopes genus zero follows from our main theorems (\zcref{lem:small-polytope}). Stacking a vertex on a facet preserves genus zero by the recursion of \zcref{sec:smaller-boundaries}. This gives an explicit recursion for the canonical form. For stacked polytopes the combinatorial rank and the real combinatorial rank are one.

\subsection{Realization spaces of polytopes}

A \emph{combinatorial polytope} is the face lattice of a convex polytope. 
Let $K$ be a combinatorial $(k-1)$-dimensional polytope with $n$ labeled vertices, with facets $F_1,\ldots,F_l$. We identify a facet $F_i \subset [n]$ with the set of vertices it contains. A \emph{realization} of $K$ is a tuple $(q_1,\ldots,q_n)$ of points of $\R^{k-1}$ such that $\operatorname{conv}(q_1,\ldots,q_n)$ has the face lattice $K$. We are interested in the Grassmannian model of the space of all realizations. Given a realization we form the $k \times n$ matrix $Q$ whose $i$-th column is the homogenized vertex $\hat{q_i} = (q_i,1) \in \R^k$. Since the polytope is full-dimensional, the matrix $Q$ has rank $k$, so its row span is a point of $\Gr_{k,n}(\R)$, and
\[
p_I(Q) = \det(\hat{q_i}_{i \in I}) = (k-1)! \cdot \operatorname{vol}(\operatorname{conv}(q_i: i \in I))
\]
is the signed volume of the simplex spanned by the vertices in $I$. The combinatorics of the polytope $K$ is reflected in these signs. If $J = \{j_1 < \dots < j_{k-1}\}$ is a set of vertices that affinely spans a facet $F$, the linear form $\ell_J(x) = \det(\hat{q_{j_1}},\dots,\hat{q_{j_{k-1}}},x)$ on $\R^k$ vanishes exactly on the hyperplane of $F$ and by convexity it has the same sign at $\hat{q}_v$ for every vertex $v \notin F$. So for $I = J \cup \{v\}$ the sign of $p_I(Q)$ records on which side of the facet $F$ the vertex $v$ lies. 

\begin{definition}[\cite{GMW-Polytope}]
    A set $I \in \binom{[n]}{k}$ is a \emph{facet extension} of $K$ if $I = J \cup \{v\}$ where $J$ consists of $k-1$ vertices that affinely span a facet $F$ of $K$ and $v \notin F$. Denote the set of facet extensions of $K$ by $\overline{\mathcal F}_K$.
\end{definition}

\begin{definition}[\cite{GMW-Polytope}]
	Fix a realization $Q_0$ of $K$ and set $\Delta_I = \sign p_I(Q_0)$ for $I \in \overline{\mathcal F}_K$. 
	The \emph{realization space model} of $K$ is $\Gr^{>0}_K \subset \Gr_{k,n}(\R)$ defined by vanishing of Pl\"ucker coordinates for which all indices are contained in a single facet, and fixed signs for Pl\"ucker coordinates for which all but one index is contained in a single facet,
	\[
		\Gr^{>0}_K = \left\{Q \in \Gr_{k,n}(\R) : p_I(Q) = 0 ~ \forall ~ I \in \binom{F_i}{k}\text{ and all $i \in [l]$}\right\} \cap \left\{\Delta_I p_I(Q) > 0 ~ \forall ~ I \in \overline{\mathcal F}_K\right\}.
	\]
\end{definition}

The idea behind this realization space model is that the vertices of each facet lie on a common hyperplane, and every other vertex lies strictly on the inner side of it; this is convex position with the prescribed face lattice.

Recall that the torus $(\R^*)^n$ acts on $\Gr_{k,n}(\R)$ by rescaling columns. Let $T_{>0} = \R^n_{>0}/\R_{>0}$ be its positive part which acts on $\Gr_K^{>0}$. Two realizations lie in the same $T_{>0}$-orbit if and only if they differ by a projective transformation. 

For polygons the realization space is the positive Grassmannian.

\begin{lemma}
	Let $K$ be the $n$-gon with vertices labeled $1,\dots,n$ in cyclic order. Then $\Gr^{>0}_K = \Gr^{>0}_{3,n}$.
\end{lemma}
\begin{proof}
    The facets are the edges $\{i,i+1\}$, indices modulo $n$, and the facet extensions are the sets $\{i,i+1,j\}$ that contain two cyclically consecutive indices. Orienting the boundary counterclockwise, all signs $\Delta_I$ are $+1$. 
    Hence the realization space model is
    \[
        \Gr_K^{>0} = \bigl\{ p_{i,i+1,j} > 0 ~ \forall ~ i \in [n], j \in [n] \setminus \{i,i+1\}\bigr\},
    \]
    so clearly $\Gr_{3,n}^{>0} \subset \Gr_K^{>0}$.
    For the reverse inclusion, we are left to show that the positivity of the Pl\"ucker coordinates with two consecutive indices is sufficient to determine the positivity of all Pl\"ucker coordinates. Let $Q \in \Gr_K^{>0}$ and $r < s < t$. We show $p_{rst} > 0$, by induction on $t-r$. If $\{r,s,t\}$ contains two cyclically consecutive indices, there is nothing to show, and this covers the case $t-r = 2$. Otherwise $r+1 < s < t-1$, and the three-term Plücker relation for the four indices $r+1 < s < t-1 < t$, with $r$ fixed, reads 
    \[
        p_{rst} p_{r(r+1)(t-1)} = p_{r(r+1)s}p_{r(t-1)t} + p_{r(r+1)t} p_{rs(t-1)}.
    \]
    The coordinates $p_{r(r+1)(t-1)}, p_{r(r+1)s}, p_{r(t-1)t}, p_{r(r+1)t}$ contain consecutive indices, and $p_{rs(t-1)} > 0$ by induction. Hence $p_{rst} > 0$.
\end{proof}

\subsection{Simplicial polytopes and the polytopal divisor}
From now on $K$ is simplicial: every facet has exactly $k-1$ vertices. In this case, the realization space model is cut out by strict inequalities alone and thus $\Gr_K^{>0}$ is an open subset of $\Gr_{k,n}(\R)$ of dimension $k(n-k).$ For two adjacent facets $F,F'$ their union $F \cup F'$ has $k$ elements, and it is the facet extension of $F$ by the vertex of $F'$ not on $F$.

\begin{definition}\label{def:polytopal-divisor}
	For a simplicial polytope $K$, the \emph{polytopal divisor} $Z_K$ is a divisor in $\Gr_{k,n}$ defined by $Z_K := Z_{S_K}$ for
    \[
		S_K := \{F \cup F' \text{ for all adjacent facets $F$ and $F'$}\} \subseteq \overline{\mathcal F}_K.
	\]
\end{definition}

For a configuration of points, $p_{F \cup F'} = 0$ exactly when the adjacent facets $F$ and $F'$ become coplanar: the polytope flattens along the ridge $F \cap F'$. For a $3$-polytope, $Z_K$ is the locus where two triangles sharing an edge become coplanar. Since $S_K \subseteq \overline{\mathcal F}_K,$ $Z_K(\R)$ is disjoint from $\Gr_K^{>0}$. For a semialgebraic set $A \subseteq \Gr_{k,n}(\R)$ we write $\partial_a A$ for its \emph{algebraic boundary}, the Zariski closure in $\Gr_{k,n}$ of its Euclidean boundary. We write $\Gr_K^{\geq 0}$ for the Euclidean closure of $\Gr_K^{>0}$ in $\Gr_{k,n}(\R)$.

Now we show that the only way to fall out of a realization space is when two adjacent facets become coplanar, and all the other facet-extension inequalities are redundant. In other words, a vertex cannot reach the plane of a non-adjacent facet, and two vertices cannot collide, before some edge has become flat. 

\begin{lemma}\label{lem:polytope-boundary}
	Let $K$ be a simplicial $(k-1)$-polytope with $n$ vertices. Then the Euclidean boundary of $\Gr_K^{\geq 0}$ lies in $Z_K(\R)$, and $\partial_a \Gr^{\geq 0}_K \subseteq Z_K$. In particular, if $\Gr_K^{>0}$ is connected, then it is a connected component of $\Gr_{k,n}(\R) \setminus Z_K(\R)$.
\end{lemma}

\begin{proof}
	The set $\Gr_K^{>0}$ is open and disjoint from $Z_K(\R)$. We show that every point $Q \in \Gr_K^{\geq 0} \setminus Z_K(\R)$ lies in $\Gr_K^{>0}$. Then the Euclidean boundary $\Gr_K^{\geq 0} \setminus \Gr_K^{>0}$ lies in $Z_K(\R)$, so does its Zariski closure, and $\Gr_K^{>0}$ is open and closed in $\Gr_{k,n}(\R) \setminus Z_K(\R)$, hence a connected component if it is connected.
	
	Fix $J \in S_K$ and represent the points of $\{p_J \neq 0\} \supseteq \Gr_K^{\geq 0} \setminus Z_K(\R)$ by matrices with columns $q_1,\dots,q_n \in \R^k$ and $\Delta_J\, p_J > 0$, so that the signs of all $p_I$ are well defined. For a facet $G = \{g_1 < \dots < g_{k-1}\}$ let $\ell_G = \pm\det(q_{g_1},\dots,q_{g_{k-1}},\,\cdot\,)$, a linear form on $\R^k$, with the sign such that $\ell_G(q_u) = \Delta_{G \cup u}\, p_{G \cup u}$ for all vertices $u \notin G$; it exists because in $Q_0$ all vertices off $G$ lie on the same side of the hyperplane of $G$. If $Q \in \Gr_K^{\geq 0} \setminus Z_K(\R)$, then by continuity
	\begin{equation}\label{eq:weak}
		\ell_G(q_u) \geq 0 \ \text{ for all facets $G$ and all vertices $u$,} \qquad \ell_G(q_u) > 0 \ \text{ if } G \cup u \in S_K.
	\end{equation}
	That is, every vertex lies weakly on the inner side of every facet hyperplane, and strictly so if it lies on an adjacent facet.
	
	\emph{Step 1: if one facet-extension coordinate vanishes at $Q$, then so does a second one.} Suppose $p_{F \cup w} = 0$ for a facet $F$ and a vertex $w \notin F$. For $a \in F$ let $F_a = (F \setminus a) \cup \tilde a$ be the facet adjacent to $F$ along the ridge $F \setminus a$. Since $F \cup \tilde a \in S_K$, we have $\ell_{F_a}(q_a) > 0$ and $w \neq \tilde a$. In particular the columns $q_f$, $f \in F$, are linearly independent, and $p_{F \cup w} = 0$ gives
	\[
	q_w = \sum_{f \in F} \alpha_f\, q_f .
	\]
	Evaluating $\ell_{F_a}$, which vanishes on $q_f$ for $f \in F \setminus a$, gives $\alpha_a\, \ell_{F_a}(q_a) = \ell_{F_a}(q_w) \geq 0$, so $\alpha_a \geq 0$ for all $a \in F$: the vertex $w$ lies in the simplex $F$. As $q_w \neq 0$, because $w$ lies in some $I \in S_K$, some $\alpha_a$ is positive. There is a facet $G$ with $w \in G$ and $a \notin G$, since the facets through $w$ cannot all contain $a$: the link of $w$ in $\partial K$ is a sphere, not a cone. Then
	\[
	0 = \ell_G(q_w) = \sum_{f \in F} \alpha_f\, \ell_G(q_f)
	\]
	is a sum of nonnegative terms by \eqref{eq:weak}, so $\ell_G(q_a) = 0$, that is, $p_{G \cup a} = 0$. This is a second coordinate: if $G \cup a = F \cup w$, then $G = (F \setminus a) \cup w$ is the facet $F_a$, so $w = \tilde a$.
	
	\emph{Step 2: codimension two is not enough to form a boundary.} Let $W \subset \Gr_{k,n}$ be the locus where at least two of the coordinates $p_I$, $I \in \overline{\mathcal F}_K$, vanish. It has codimension two, since distinct Pl\"ucker coordinates cut out distinct irreducible divisors, so $W(\R)$ is a closed semialgebraic set of dimension at most $k(n-k)-2$ \cite[\S 2.8]{BCR}. As $Q$ was arbitrary, Step~1 says
	\begin{equation}\label{eq:bdry}
		\Gr_K^{\geq 0} \setminus Z_K(\R) \;\subseteq\; \Gr_K^{>0} \cup W(\R).
	\end{equation}
	Now suppose, for contradiction, that some $Q \in \Gr_K^{\geq 0} \setminus Z_K(\R)$ is not in $\Gr_K^{>0}$, and let $B$ be a connected open neighborhood of $Q$ in $\Gr_{k,n}(\R) \setminus Z_K(\R)$. Then $B \setminus W(\R)$ is connected \cite[Thm.~IV~4, Cor.~1]{HurewiczWallman} and dense in $B$. The set $\Gr_K^{>0} \cap B$ is a nonempty open subset of $B \setminus W(\R)$. It is also closed in $B \setminus W(\R)$: a limit point $P \in B \setminus W(\R)$ of $\Gr_K^{>0}$ lies in $\Gr_K^{\geq 0} \setminus Z_K(\R)$, because $B$ does not meet $Z_K(\R)$, and then $P \in \Gr_K^{>0}$ by \eqref{eq:bdry}, because $P \notin W(\R)$. A connected space has no nonempty open and closed subset other than itself, so $\Gr_K^{>0} \cap B$ is all of $B \setminus W(\R)$. As this set is dense in $B$, the weak inequalities $\Delta_I\, p_I \geq 0$, $I \in \overline{\mathcal F}_K$, hold on all of $B$ by continuity.
	
	But since $Q \notin \Gr_K^{>0}$, some $p_{F \cup w}$ with $F \cup w \in \overline{\mathcal F}_K$ vanishes at $Q$, and it changes sign near $Q$: the columns of $Q$ indexed by $F$ are independent, so we can move the column $q_w$ to $q_w + \varepsilon d$ with $\ell_F(d) < 0$, which for small $\varepsilon > 0$ gives a point of $B$ with $\Delta_J\, p_J > 0$ and $\Delta_{F \cup w}\, p_{F \cup w} = \varepsilon\, \ell_F(d) < 0$. This is a contradiction. Hence every point of $\Gr_K^{\geq 0} \setminus Z_K(\R)$ lies in $\Gr_K^{>0}$, which proves the lemma.
\end{proof}

\begin{remark}\label{rem:connected-realization}
	By \cite[Thm.~3.7]{GMW-Polytope}, every point of $\Gr_K^{>0}$ is represented, up to a positive rescaling of the columns, by a homogenized realization of $K$. Hence $\Gr_K^{>0}$ is connected whenever the realization space of $K$ is. This is the case for simplicial $3$-polytopes, whose realization spaces are contractible \cite[Thm.~13.3.3]{RichterGebert}, and for stacked polytopes of every dimension: a realization of the polytope obtained from $K$ by stacking a vertex on a facet $F$ is a realization of $K$ together with a point beyond $F$ and beneath all other facets, and these points form an open convex set, so connectedness propagates from the simplex. In dimension four and higher the realization space of a simplicial polytope can be disconnected, so the hypothesis of the last statement of \zcref{lem:polytope-boundary} is not automatic.
\end{remark}

\subsection{Small 3-polytopes}

For a simplicial $3$-polytope $K$ the pair $(\Gr_{4,n}, Z_K)$ satisfies the reality conditions of \zcref{def:posarr}. By \zcref{lem:polytope-boundary}, $\Gr_K^{>0}$ is a connected component of $\Gr_{4,n}(\R) \setminus Z_K(\R)$ and it is orientable, since it lies in the affine chart $\{p_I \neq 0\} \cong \R^{4(n-4)}$ for any $I \in \overline{\mathcal F}_K$. Hence, as soon as $g(\Gr_{4,n}, Z_K) = 0$, the pair is a positive arrangement and the closure of $\Gr_K^{>0}$ is one of its regions, with a canonical form $\omega_K = \omega(\overline{\Gr_K^{>0}})$.  

\begin{lemma}\label{lem:small-polytope}
    Let $K$ be a simplicial $3$-polytope with $n \leq 7$ vertices. Then $(\Gr_{4,n}, Z_K)$ is a positive arrangement, and the closure of $\Gr^{>0}_K$ is one of its regions.
\end{lemma}
\begin{proof}
    It suffices to show that the genus vanishes. Since $Z_K \subseteq Z_{4,n}$, we have $g(\Gr_{4,n}, Z_K) \leq g(\Gr_{4,n}, Z_{4,n})$ and the duality $\Gr_{4,n} \cong \Gr_{n-4,n}$ identifies $(\Gr_{4,n}, Z_{4,n})$ with $(\Gr_{n-4,n}, Z_{n-4,n})$. For $n=4$ vertices, $\Gr^{>0}_K = \Gr_{4,4}$ is a single point. For $n=5$ vertices, it is $\PP^4$ with five hyperplanes in general position whose complement is the torus $(\C^*)^4$. For $n=6$ vertices, it is $(\Gr_{2,6}, Z_{2,6})$ which has genus zero by \zcref{thm:guide-M0n}, since $\cX_{2,6} \cong M_{0,6}$ is the complement of a hyperplane arrangement. For $n=7$ vertices, it is $(\Gr_{3,7}, Z_{3,7})$, which has genus zero by \zcref{thm:posarr-3n}.
\end{proof}

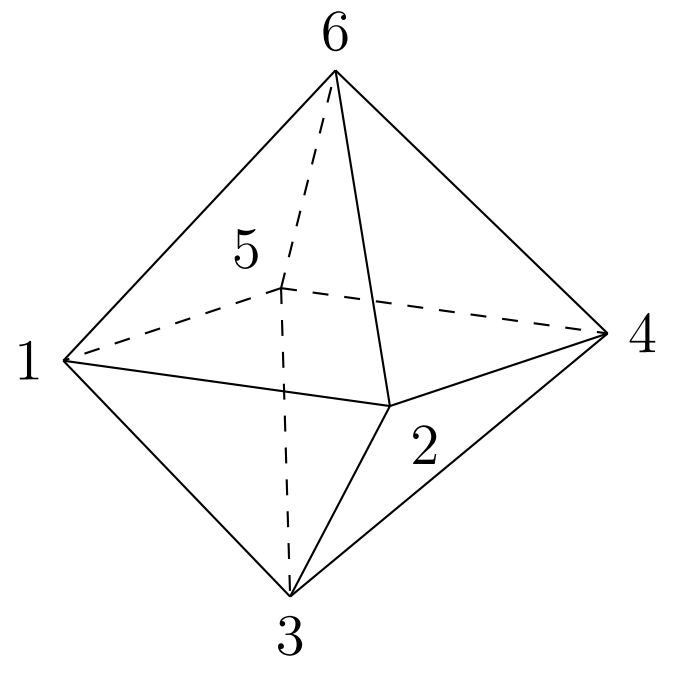
\begin{figure}
    \begin{tikzpicture}[scale=0.6]
          \coordinate (1) at (-3,0);
          \coordinate (2) at (0.6,-0.5);
          \coordinate (3) at (-0.5,-2.6);
          \coordinate (4) at (3,0.3);
          \coordinate (5) at (-0.6,0.8);   %
          \coordinate (6) at (0,3.2);
        
          \draw[dashed] (5) -- (1);
          \draw[dashed] (5) -- (4);
          \draw[dashed] (5) -- (6);
          \draw[dashed] (5) -- (3);
        
          \draw (1) -- (6);
          \draw (6) -- (4);
          \draw (6) -- (2);
          \draw (1) -- (2);
          \draw (1) -- (3);
          \draw (3) -- (2);
          \draw (2) -- (4);
          \draw (3) -- (4);
        
          \node[left]  at (1) {1};
          \node[below right] at (2) {2};
          \node[below] at (3) {3};
          \node[right] at (4) {4};
          \node[above left] at (5) {5};
          \node[above] at (6) {6};
    \end{tikzpicture}
    \caption{An octahedron.}\label{fig:octahedron}
\end{figure}

\begin{example}\label{eg:octahedron}
    Let $K$ be the octahedron, with vertices labeled as in \zcref{fig:octahedron} such that the following pairs do not have an edge between them: $1$ and $4$, $2$ and $5$, $3$ and $6$. Carefully checking the orientations of each tetrahedron in \zcref{fig:octahedron}, one may identify the necessary signs to state the realization space model $\Gr_K^{>0}$ as
    \[
        \Gr_K^{>0} = \left\{ \begin{matrix} p_{3456} > 0, & p_{2456} > 0, & & p_{2346} < 0, & p_{2345} < 0 \\ & p_{1456} > 0, & p_{1356} >0, &  & p_{1345} < 0 \\ & & p_{1256} > 0 , & p_{1246} > 0 , & \\ & & & p_{1236} > 0, & p_{1235} > 0 \\ & & & & p_{1234} > 0\end{matrix} \right\}.
    \]
    Using the duality $\Gr_{4,6} \cong \Gr_{2,6}$ where the Pl\"ucker coordinates on $\Gr_{2,6}$ are $\tilde p_{ij} = (-1)^{i+j} p_{[6] \setminus \{i,j\}}$ up to a common sign, we find the dual of the realization space model as
    \[
        \widetilde \Gr_K^{>0} = \left\{ \begin{matrix} \tilde p_{12} < 0, & \tilde p_{13} > 0, & & \tilde p_{15} < 0, & \tilde p_{16} > 0 \\ & \tilde p_{23} < 0, & \tilde p_{24} >0, &  & \tilde p_{26} < 0 \\ & & \tilde p_{34} < 0 , & \tilde p_{35} > 0 , & \\ & & & \tilde p_{45} < 0, & \tilde p_{46} > 0 \\ & & & & \tilde p_{56} < 0\end{matrix} \right\}.
    \]

    The signs of $\tilde p_{14}$, $\tilde p_{25}$ and $\tilde p_{36}$ decompose $\widetilde \Gr_K^{>0}$ into eight regions of $(\Gr_{2,6},Z_{2,6})$. Summing the canonical forms of these $8$ regions, one finds by linearity of canonical forms, and since the boundary of $\widetilde \Gr_K^{>0}$ lies on the twelve divisors $\{\tilde p_{ij}=0\}$ with $ij\notin\{14,25,36\}$, the expression
    \[
        \begin{aligned}
            \omega_{\widetilde \Gr_K^{>0}} &= \frac{dC}{\tilde p_{14}\tilde p_{42}\tilde p_{25}\tilde p_{53}\tilde p_{36}\tilde p_{61}} + \frac{dC}{\tilde p_{41} \tilde p_{12} \tilde p_{25} \tilde p_{53} \tilde p_{36} \tilde p_{64}} + \cdots \\
            &= \frac{(\tilde p_{16} \tilde p_{24} \tilde p_{35} - \tilde p_{15} \tilde p_{26} \tilde p_{34} )^2}{\tilde p_{12}\tilde p_{13} \tilde p_{23} \tilde p_{34} \tilde p_{35} \tilde p_{45} \tilde p_{15} \tilde p_{16} \tilde p_{56} \tilde p_{24} \tilde p_{26} \tilde p_{46}}dC.
        \end{aligned}
    \]
    One can recover the canonical form on $\Gr_K^{>0}$ by returning to the original Pl\"ucker coordinates.
\end{example}

\subsection{Stacked polytopes}

Stacking a vertex onto a facet is the simplest way to enlarge a simplicial polytope. In the Grassmannian model the new vertex is a new column, and its position relative to the old polytope is recorded by $k$ Pl\"ucker coordinates, one for the facet and one for each of its neighbors. These coordinates come from a nesting, so the recursion of \zcref{sec:nestings} applies, and it shows that stacked polytopes give positive arrangements of combinatorial rank one.

\begin{definition}
    Let $K$ be a simplicial $(k-1)$-polytope with vertex set $[n-1]$ and let $F = \{a_1,\ldots,a_{k-1}\}$ be a facet of $K$. \emph{Stacking} a new vertex $n$ on $F$ produces the simplicial polytope $K'$ with vertex set $[n]$ whose facets are the facets of $K$ other than $F$ together with the $k-1$ sets $(F \setminus a_r) \cup n$. Geometrically, one places the vertex $n$ beyond the facet $F$ and beneath all other facets of a realization of $K$ and takes the convex hull. A simplicial polytope is \emph{stacked} if it is obtained from a simplex by repeated stacking. 
\end{definition}

Throughout this subsection $K$, $F$ and $K'$ are as in the definition. For $r\in[k-1]$ let
\[
G_r=(F\setminus a_r)\cup b_r
\]
be the facet of $K$ adjacent to $F$ along the ridge $F\setminus a_r$. 

We now show how the polytopal divisor changes under stacking in terms of nesting.

\begin{lemma}\label{lem:stacking-nesting}
	The tuple $\Gamma(K,F):=(G_1,\dots,G_{k-1},F)$ is a $k$-nesting on $[n-1]$, with nesting indices $j_r=a_r$ for $r\le k-1$ and $j_k=b_{k-1}$. Its numerator and denominator are
	\[
	N(\Gamma)=\{F\cup b_r:\ r\in[k-1]\},\qquad
	D(\Gamma)=\{G_r\cup n:\ r\in[k-1]\}\cup\{F\cup n\},
	\]
	the elements of $S_K$ that come from $F$, and the elements of $S_{K'}$ that come from the new facets. In particular $S_{K'}\cup N(\Gamma)=S_K\cup D(\Gamma)$, and \zcref{lem:nesting-isomorphism} with $S=S_K$ gives an isomorphism
	\begin{equation}\label{eq:stacking-iso}
		f_\Gamma\colon\ (\Gr_{k,n-1}\setminus Z_K)\times(\C^*)^k\ \longrightarrow\ \Gr_{k,n}\setminus Z_{S_{K'}\cup N(\Gamma)} .
	\end{equation}
\end{lemma}

\begin{proof}
	For $r\le k-1$ we have $a_r\notin G_r$, while $a_r\in G_s$ for $s\ne r$ and $a_r\in F$. So $j_r=a_r$ satisfies \zcref{def:nesting}, and the element of $G_{k-1}$ not in $F$ is $b_{k-1}$. Hence $N(\Gamma)=\{G_r\cup a_r\}=\{F\cup b_r\}$, and $D(\Gamma)$ is as stated. The adjacent pairs of facets of $K'$ are of three kinds: pairs of facets of $K$ not involving $F$, which are adjacent in $K$ as well; the pairs $G_r$ and $(F\setminus a_r)\cup n$ along the ridge $F\setminus a_r$, with union $G_r\cup n$; and the pairs of new facets along $(F\setminus\{a_r,a_s\})\cup n$, with union $F\cup n$. The adjacent pairs of $K$ involving $F$ are the pairs $F,G_r$, with union $F\cup b_r$. Hence $S_{K'}\setminus D(\Gamma)\subseteq S_K$ and $S_K\setminus N(\Gamma)\subseteq S_{K'}$, which gives $S_{K'}\cup N(\Gamma)=S_K\cup D(\Gamma)$. An element $F\cup b_r$ of $N(\Gamma)$ may survive in $S_{K'}$, namely when $b_r=b_s$ for some $s\ne r$, because then $F\cup b_r=G_r\cup G_s$. Finally $N(\Gamma)\subseteq S_K$, so \zcref{lem:nesting-isomorphism} with $S=S_K$ gives \eqref{eq:stacking-iso}.
\end{proof}

\begin{remark}\label{rem:stacking-polygon}
	For the $(n-1)$-gon with vertices $1,\dots,n-1$ in cyclic order and $F=\{n-1,1\}$, stacking gives the $n$-gon, and $\Gamma(K,F)=(\{1,2\},\{n-2,n-1\},\{n-1,1\})$. Its numerator $\{\{n-1,1,2\},\{n-2,n-1,1\}\}$ and its denominator are those of the nesting of \zcref{eg:positive_grassmannian_isomorphism} for $k=3$, of which $\Gamma(K,F)$ is a reordering. So \eqref{eq:stacking-iso} is an isomorphism between the same spaces as the recursion for $\Gr^{\ge0}_{3,n}$ in \zcref{cor:positive-grassmannian-rank}.
\end{remark}

\begin{theorem}\label{thm:stacked}
	Let $K$ be a stacked $(k-1)$-polytope with $n$ vertices. Then $(\Gr_{k,n},Z_K)$ is a positive arrangement, and $\Gr_K^{\geq 0}$ is one of its regions, and
	\[
	\cromb_\R(\Gr_{k,n},Z_K)\ =\ \cromb(\Gr_{k,n},Z_K)\ =\ 1 .
	\]
\end{theorem}

\begin{proof}
	After relabeling the vertices we may assume that $K$ arises from the simplex on $[k]$ by stacking the vertices $k+1,\dots,n$. For $n=k$ the Grassmannian $\Gr_{k,k}$ is a point and there is nothing to show. For $n>k$, write the polytope of the theorem as $K'$, obtained by stacking $n$ on a facet $F$ of a stacked polytope $K$ with vertex set $[n-1]$, and let $\Gamma=\Gamma(K,F)$. We show by induction on $n$ that $(\Gr_{k,n},Z_{K'})$ is a positive arrangement with $\cromb=1$.

	For $n=k+1$, $S_{K'}$ contains the $k$ sets $D(\Gamma)$, which are the $k$-subsets containing $k+1$, and also $[k]=G_r\cup G_s$ for $r\ne s$. So $S_{K'}=\binom{[k+1]}{k}$, and the Pl\"ucker coordinates identify $(\Gr_{k,k+1},Z_{K'})$ with $\PP^k$ and its coordinate hyperplanes. This is a toric pair with $\cromb_\R=\cromb=1$ \cite[\S5.3]{Brown:2025jjg}.

	For $n\ge k+2$, \zcref{lem:stacking-nesting} gives $N(\Gamma)\subseteq S_K$ and $S_{K'}\cup N(\Gamma)=S_K\cup D(\Gamma)$ with $S_K\setminus N(\Gamma)\subseteq S_{K'}$, so $(S_K\setminus N(\Gamma))\cup D(\Gamma)\subseteq S_{K'}\subseteq S_K\cup D(\Gamma)$. By induction, $(\Gr_{k,n-1},Z_K)$ is a positive arrangement with $\cromb=1$. So part (3) of \zcref{thm:recursion} with $S=S_K$ and $S'=S_{K'}$ shows that $(\Gr_{k,n},Z_{K'})$ is a positive arrangement with $\cromb=1$.

	It remains to show that $\Gr_{K'}^{\geq0}$ is a region and that $\cromb_\R(\Gr_{k,n},Z_{K'})=1$, which we also do by induction on $n$. By \zcref{lem:polytope-boundary} the realization space $\Gr_{K'}^{>0}$ is open and closed in $\Gr_{k,n}(\R)\setminus Z_{K'}(\R)$. It is connected because $K'$ is stacked (\zcref{rem:connected-realization}), so it is a connected component, and its closure $\Gr_{K'}^{\geq0}$ is a region. Since $\cromb_\R\le\cromb=1$, it suffices to show that the canonical form $\omega_{K'}$ of $\Gr_{K'}^{\geq0}$ is nonzero.

	For $n=k+1$, the region $\Gr_{K'}^{\geq0}$ is a simplex in $\PP^k(\R)$ bounded by coordinate hyperplanes, and its canonical form is nonzero \cite[\S5.3]{Brown:2025jjg}. For $n\ge k+2$, we have $\omega_K\neq0$ by induction. Let $W$ be the target of \eqref{eq:stacking-iso}. By \zcref{rem:connected-realization}, the first $n-1$ columns of a point of $\Gr_{K'}^{>0}$ lie in $\Gr_K^{>0}$ (for compatible choices of the reference realizations), where $p_{F\cup b_r}\neq0$. So the connected set $\Gr_{K'}^{>0}$ lies in $W(\R)$, and in $f_\Gamma(\Gr_K^{>0}\times O)$ for an orthant $O$ of $(\R^*)^k$. The latter is a connected component of $W(\R)$, and $\Gr_{K'}^{>0}$ is open and closed in $W(\R)$, so the two are equal.

	Canonical forms are compatible with products \cite[Prop.~2.12]{Brown:2025jjg} and with pushforwards \cite[\S2.3.4]{Brown:2025jjg}, here along the isomorphism $f_\Gamma$ and along the identity of $\Gr_{k,n}$, which enlarges the boundary from $Z_{K'}$ to $Z_{S_K\cup D(\Gamma)}$; all these pairs have genus zero by \zcref{thm:recursion}. As the closed orthant $\overline O$ has canonical form $\pm\dlog z_1\wedge\cdots\wedge\dlog z_k$, we get
	\[
	\omega_{K'}|_W=\pm(f_\Gamma)_*\bigl(\omega_K\wedge\dlog z_1\wedge\cdots\wedge\dlog z_k\bigr),
	\]
	the image of $\omega_K\neq0$ under the injective map of part (2) of \zcref{thm:recursion}. Hence $\omega_{K'}\neq0$.

	By induction, for every stacked polytope, in particular for the polytope of the theorem, the closure of the realization space is a region with nonzero canonical form, and $\cromb_\R=\cromb=1$.
\end{proof}

\begin{figure}
	\[
	\begin{tikzpicture}[scale=1.1]
		\coordinate (1) at (0.00,1.57);
		\coordinate (2) at (-1.29,-0.17);
		\coordinate (3) at (-0.17,-1.12);
		\coordinate (4) at (1.46,-0.28);
		\draw[dashed] (2) -- (4);
		\draw (1) -- (2);
		\draw (1) -- (3);
		\draw (1) -- (4);
		\draw (2) -- (3);
		\draw (3) -- (4);
		\node[above] at (1) {1};
		\node[left] at (2) {2};
		\node[below] at (3) {3};
		\node[right] at (4) {4};
	\end{tikzpicture}
	\hspace{2em}
	\begin{tikzpicture}[scale=1.1]
		\coordinate (1) at (0.00,1.57);
		\coordinate (2) at (-1.29,-0.17);
		\coordinate (3) at (-0.17,-1.12);
		\coordinate (4) at (1.46,-0.28);
		\coordinate (5) at (-0.73,0.26);
		\draw[dashed] (2) -- (4);
		\draw (1) -- (2);
		\draw (1) -- (3);
		\draw (1) -- (4);
		\draw (2) -- (3);
		\draw (3) -- (4);
		\draw (1) -- (5);
		\draw (2) -- (5);
		\draw (3) -- (5);
		\node[above] at (1) {1};
		\node[left] at (2) {2};
		\node[below] at (3) {3};
		\node[right] at (4) {4};
		\node[right] at (5) {5};
	\end{tikzpicture}
	\hspace{2em}
	\begin{tikzpicture}[scale=1.1]
		\coordinate (1) at (0.00,1.57);
		\coordinate (2) at (-1.29,-0.17);
		\coordinate (3) at (-0.17,-1.12);
		\coordinate (4) at (1.46,-0.28);
		\coordinate (5) at (-0.73,0.26);
		\coordinate (6) at (0.61,0.21);
		\draw[dashed] (2) -- (4);
		\draw (1) -- (2);
		\draw (1) -- (3);
		\draw (1) -- (4);
		\draw (2) -- (3);
		\draw (3) -- (4);
		\draw (1) -- (5);
		\draw (2) -- (5);
		\draw (3) -- (5);
		\draw (1) -- (6);
		\draw (3) -- (6);
		\draw (4) -- (6);
		\node[above] at (1) {1};
		\node[left] at (2) {2};
		\node[below] at (3) {3};
		\node[right] at (4) {4};
		\node[right] at (5) {5};
		\node[left] at (6) {6};
	\end{tikzpicture}
	\]
	\caption{A tetrahedron (left), with the vertex $5$ stacked on the face $123$ (middle), and with the vertex $6$ then stacked on the face $134$ (right).}\label{fig:stacked-tetrahedron}
\end{figure}
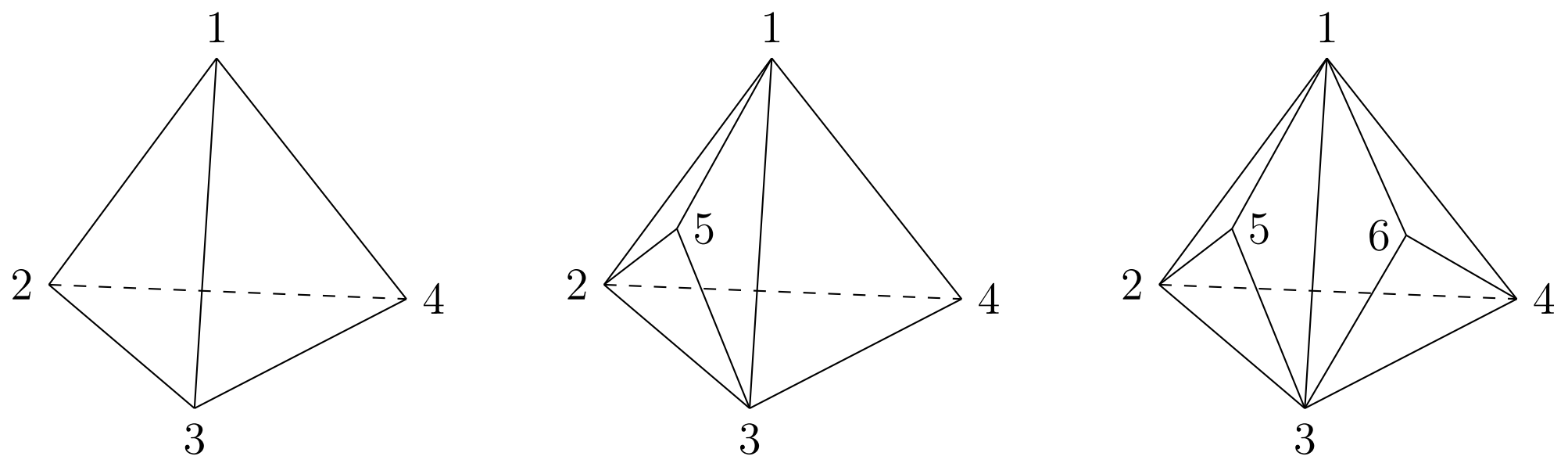

\begin{example}
	Let $K'$ be the stacked $3$-polytope with $6$ vertices obtained from the tetrahedron by stacking the vertex $5$ on the face $123$, obtaining $K$ with $5$ vertices, and then the vertex $6$ on the face $134$, see \zcref{fig:stacked-tetrahedron}.

	The facets of the tetrahedron adjacent to $123$ are $234$, $134$ and $124$, and the facets of $K$ adjacent to $134$ are $234$, $124$ and $135$, so the $4$-nestings of \zcref{lem:stacking-nesting} are
	\[
		\Gamma_1 = (234,134,124,123), \quad \Gamma_2 = (234,124,135,134),
	\]
	with nesting indices $(1,2,3,4)$ and $(1,3,4,5)$, and numerators and denominators
	\[
		\begin{aligned}
			N(\Gamma_1) &= \{1234\}, & D(\Gamma_1) &= \{2345,1345,1245,1235\}, \\
			N(\Gamma_2) &= \{1234,1345\}, & D(\Gamma_2) &= \{2346,1246,1356,1346\}.
		\end{aligned}
	\]
	They give a sequence of isomorphisms
	\[
		\Gr_{4,4} \times (\mathbb C^*)^4 \times (\mathbb C^*)^4 \overset{f_{\Gamma_1} \times \, \mathrm{id}}{\longrightarrow} (\Gr_{4,5} \setminus Z_{K}) \times (\mathbb C^*)^4 \overset{f_{\Gamma_2}}{\longrightarrow} \Gr_{4,6} \setminus (Z_{K'} \cup Z_{\{1345\}}).
	\]
	Using \zcref{lem:form-recursion} and starting with $\omega_{\Gr_{4,4}} = 1 = p_{1234}^{-4} dC$, one obtains the form
	\[
		\begin{aligned}
			\omega_{\Gr_{K'}^{\geq 0}} &= \pm \frac{p_{1234}^2p_{1345}}{p_{2346}p_{1246}p_{1356}p_{1346}} \cdot \frac{p_{1234}^3}{p_{2345}p_{1345}p_{1245}p_{1235}} \cdot p_{1234}^{-4} dC \\
			&= \pm \frac{p_{1234}}{p_{2345}p_{1245}p_{1235}p_{2346}p_{1246}p_{1356}p_{1346}} dC\,.
		\end{aligned}
	\]
\end{example}

\begin{remark}
	Even though \zcref{thm:stacked} refers only to polytopes where the inductive step uses the simplex as a base case, one can use the same argument to stack onto polytopes that are not themselves stacked, such as those from \zcref{lem:small-polytope}.
\end{remark}

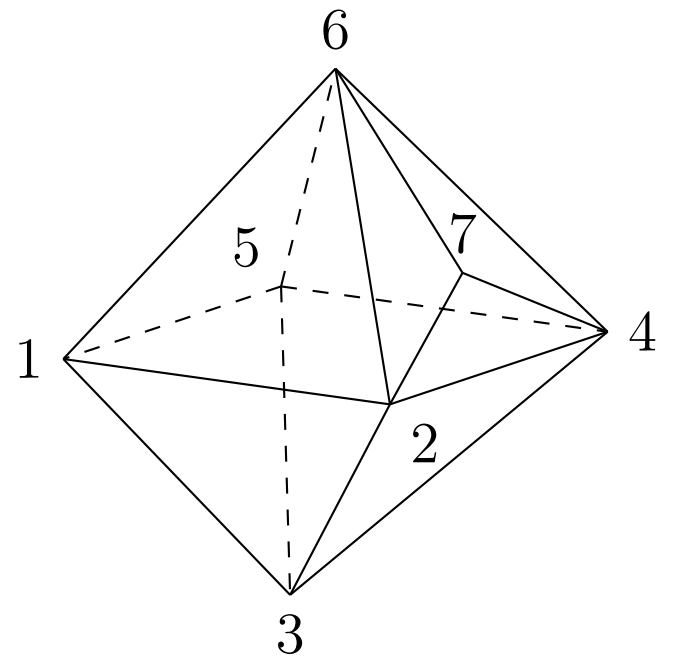
\begin{figure}[ht]
	\[
	\begin{tikzpicture}[scale=0.6]
		\coordinate (1) at (-3,0);
		\coordinate (2) at (0.6,-0.5);
		\coordinate (3) at (-0.5,-2.6);
		\coordinate (4) at (3,0.3);
		\coordinate (5) at (-0.6,0.8);   %
		\coordinate (6) at (0,3.2);
		\coordinate (7) at (1.4,0.95);  %
		\draw[dashed] (5) -- (1);
		\draw[dashed] (5) -- (4);
		\draw[dashed] (5) -- (6);
		\draw[dashed] (5) -- (3);
		\draw (1) -- (6);
		\draw (6) -- (4);
		\draw (6) -- (2);
		\draw (1) -- (2);
		\draw (1) -- (3);
		\draw (3) -- (2);
		\draw (2) -- (4);
		\draw (3) -- (4);
		\draw (7) -- (2);
		\draw (7) -- (4);
		\draw (7) -- (6);
		\node[left]  at (1) {1};
		\node[below right] at (2) {2};
		\node[below] at (3) {3};
		\node[right] at (4) {4};
		\node[above left] at (5) {5};
		\node[above] at (6) {6};
		\node[above] at (7) {7};
	\end{tikzpicture}
	\]
	\caption{The octahedron of \zcref{fig:octahedron} with the vertex $7$ stacked on the face $246$.}\label{fig:stacked-octahedron}
\end{figure}

\begin{example}
	Let $K'$ be the polytope with $7$ vertices obtained from the octahedron $K$ of \zcref{fig:octahedron} by stacking the vertex $7$ on the face $246$, see \zcref{fig:stacked-octahedron}. The facets of $K$ adjacent to $246$ are $456$, $126$ and $234$, so the $4$-nesting of \zcref{lem:stacking-nesting} is
	\[
		\Gamma = (456,126,234,246),
	\]
	with nesting indices $(2,4,6,3)$, numerator $N(\Gamma) = \{2456,1246,2346\}$ and denominator $D(\Gamma) = \{4567,1267,2347,2467\}$. It gives an isomorphism
	\[
		(\Gr_{4,6} \setminus Z_K) \times (\mathbb C^*)^4 \xrightarrow{f_\Gamma} \Gr_{4,7} \setminus (Z_{K'} \cup Z_{N(\Gamma)}).
	\]
	In the Pl\"ucker coordinates of $\Gr_{4,6}$, the canonical form of \zcref{eg:octahedron} reads
	\[
		\omega_{\Gr_K^{\geq 0}} = \pm \frac{(p_{1246}p_{1356}p_{2345} - p_{1256}p_{1345}p_{2346})^2}{p_{1234}p_{1235}p_{1236}p_{1246}p_{1256}p_{1345}p_{1356}p_{1456}p_{2345}p_{2346}p_{2456}p_{3456}} dC .
	\]
	Using \zcref{lem:form-recursion}, the numerator $p_{2456}p_{1246}p_{2346}$ cancels the three poles along $Z_{N(\Gamma)}$, and one obtains
	\[
		\omega_{\Gr_{K'}^{\geq 0}} = \pm \frac{(p_{1246}p_{1356}p_{2345} - p_{1256}p_{1345}p_{2346})^2}{p_{1234}p_{1235}p_{1236}p_{1256}p_{1345}p_{1356}p_{1456}p_{2345}p_{3456} \cdot p_{1267}p_{2347}p_{2467}p_{4567}} dC .
	\]
	The $13$ factors in the denominator are exactly the elements of $S_{K'}$.
\end{example}

\bibliographystyle{amsalpha}
\bibliography{biblio}{}

\end{document}